\documentclass[a4paper,11pt]{article}

\usepackage{hyperref}
\usepackage{amssymb}
\usepackage{amsfonts}
\usepackage{amsmath}
\usepackage{a4wide}

\def\poly{{\tt poly.x}}    \def\class{{\tt class.x}} 	\def\nef{{\tt nef.x}}
\def\glo{{\tt Global.h}}   \def\mori{{\tt mori.x}}      \def\cws{{\tt cws.x}}
\long\def\new#1\endnew{{\bf #1}}
\long\def\del#1\enddel{}

\makeatletter 
 \g@addto@macro\@verbatim\small 
\makeatother

\def\IR{{\mathbb R}}  \def\IP{{\mathbb P}}
  
\def\IZ{{\mathbb Z}} 
 
\def\opt#1{\subsubsection{#1}}

\def\cO{\mathcal{O}}

\def\mP{\mathbb{P}}
\def\mR{\mathbb{R}}
\def\mZ{\mathbb{Z}}
\numberwithin{equation}{section}
\numberwithin{table}{section}
\usepackage{amsxtra}
\DeclareFontFamily{U}{mathx}{\hyphenchar\font45}
\DeclareFontShape{U}{mathx}{m}{n}{
      <5> <6> <7> <8> <9> <10>
      <10.95> <12> <14.4> <17.28> <20.74> <24.88>
      mathx10
      }{}
\DeclareSymbolFont{mathx}{U}{mathx}{m}{n}
\DeclareFontSubstitution{U}{mathx}{m}{n}
\DeclareMathAccent{\widecheck}{0}{mathx}{"71}

\usepackage{alltt}

\begin{document}
\thispagestyle{empty}
\begin{flushright}
TUW-12-10\\
IPMU12-0094\\
\end{flushright}
\vspace{1cm}
\begin{center}
{\LARGE \bf PALP -- a User Manual}
\end{center}
\vspace{8mm}
\begin{center}
Andreas P. Braun${}^1$, Johanna Knapp${}^2$, Emanuel Scheidegger${}^3$,\\ Harald Skarke${}^1$ and Nils-Ole Walliser${}^1$
\end{center}
\vspace{3mm}
\begin{center}
${}^1$ Institute for Theoretical Physics, Vienna University of Technology,\\
Wiedner Hauptstrasse 8-10/136, 1040 Vienna, Austria\\
{\tt abraun, skarke, walliser@hep.itp.tuwien.ac.at}\\\vspace{1mm}
${}^2$ Kavli IPMU (WPI), The University of Tokyo, \\
5-1-5 Kashiwanoha, Kashiwa, Chiba 277-8583, Japan\\
{\tt johanna.knapp@ipmu.jp}\\\vspace{1mm}
${}^3$ Institute for Mathematics, University of Freiburg,\\
Eckerstrasse 1, 79104 Freiburg, Germany \\
{\tt emanuel.scheidegger@math.uni-freiburg.de}
\end{center}
\vspace{15mm}
\begin{abstract}
\noindent This article provides 
a complete user's guide to version 2.1 
of the toric geometry package PALP by Maximilian Kreuzer and others. 
In particular, previously undocumented applications such as the program 
{\tt nef.x} are discussed in detail. 
New features of PALP 2.1 include an extension of the program {\tt mori.x} 
which can now compute Mori cones and intersection rings 
of arbitrary dimension and can also take specific triangulations of  
reflexive polytopes as input. 
Furthermore,
the program {\tt nef.x} is enhanced by an option that allows 
the user 
to enter reflexive Gorenstein cones as input. 
The present documentation is complemented by a Wiki which is available 
online.
\end{abstract}
\vspace{2cm}

\newpage
\setcounter{tocdepth}{2}
\tableofcontents

%%%%%%%%%%%%%%%%%%%%%%%%%%%%%%%%%%%%%%%%%%%%%%%%%%%%%%%%%%%%%%%%%%%%%%%%%%%%%%
\section{Introduction}
\label{sec-intro}
\subsection{A brief history of PALP}
The first lines of code that would eventually become a part of PALP were probably written in 1992.
At that time Max Kreuzer worked together with one of us (HS) on certain quasihomogeneous functions relevant to the description of Landau--Ginzburg models that also had interpretations in terms of Calabi--Yau hypersurfaces in weighted projective spaces.
This culminated in the classification (also found, independently, by Klemm and Schimmrigk \cite{Klemm:1992bx}) of all such functions relevant to standard string compactifications.
As the title of the paper \cite{Kreuzer:1992da}, `No mirror symmetry in Landau--Ginzburg spectra!', suggests, mirror symmetry was incomplete in that class of models and it was necessary to look for more general scenarios.
These were indeed provided by Batyrev's elegant construction of mirror pairs of Calabi--Yau spaces via dual pairs of reflexive polytopes \cite{Batyrev:1994hm}.

After the proposal \cite{Kreuzer:1995cd,Skarke:1996hq} of an algorithm for the classification of reflexive polytopes, work on the implementation of the required routines commenced.
The expertise gained in this project and parts of the code could be used to
consider questions like the manifestation of fibration structures in the toric context \cite{Avram:1996pj,Kreuzer:1997zg} or the connectedness of the moduli space of Calabi--Yau hypersurfaces described by reflexive polytopes \cite{Avram:1997rs}, and these projects in turn enhanced the stock of available C routines.
A first implementation of the whole algorithm led to the generation of the complete list of reflexive 3-polytopes \cite{Kreuzer:1998vb}, but only a thoroughly revised and optimized version of the code could generate all 473,800,776 reflexive polytopes in four dimensions \cite{Kreuzer:2000xy}.

By that time Max was also working with his graduate student Erwin Riegler on an extension to include nef partitions (leading eventually to 
\cite{Kreuzer:2001fu,Klemm:2004km,rieglerPhD} and to {\tt nef.x}).
The collection of available routines had reached a number, a level of complexity and a lack of documentation that would have rendered them useless within a very short time without any efforts at preservation.
Besides, it was clear that the programs might be useful to other people as well.
So it was decided to work on polishing and documenting the existing
routines with the aim of combining them into a publicly available package.
After some time and several candidates (among them `lpoly') the name of the package became PALP, containing {\tt poly.x}, {\tt class.x}, {\tt cws.x}. This is an acronym for `Package for Analysing Lattice Polytopes', but we find it quite appropriate that it shares this name with rather obscure body parts of arachnidae \cite{pedipalp}.

During the period when Max and Erwin were starting to compile nef partitions, one of us (ES) joined Max' group in Vienna as a postdoc. This led to a shift of the focus from the classification of nef partitions more towards applications in mirror symmetry and led to a number of new options in {\tt nef.x}.  

A refinement of polytope data by triangulations and the corresponding Mori cones was made desirable by the following well known facts.
Different triangulations of a polytope,
hence different intersection rings, may lead to topologically distinct
Calabi--Yau manifolds, while non--isomorphic polytopes can give rise to equivalent Calabi--Yau manifolds; the intersection ring is an essential ingredient in Wall's theorem on the
classification of 6-manifolds~\cite{wall1966classification}.
From the point of view of mirror symmetry, the intersection ring and the
Mori cone are important as they enter the GKZ hypergeometric system of
differential equations governing the periods of the mirror hypersurface. 

This was enough motivation to extend the existing routines to the computation of the Mori cone which can be defined entirely in terms of combinatorial data. At that time the triangulation was viewed as an external input determined by some other specialized program such as TOPCOM~\cite{Rambau:TOPCOM-ICMS:2002}. After the initial success one of Max' graduate students (JK) started to develop a code in SINGULAR~\cite{DGPS} that computes from this combinatorial data the intersection rings of the toric variety and the Calabi--Yau hypersurface. This spawned what later would become {\tt mori.x}. 

A couple of years later Max, together with another graduate student (NW) took this up with the goal of creating a routine which determines all the unimodular coherent star triangulations within PALP without having to rely on any external input. 

Despite the fact that PALP was originally designed for the specific purposes 
mentioned above it has become a versatile tool for both mathematics and physics applications. One indicator for the success of PALP is that it has been included into the Sage package \cite{sage} and the Debian repositories.  

\subsection{How to use this manual}
One of the biggest drawbacks of PALP is the combination of complicated syntax and lack of concise documentation. While we decided to keep the syntax and its oddities (cf. section \ref{sec-general}) for the sake of continuity, we would like to overcome the documentation issue with this article and a PALP Wiki which is available at \cite{wiki}. Some parts of PALP have already been discussed previously. 
The original paper accompanying the first version of PALP is \cite{Kreuzer:2002uu}. It contains documentation on the programs {\tt poly.x}, {\tt cws.x} and {\tt class.x}. The program {\tt mori.x} has been presented in \cite{Braun:2011ik}. The program {\tt nef.x} for analyzing complete intersections in toric ambient spaces has been written by Erwin Riegler as part of his PhD thesis \cite{rieglerPhD} but has never been documented. 
In writing the present manual we have tried to cover all applications, i.e. there should be no need to read the older papers as well, except for few passages that we cite at the appropriate points.

In general we do not explain concepts from the theory of polytopes or from toric geometry, except where this serves to fix notation or where we use non--standard terminology. 
The reader is referred to the standard textbooks \cite{Fulton:1993,oda1988convex,Cox:2011ab} or any of a number of reviews (those written by PALP programmers \cite{Skarke:1998yk,Kreuzer:2006ax} are probably closest to the style of this manual).

We recommend that everyone interested in using PALP read section \ref{sec-general} on general aspects of the package, which may hold some surprises even for reasonably experienced users.
The next step is to choose some application of PALP
(consulting the following paragraph should help to decide which program provides this application).
Then one can jump to the section describing that program and read the general part of that section.
Finally one should consult the subsections where the required options are described.

This article is organized as follows. In section \ref{sec-general} we give a general overview of the PALP package and discuss generic properties such as the input of polytope data and error handling. Furthermore we point out some peculiarities of PALP. 
The remaining sections each correspond to one of the executable programs, with a brief general introduction followed by descriptions for all the available options.
Section \ref{sec-poly} is devoted to the program {\tt poly.x} which contains mainly general purpose routines for analyzing lattice polytopes but also some specialized routines related to applications in string theory and algebraic geometry that do not fit into other parts of the package.
In sections \ref{sec-cws} and \ref{sec-class} we describe the programs {\tt cws.x} and {\tt class.x} which have been essential for the classification of reflexive polytopes.
Section \ref{sec-nef} contains the documentation of the program {\tt nef.x} which provides routines to analyze nef partitions of reflexive polytopes. In section \ref{sec-mori} we discuss PALP's most recent application {\tt mori.x} which computes the Mori cone of a toric variety and, with the help of the program SINGULAR, topological data such as intersection numbers of (not necessarily Calabi-Yau) hypersurfaces in those ambient spaces. 
%%%%%%%%%%%%%%%%%%%%%%%%%%%%%%%%%%%%%%%%%%%%%%%%%%%%%%%%%%%%%%%%%%%%%%%%%%%%%%
\section{{General aspects of using PALP}}
\label{sec-general}
In this section we treat aspects of PALP that are common to most or all of its
applications.
The first step is to download the package from the website 
\cite{cypalp}
and follow the compilation instructions given there, which should result in
the existence of a directory `palp' containing the program as well as the 
executable files.

\subsection{Polytope input}
\label{sec:polytope-input}
The majority of applications requires input in the form of a list of polytopes.
There are essentially two ways of entering the data of a polytope.
Matrix input starts with a line containing two numbers $n_{\rm lines}$ and 
$n_{\rm columns}$ (which may be followed by text which is simply ignored
by the program) and proceeds with a matrix with the corresponding numbers
of lines and columns.
PALP requires $n_{\rm lines} \not=n_{\rm columns}$ and interprets the smaller
of the two numbers as the dimension of the polytope and the other one as the 
number of polytope points entered as lines or columns of input.

{\small
\begin{verbatim}
palp$ poly.x
Degrees and weights  `d1 w11 w12 ... d2 w21 w22 ...'
  or `#lines #columns' (= `PolyDim #Points' or `#Points PolyDim'):
3 2  This text is ignored by PALP
Type the 6 coordinates as #pts=3 lines with dim=2 columns:
2 0
0 2
0 0
M:6 3 F:3
Degrees and weights  `d1 w11 w12 ... d2 w21 w22 ...'
  or `#lines #columns' (= `PolyDim #Points' or `#Points PolyDim'):
2 3 The same example with transposed input
Type the 6 coordinates as dim=2 lines with #pts=3 columns:
2 0 0 
0 2 0
M:6 3 F:3
\end{verbatim}}
In both cases the input specifies the polygon (2-polytope) that is the convex 
hull of the 3 points $\{(2,0), (0,2), (0,0)\}$ in $M=\IZ^2$.
The output just means that this polygon has 6 lattice points, 3 vertices and 3
facets (here, edges).
The possibility of ignored text in the input is useful because PALP's output
can often be used as input for further applications; thereby extra information
can be displayed without destroying the permissible format.

For applications in the context of toric geometry one should be aware of the 
fact that there are two relevant, mutually dual lattices $M$ and $N$ whose 
toric interpretations are quite different.
By default PALP interprets the input polytope as $\Delta\subset M_\IR$.
Note that PALP refers to this polytope as $P$; in this paper we shall 
use both notations.
If $\Delta$ ($=P$) is reflexive, it is very natural (and, for some 
applications, more natural) to consider its dual $\Delta^*\subset N_\IR$ as 
well.
If PALP should interpret the input as $\Delta^*$, it must be instructed to
do that by an option ({\tt -D} for \poly\ and \mori, {\tt -N} for \nef).
In fact, in the case of \mori\ it would be extremely unnatural to use 
$\Delta$ as input; therefore matrix input is allowed only with {\tt -D} to
avoid errors.

A second input format uses the fact that many polytopes (in particular the ones
related to the toric description of weighted projective spaces) afford a 
description as the convex hull of all points $X$ that lie in the 
$(n-1)$-dimensional sublattice 
$M\subset\IZ^n$ determined by $\sum_{i=1}^n w_iX_i=0$ and satisfy the 
inequalities $X_i\ge -1$ for $i\in\{1,\ldots,n\}$.
Given such a weight system in the format {\tt d w1 w2 ... wn} where the $w_i$
must be positive integers and
$d=\sum_{i=1}^n w_i$, PALP computes the corresponding list of points and
makes a transformation to $M\simeq \IZ^{n-1}$.
The following example corresponds to the Newton polytope of the quintic 
threefold in $\IP^4$.

{\small
\begin{verbatim}
palp$ poly.x -v
Degrees and weights  `d1 w11 w12 ... d2 w21 w22 ...'
  or `#lines #columns' (= `PolyDim #Points' or `#Points PolyDim'):
5 1 1 1 1 1
4 5  Vertices of P
   -1    4   -1   -1   -1
   -1   -1    4   -1   -1
   -1   -1   -1    4   -1
   -1   -1   -1   -1    4
\end{verbatim}}
As the first line of the prompt indicates, this format can be generalized 
to the case
of $k$ weight systems describing a polytope in $M\simeq \IZ^{n-k}$.
We call the corresponding data, which should satisfy 
$w_{ij}\ge 0$ and $(w_{1j},\ldots, w_{kj})\ne (0,\ldots,0)$, 
a CWS (`combined weight system').

It is also possible to specify a sublattice of finite index corresponding to
the condition $\sum_{i=1}^nl_ix_i=0$ mod $r$ by writing {\tt /Zr: l1 ...ln} 
after the specification of the (C)WS.

In the following example, {\tt 21100 20011} describes a square whose edges 
have lattice length 2, whereas the condition indicated by {\tt Z2: 1 0 1 0}
eliminates the interior points of the edges.
The particular output arises because PALP transforms the original and
the reduced lattice to $\IZ^2$ in different ways.
{\small
\begin{verbatim}
palp$ poly.x -v
Degrees and weights  `d1 w11 w12 ... d2 w21 w22 ...'
  or `#lines #columns' (= `PolyDim #Points' or `#Points PolyDim'):
2 1 1 0 0  2 0 0 1 1
2 4  Vertices of P
   -1    1   -1    1
   -1   -1    1    1
Degrees and weights  `d1 w11 w12 ... d2 w21 w22 ...'
  or `#lines #columns' (= `PolyDim #Points' or `#Points PolyDim'):
2 1 1 0 0  2 0 0 1 1  /Z2: 1 0 1 0
2 4  Vertices of P
   -1    0    0    1
    1   -1    1   -1
\end{verbatim}}
For a reconstruction of the CWS given a polytope as matrix input see
the option {\tt cws.x -N} in Section~\ref{sec:cws-n}.

If PALP is used interactively, it can be terminated by entering an
empty line instead of the data of a polytope. In the case of file input the 
end of the file results in the termination.

\subsection{Error handling}\label{error-handling}
PALP is designed in such a way that it should exit with
an error message rather than crash or display wrong results.
The main sources for problems are inappropriately set
parameters, lack of memory and numerical overflows.
The most important settings of parameters all occur at the beginning
of \glo, which is probably the only file that a user may want to modify.

Here are some typical error messages.
If we want to analyze the Calabi--Yau sixfold that is a
hypersurface in $\IP^7$ with ${\tt poly.x}$, the following will happen
if PALP has been compiled with the default settings. 

{\small
\begin{verbatim}
8  1 1 1 1 1 1 1 1
Please increase POLY_Dmax to at least 7
\end{verbatim}}
In this case one should edit \glo\ (see also section \ref{polydmax}), setting
{\small
\begin{verbatim}
#define          POLY_Dmax       7       /* max dim of polytope    */
\end{verbatim}}
\noindent
and compile again.
Similarly the program may ask for changes of other basic parameters, all of
which are defined within the first 52 lines of \glo.

In many cases we have implemented checks with the help of the
`assert' routine, leading to error messages such as the following.
{\small
\begin{verbatim}
poly.x: Vertex.c:613: int Finish_IP_Check(PolyPointList *, ...
EqList *, CEqList *, INCI *, INCI *): Assertion `_V->nv<32' failed.
Abort
\end{verbatim}}
\noindent
In this case one should look up line 613 of {\tt Vertex.c},
{\small
\begin{verbatim}
    assert(_V->nv<VERT_Nmax);
\end{verbatim}}
\noindent
This indicates that the value of {\tt $_-$V->nv} has risen above the value 32
assigned to {\tt VERT$_-$Nmax} in \glo\ and that the value of
{\tt VERT$_-$Nmax} should be changed correspondingly.
At this point it is important to note that the setting of parameters in \glo\ 
depends on the setting of {\tt POLY$_-$Dmax}:
{\small
\begin{verbatim}
#define         POLY_Dmax        6      /* max dim of polytope     */
...
#if	(POLY_Dmax <= 3)
#define         POINT_Nmax      40      /* max number of points    */
#define         VERT_Nmax       16      /* max number of vertices  */
#define         FACE_Nmax       30      /* max number of faces     */
#define         SYM_Nmax        88      /* cube: 2^D*D! plus extra */

#elif	(POLY_Dmax == 4)
#define         POINT_Nmax     700      /* max number of points    */
#define         VERT_Nmax       64      /* max number of vertices  */
...
\end{verbatim}}
Of course one should then modify a parameter such as {\tt VERT$_-$Nmax} or
{\tt SYM$_-$Nmax} at the position corresponding to the chosen value of 
{\tt POLY$_-$Dmax}.
For {\tt POLY$_-$Dmax} taking values up to 4, the default parameters in \glo\
are chosen in such a way that they work for any reflexive polytope.

While the error messages mentioned above are related to parameter values that 
are too low, excessively high values may also lead to problems such as
slowing down the calculation.
In particular, computation time depends very sensitively upon whether {\tt
VERT$_-$Nmax} is larger than 64.
Very high parameter values may also lead to troubles with memory which 
manifest themselves in error messages such as
{\tt Unable to alloc space for ...} or {\tt Allocation failure in ...} or
even {\tt Segmentation fault}.
In such a case one can only check whether there are possibilities for making
parameters smaller, or use a computer with more RAM. 

An assertion failure that
does not refer to an inequality involving a parameter or an allocation
failure, such as\\[2mm]
{\tt NFX\_Limit in GL -> 1074575416 !!}\\[1mm]
is very likely to point to a numerical overflow.
In such a case it might help to change line 12 of \glo\ from\\[2mm]
{\tt \#define ~~~~~~~~~~Long~~~~~~~~~long}\\[1mm]
to\\[1mm]
{\tt \#define ~~~~~~~~~~Long~~~~~~~~~long long}

These issues are particularly relevant to the analysis of
high-dimensional polytopes, e.g. in the case of {\tt nef.x} with nef
partitions of large length. In this case, it may happen that 
certain parameters
in the header file ${\tt Nef.h}$ may also need to
be modified. Here we give a particularly nasty
example: 
\begin{verbatim}
palp$ nef.x -Lp -N -c6 -P
Degrees and weights  `d1 w11 w12 ... d2 w21 w22 ...'
  or `#lines #columns' (= `PolyDim #Points' or `#Points PolyDim'):
7 9
Please increase POLY_Dmax to at least 12 = 7 + 6 - 1
(nef.x requires POLY_Dmax >= dim N + codim - 1)
\end{verbatim}
This means that in {\tt Global.h} we need to set {\tt POLY\_Dmax} to at least $12$: 
\begin{verbatim}
#define         POLY_Dmax       12       /* max dim of polytope  */
\end{verbatim}
After recompiling PALP we get further but not far enough:
\begin{verbatim}
palp$ nef.x -Lp -N -c6 -P
Degrees and weights  `d1 w11 w12 ... d2 w21 w22 ...'
  or `#lines #columns' (= `PolyDim #Points' or `#Points PolyDim'):
7 9
Type the 63 coordinates as dim=7 lines with #pts=9 columns:
 1  0  0  0  0 -1  0  0 -1
 0  1  0  0  0 -1  0  0 -1
 0  0  1  0  0 -1  0  0 -1
 0  0  0  1  0 -1  0  0  0
 0  0  0  0  1 -1  0  0  0
 0  0  0  0  0  0  1  0 -1
 0  0  0  0  0  0  0  1 -1
M:5214 12 N:10 9  codim=6 #part=1
7 10  Points of Poly in N-Lattice:
    1    0    0    0    0   -1    0    0   -1    0
    0    1    0    0    0   -1    0    0   -1    0
    0    0    1    0    0   -1    0    0   -1    0
    0    0    0    1    0   -1    0    0    0    0
    0    0    0    0    1   -1    0    0    0    0
    0    0    0    0    0    0    1    0   -1    0
    0    0    0    0    0    0    0    1   -1    0
--------------------------------------------------
    1    1    1    1    1    1    0    0    0  d=6  codim=2
    1    1    1    0    0    0    1    1    1  d=6  codim=2
nef.x: Vertex.c:613: Finish_Find_Equations: 
                     Assertion `_V->nv<64' failed.
Aborted
\end{verbatim}
This can be remedied by adjusting the global variable {\tt VERT\_Nmax} in {\tt Global.h} as follows (it should not be too large):
\begin{verbatim}
#define           VERT_Nmax    128    /* !! use optimal value !!  */
\end{verbatim}
After recompilation it works for a while. Then the following error
occurs
\begin{verbatim}
Unable to alloc space for _BL
\end{verbatim}
This means that the program has run out of memory.

\subsection{Some peculiarities of PALP}
Much of PALP's code was written originally with a very specific aim 
(the classification problem) in mind, and not with the intention of 
designing a package that would be immediately accessible to many users.
Other applications were added later by different people.
This has resulted in several peculiar features whose comprehension might
help to avoid errors. 
In the following we list a few of them.

\subsubsection{Option names}
PALP does not have any clear convention on how options are named.
In fact, it can happen that options with the same effects have different 
names in different parts of the package: for instance, in order to
make PALP interpret an input polytope as the (dual) $N$ lattice polytope,
one has to use {\tt -D} with \poly\ and \mori\ but {\tt -N} with \nef.
Ironically, \poly\ also has an option {\tt -N} whose effect is completely 
different.
It is also worth noting that \poly\ and \mori\ admit concatenating
several options into one string, e.g. {\tt poly.x -gve}, {\tt mori.x -Hb}
whereas other programs require them to be separate, e.g. 
{\tt nef.x -Lp -N -c6 -P}.

\subsubsection{Indexing conventions}
Most of PALP, being programmed in C, follows the convention of using 
$\{0,1,\ldots,n-1\}$ as the standard $n$ element set.
So lists of vertices take the form $v_0, v_1,\ldots,v_{n_v-1}$, a list
of points is given as $p_0, p_1,\ldots,p_{n_p-1}$, and so on.
The exception is \mori\ which uses $\{1,\ldots,n\}$ as the standard set.

\subsubsection{Binary representation of incidences}\label{inci}
PALP represents incidences as bit patterns both internally and in its output;
in the case of {\tt mori.x -M} even input is required in that format.
This works in the following manner.
To any face $\phi$ and vertex $v_i$ of $P$ a bit $b_i$
is assigned via $b_i=1$ if $v_i\in\phi$ and $b_i=0$ otherwise;
$v_0,\ldots,v_{n_v-1}$ are ordered as in the 
output of {\tt poly.x -v}. 
This results in a bit pattern $B_v(\phi)$ which
can be written as if it represented a binary number, 
$B_v=b_{n_v-1}b_{n_v-2}\ldots b_0$.
This is the convention implemented in \poly, whereas \mori\ writes it as a 
sequence from left to right, $B_v=b_1\ldots b_{n_v-1}b_{n_v}$ (remember the
last subsection about indexing conventions!).
Furthermore, a bit pattern $B_f=\{\tilde b_j\}$ related to $P$'s facets 
$f_j$ (ordered as in the output of {\tt poly.x -e}) can be assigned via 
$\tilde b_j=1$ if $\phi\subseteq f_j$ and $\tilde b_j=0$ otherwise.

\subsubsection{The parameter {\tt POLY\_Dmax}}\label{polydmax}
{\tt POLY\_Dmax} (line 22 of \glo) is the most important parameter to set 
before compilation;
for most applications it is probably the only parameter one has to care about.
While it usually suffices to have {\tt POLY\_Dmax} not smaller than the 
dimension of any polytope that one wants to analyze, there are the following
important exceptions.\\
\nef\ normally requires {\tt POLY\_Dmax} 
$\ge\mathrm{dim}(N) + \mathrm{codim} - 1$,
which is the dimension of the support polytope of the 
corresponding Gorenstein cone which \nef\ analyzes;\\
{\tt nef.x -G} requires {\tt POLY\_Dmax} $\ge$ dim(input-polytope) $+1$
because the input-polytope is interpreted as the support of the cone and the
full dimension of the cone is required for technical reasons;\\
{\tt mori.x -M} requires {\tt POLY\_Dmax} $\ge n_p -\mathrm{dim}(N) - 1$,
the dimension of the Mori cone; $n_p$ is the total number of 
input points including the lattice origin, hence the subtraction of 1.

\subsubsection{IP property and IP simplices}
\label{ipps}
In PALP's help screens and output messages every now and then the 
abbreviation IP occurs.
A priori this is just a shortcut for writing `interior point', as in
`the generic CY hypersurface does not intersect the divisors corresponding
to IPs of facets'.
However, when we say that a polytope has the IP property, we mean that the 
polytope has the lattice origin in its interior.
IP simplices are simplices with this property, usually with vertices 
that are points or vertices of some given polytope.
This concept played an important role in the classification scheme of
\cite{Kreuzer:1995cd,Skarke:1996hq,Kreuzer:2000qv}.
The lattice vectors corresponding to the vertices of an IP simplex define a 
linear relation with positive coefficients which is unique up to scaling.
Conversely any positive linear relation among lattice points that cannot be 
written as the sum of two other such relations defines an IP simplex.
Of course, the set of coefficients 
is just the weight system defined by the IP simplex.
%%%%%%%%%%%%%%%%%%%%%%%%%%%%%%%%%%%%%%%%%%%%%%%%%%%%%%%%%%%%%%%%%%%%%%%%%%%%%%
\section{{\tt poly.x}}
\label{sec-poly}
\subsection{General description of {\tt poly.x}}
{\tt poly.x} is the program that provides an interface for PALP's general 
purpose routines as well as certain specialized applications that do not
belong to any of the other executables.
In other words, \poly\ deals with all applications that are not related to
nef partitions, Mori cones or the classification of reflexive polytopes.
As for all of PALP's programs, a rough guide can be 
obtained with the help option:

\begin{verbatim}
palp$ poly.x -h

This is 'poly.x':  computing data of a polytope P
Usage:   poly.x [-<Option-string>] [in-file [out-file]]

Options (concatenate any number of them into <Option-string>):
h  print this information           
f  use as filter                    
g  general output:                  
   P reflexive: numbers of (dual) points/vertices, Hodge numbers 
   P not reflexive: numbers of points, vertices, equations    
p  points of P                      
v  vertices of P                    
e  equations of P/vertices of P-dual
m  pairing matrix between vertices and equations                  
d  points of P-dual (only if P reflexive)          
a  all of the above except h,f      
l  LG-`Hodge numbers' from single weight input                   
r  ignore non-reflexive input       
D  dual polytope as input (ref only)
n  do not complete polytope or calculate Hodge numbers        
i  incidence information            
s  check for span property (only if P from CWS)           
I  check for IP property            
S  number of symmetries             
T  upper triangular form	       
N  normal form                      
t  traced normal form computation   
V  IP simplices among vertices of P*
P  IP simplices among points of P* (with 1<=codim<=# when # is set)
Z  lattice quotients for IP simplices
#  #=1,2,3  fibers spanned by IP simplices with codim<=#        
## ##=11,22,33,(12,23): all (fibered) fibers with specified codim(s) 
   when combined: ### = (##)#       
A  affine normal form
B  Barycenter and lattice volume [# ... points at deg #]
F  print all facets
G  Gorenstein: divisible by I>1
L  like 'l' with Hodge data for twisted sectors
U  simplicial facets in N-lattice
U1 Fano (simplicial and unimodular facets in N-lattice)
U5 5d fano from reflexive 4d projections (M lattice)
C1 conifold CY (unimodular or square 2-faces)
C2 conifold FANO (divisible by 2 & basic 2 faces)
E  symmetries related to Einstein-Kaehler Metrics

Input:    degrees and weights `d1 w11 w12 ... d2 w21 w22 ...'
          or `d np' or `np d' (d=Dimension, np=#[points]) and
              (after newline) np*d coordinates
Output:   as specified by options
\end{verbatim}
If an input file is indicated, \poly\ reads its input data from there;
otherwise it asks for input interactively.
The output is displayed to the screen unless an output file is specified.
The following subsection will explain all of the possible options, in the order 
in which they appear in the help screen.
Here is a rough guide in terms of specific topics:
\begin{itemize}
\item General polytope analysis without reference to toric geometry or 
string theory: 
{\tt -g, -p, -v, -e, -m, -d, -a, -i, -s, -I, -T, -N, -t, -A, -B, -F, -G, -U, -U1, -E};
\item Conifold singularities: {\tt -C1, -C2};
\item Fano varieties: {\tt -U1, -U5, -C2};
\item Fibration structures: {\tt -[number], -V, -P, -Z};
\item IP property (see section \ref{ipps}): {\tt -I};
\item IP simplices (see section \ref{ipps}): {\tt -V, -P, -Z, -[number]};
\item Landau-Ginzburg type superconformal field theories: {\tt -l, -L};
\item Normal forms: {\tt -N, -A};
\item Sublattices and quotient actions: {\tt -S, -Z, -G};
\item Symmetries of a polytope: {\tt -S, -t}.
\end{itemize}
\subsection{Options of {\tt poly.x}}
As many of the following options are very simple to use, we do not always
provide examples.
In such cases we highly recommend to try using the option with 
simple input, e.g. the weight system {\tt 5 1 1 1 1 1} corresponding to the 
quintic, and a non-reflexive example such as

\begin{verbatim}
3 2
Type the 6 coordinates as #pts=3 lines with dim=2 columns:
2 0
0 2
0 0
\end{verbatim}

\opt{no option set}
In this case the program behaves as if the {\tt -g} option (see below) were set.
This is also the case if no option  other than {\tt -r}, {\tt -D}, {\tt -n}, {\tt -Z}, {\tt -U}
or {\tt -U1}, which do not generate any output per se, is applied.
\opt{-h}
The help screen is displayed (see above).
\opt{-f}
The filter flag switches off the prompt for input data.
This is useful for building pipelines.
\opt{-g}
The following output is generated.
First, the input is repeated if it is of weight/CWS type, 
but not otherwise.
Then the numbers $\# p$ and $\# v$ of lattice points and vertices,
respectively, are displayed in the format `M: $\# p$ $\# v$'.
The remaining output depends on whether $P$ is reflexive.
In this case the numbers $\# d$ and $\# e$ of dual lattice points
and vertices are displayed in the format `N: $\# d$ $\# e$',
followed by information on the Hodge numbers of the corresponding
Calabi--Yau manifold if dim$(P)\ge 4$;
in the case of a three dimensional polytope corresponding to a K3
surface, where the Hodge numbers are determined anyway, information on
the Picard number and the `correction term' is given instead
(the latter is the non-linear term in 
Batyrev's formula for the Picard number \cite{Batyrev:1994hm};
the Picard numbers of a K3 and its mirror add up to $20 + Cor$).
For non-reflexive $P$ the number $\# e$ of facets is shown as `F: $\# e$'.

Using this option implies the completion of the set of lattice
points in the convex hull (`points' in the help screen always means
`lattice points').
In the reflexive case it also leads to the completion of the dual
polytope and the computation of the complete incidence structure which
is required for the calculation of the Hodge numbers.
For large dimensions
these tasks may result in a long response time or in a crash of the program.
In such a case
one should use other options, e.g. {\tt -nve},
if information on the number of lattice points or Hodge numbers is not
required.
\opt{-p}
The lattice points of the polytope are displayed.
\opt{-v}
The vertices of the polytope are displayed.
\opt{-e}
The equations of the hyperplanes bounding the polytope are displayed.
If the polytope is not reflexive, these facet equations are given as lines
`$a_1~\ldots~a_n~~c$' normalized
such that the $a_i$ have no common divisor and the
inequalities $\vec a\cdot\vec x + c \ge 0$ are satisfied for all points of $P$.

Reflexivity of a lattice polytope $P$ is equivalent to $P^*$ being a
lattice polytope, i.e. to $c=1$ for all facet equations. In that situation
the lines $\vec a$ can be interpreted as vertices of $P^*$ and
\poly\ omits the final column of 1's, indicating that the resulting matrix 
can be interpreted as the list of vertices of the dual polytope.
This has the advantage that the output can be used as input for further 
computations. 

\begin{verbatim}
palp$ poly.x -e
Degrees and weights  `d1 w11 w12 ... d2 w21 w22 ...'
  or `#lines #columns' (= `PolyDim #Points' or `#Points PolyDim'):
3 2
Type the 6 coordinates as #pts=3 lines with dim=2 columns:
1 0
0 1
0 0
3 2  Equations of P
   1   0     0
   0   1     0
  -1  -1     1
Degrees and weights  `d1 w11 w12 ... d2 w21 w22 ...'
  or `#lines #columns' (= `PolyDim #Points' or `#Points PolyDim'):
3 2
Type the 6 coordinates as #pts=3 lines with dim=2 columns:
1 0
0 1
-1 -1
3 2  Vertices of P-dual <-> Equations of P
   2  -1
  -1   2
  -1  -1
\end{verbatim}
In the first case the output indicates that $P$ can be described by
$x_1\ge 0$, $x_2\ge 0$, $-x_1-x_2+1\ge 0$; in the second case, where a last 
output column of 1's only is implicit, $P$ corresponds to
$2x_1-x_2+1\ge 0$, $-x_1+2x_2+1\ge 0$, $-x_1-x_2+1\ge 0$.

\opt{-m}\label{vpm}
One gets the $n_v\times n_e$ matrix with entries
$\vec a_j\cdot\vec v_i + c_j$, $1\le i\le n_v$, $1\le j\le n_e$, where
$n_v$, $n_e$ are the numbers of vertices and equations, respectively.
The elements of this `pairing matrix'
represent the lattice distances between the respective vertices and
facets.
The orders of vertices and facets are the same as for {\tt -v} and {\tt -e},
so it is useful to combine {\tt -m} with these options to see precisely
which vertex and facet an entry of the pairing matrix corresponds to.

\opt{-d} 
If the polytope is reflexive the lattice points of the
dual polytope are displayed.

\opt{-a}
This is a shortcut for {\tt -gpvemd}; it can be combined with any other
options.

\opt{-l}
This option is relevant to applications in the context of Landau-Ginzburg
models. 
Together with its close relative {\tt -L}, it is the only option of \poly\ that
requires non-standard input.
Rather than indicate a polytope via matrix or CWS input, one specifies 
a single weight system which need not satisfy $\sum n_i = d$, which
is interpreted as data for a superconformal field theory.
If the central charge of this SCFT is a multiple of 3 (which is required by
PALP 2.1), the analogue of the Hodge numbers 
\cite{Vafa:1989xc,Intriligator:1990ua} is computed.

\begin{verbatim}
palp$ poly.x -l
type degree and weights  [d  w1 w2 ...]: 5 1 1 1 1 1
5 1 1 1 1 1 M:126 5 N:6 5 V:1,101 [-200]
type degree and weights  [d  w1 w2 ...]: 3 1 1 1 1 1 1
3 1 1 1 1 1 1 M:56 6 F:6 LG: H0:1,0,1 H1:0,20 H2:1 RefI2
type degree and weights  [d  w1 w2 ...]: 3 1 1 1 1 1 1 /Z3: 
0 1 2 0 1 2 3 1 1 1 1 1 1 /Z3: 0 1 2 0 1 2 
M:20 6 F:6 LG: H0:1,0,1 H1:0,20 H2:1 RefI2
\end{verbatim}
Here the first example is the familiar quintic treated as the Gepner model 
$(3)^5$, the second example is the Gepner model $(1)^6$ and the third example 
an orbifold of the second one.
More information can be found in section 4.1 of \cite{Kreuzer:2002uu}.

\opt{-r}
Any input that does not correspond to a reflexive polytope will be ignored.
This is useful for filtering out reflexive polytopes from a larger list and
saves calculation time if one is interested only in reflexive polytopes.

\opt{-D}
The input is regarded as the dual polytope $P^*\subset N_\IR$.
As this makes sense only in the reflexive case there is an error
message (but no exit from the program) for non-reflexive input.
This option is useful, in particular, if one wants to have
control over the order of the points in the $N$ lattice.

\opt{-n}
The completion of the set of lattice points is suppressed.
Hence the Hodge numbers cannot be calculated and the output will look
like the one for non-reflexive polytopes even in the reflexive case.
In particular, if the input is not of the (C)WS type, the number of
points may be displayed wrongly.
If dim($M$) is large this option saves a lot of calculation time.

\opt{-i}
Information on the incidence structure is displayed in the following manner.
Remember from section \ref{inci} that any face $\phi$ can be assigned a bit
pattern $B_v(\phi)$ encoding which vertices lie on $\phi$ and a bit pattern
$B_f(\phi)$ encoding to which facets $\phi$ belongs.
{\tt poly.x -i} displays both types of bit patterns as binary 
numbers for all faces of $P$, 
with `{\tt v[i]:}' starting a line of $B_v$'s corresponding to $i$-faces 
and `{\tt f[j]:}' starting a line of $B_f$'s corresponding to $j$-faces.
Bit patterns at the same positions in the `{\tt v[i]:}' and `{\tt f[i]:}' 
lines correspond to the same $i$-faces.
The orders of faces within the lines are a consequence of the way 
PALP computes them; they conform to the output of other options in the
case of $i=n-1$ (facets) but not for $i=0$ (vertices).

\opt{-s}
This option refers to a property of (combined) weight systems that was
considered in the early stages of the classification program 
\cite{Kreuzer:1995cd}.
In the higher dimensional embedding defined by a
weight system, the polytope is bounded by the inequalities $X_i\ge -1$.
We say that the polytope has the span property if the pullbacks of
the equations $X_i=-1$ to the subspace carrying the polytope are spanned by
vertices of the polytope, i.e. if these equations correspond to
facets.
With {\tt -s} a message is given if the (combined) weight systems does
not have this property (try, for example, the weight system `8 3 3 2').

\opt{-I}
There is a message if the polytope does not have
the origin of the coordinate system in its interior.

\opt{-S}
The output contains the following two numbers.
The first is the number of lattice automorphisms (elements of
$GL(n,\IZ)$) that leave the polytope invariant; 
each such automorphism acts as a permutation on the set of vertices.
The second one is the number of permutations of the set of vertices
that leave the vertex pairing matrix (see section \ref{vpm}) invariant 
(after taking into account the induced permutations of facets).
This number can also be interpreted as the number of symmetries of
the polytope in $\IR^n$; it may be larger than the
number of symmetries in the given lattice.

\begin{verbatim}
palp$ poly.x -S
Degrees and weights  `d1 w11 w12 ... d2 w21 w22 ...'
  or `#lines #columns' (= `PolyDim #Points' or `#Points PolyDim'):
5 1 1 1 1 1
#GL(Z,4)-Symmetries=120, #VPM-Symmetries=120
Degrees and weights  `d1 w11 w12 ... d2 w21 w22 ...'
  or `#lines #columns' (= `PolyDim #Points' or `#Points PolyDim'):
5 1 1 1 1 1 /Z5: 0 1 2 3 4
#GL(Z,4)-Symmetries=20, #VPM-Symmetries=120 
\end{verbatim}
In the first case the symmetry group of the lattice polytope and that of the 
vertex pairing matrix are just the permutation group of the 5 vertices, 
of order 120.
In the second case the invariance under the $\IZ_5$ group fixes any 
permutation once the permutations of two of the five vertices have been chosen,
reducing the number of group elements to 20.
The vertex pairing matrix remains the same, namely diag(5, 5, 5, 5, 5), and 
therefore keeps all 120 symmetries.

\opt{-T}
A coordinate change is performed that makes the matrix of
coordinates of the points specified in the input upper triangular,
with minimal entries above the diagonal.
This may be useful for representing the polytope in a specific lattice
basis consisting of points of $P$, or for finding the volume of a
specific cone (if the generators of the cone are the first input
points, the volume will be the product of the entries in the diagonal
after the transformation).
Using this option only makes sense if its results are displayed.
Therefore it exits with an error if it is not combined with an output
generating option ({\tt -v} is a natural choice). 

\opt{-N}\label{normalform}
This option leads to the computation of a normal form of the polytope,
i.e. a matrix containing the vertices in a specific order in a
particular coordinate system, such that this output is the same for
any two polytopes related by a lattice automorphism (see section 3.4 of 
\cite{Kreuzer:1998vb} for the actual algorithm).
This is useful, for example, if two polytopes are suspected to be
isomorphic because they are isomorphic if and only if their normal
forms are identical.

Example: the weight systems `3402 40 41 486 1134 1701' and 
`3486 41 42 498 1162 1743' give rise to the same pair of Hodge numbers 
(491, 11). 
The suspicion that they correspond to the same polytope is confirmed by

\begin{verbatim}
palp$ poly.x -N
Degrees and weights  `d1 w11 w12 ... d2 w21 w22 ...'
  or `#lines #columns' (= `PolyDim #Points' or `#Points PolyDim'):
3402 40 41 486 1134 1701
4 5  Normal form of vertices of P    perm=43210
   1   0   0   0 -42
   0   1   0   0 -28
   0   0   1   0 -12
   0   0   0   1  -1
Degrees and weights  `d1 w11 w12 ... d2 w21 w22 ...'
  or `#lines #columns' (= `PolyDim #Points' or `#Points PolyDim'):
3486 41 42 498 1162 1743
4 5  Normal form of vertices of P    perm=43210
   1   0   0   0 -42
   0   1   0   0 -28
   0   0   1   0 -12
   0   0   0   1  -1
\end{verbatim}
The {\tt perm=43210} part indicates how the vertices of the polytope had to be 
permuted to arrive at the normal form; this is only interesting if the input is
of matrix type.
\opt{-t}
The calculation of the normal form involves determining the
pairing matrix, a normal form for the pairing matrix, an analysis of
which symmetries leave this normal form invariant, a preferred
ordering of the vertices and a conversion of the resulting vertex
coordinate matrix to upper triangular form.
The results of these steps, which may provide further useful
information on the structure of $P$, are displayed.

\opt{-V}
The IP simplices (see section \ref{ipps}) whose vertices are also vertices 
of the dual polytope are displayed. For illustrating examples look at the 
next options ({\tt -P}, {\tt -Z}) which are closely related.

\opt{-P}
The IP simplices (see section \ref{ipps}) whose vertices are lattice points 
of the dual polytope $P^*$ are displayed.
This option should only be used if it is fairly clear that $P^*$
does not have too many lattice points.

\begin{verbatim}
palp$ poly.x -P
Degrees and weights  `d1 w11 w12 ... d2 w21 w22 ...'
  or `#lines #columns' (= `PolyDim #Points' or `#Points PolyDim'):
6 1 2 3
2 7  points of P-dual and IP-simplices
    1    0   -2   -1    0   -1    0
    0    1   -3   -2   -1   -1    0
------------------------------    #IP-simp=4
    2    3    1    0    0    0   6=d  codim=0
    1    2    0    1    0    0   4=d  codim=0
    1    1    0    0    0    1   3=d  codim=0
    0    1    0    0    1    0   2=d  codim=1
\end{verbatim}
This example shows that the $N$ lattice polytope for $\IP^2_{(1,2,3)}$ contains 
3 lattice triangles with the origin in their respective interiors, as well
as a lattice line segment with that property.
The weight systems corresponding to these simplices are indicated.

\opt{-Z}
This option only has an effect if combined with {\tt -V} or {\tt -P}. 
Any IP simplex $S$ of dimension $d$ occurring there defines two potentially 
distinct $d$ dimensional sublattices of $N$: 
the lattice generated by the vertices of $S$,
and the sublattice of $N$ lying in the linear subspace that is spanned by $S$.
Then the program computes the corresponding quotient action.
The following example illustrates this for the case of a bipyramid over a
triangle whose vertices generate an index 3 sublattice of $\IZ^2$.

\begin{verbatim}
palp$ poly.x -VZD
Degrees and weights  `d1 w11 w12 ... d2 w21 w22 ...'
  or `#lines #columns' (= `PolyDim #Points' or `#Points PolyDim'):
3 5
Type the 15 coordinates as dim=3 lines with #pts=5 columns:
   -1   -1    2    0    0
   -1    2   -1    0    0
    0    0    0    1   -1
3 5  vertices of P-dual and IP-simplices
   -1   -1    2    0    0
   -1    2   -1    0    0
    0    0    0    1   -1
-------------------------   #IP-simp=2 I=3 /Z3: 2 1 0 0 0
    1    1    1    0    0    3=d  codim=1 /Z3: 2 1 0 0 0
    0    0    0    1    1    2=d  codim=2
\end{verbatim}

\opt{-[numbers]}
If one of the options {\tt -B}, {\tt -U} or {\tt -C} is set, any number is assumed to
refer to that option (see below for descriptions).
Otherwise information on fibration structures of the Calabi--Yau manifold
corresponding to a reflexive polytope is displayed.
These structures correspond to reflexive subpolytopes of the $N$
lattice polytope $P^*$ that are intersections of the dual polytope with a
linear subspace of $N_\IR$; see, e.g., 
\cite{Avram:1996pj,Kreuzer:1997zg,Kreuzer:2000qv}.
As this is a time consuming task and the desired output format may
vary, there exist two different versions.

If a single number $\#\in\{1,2,3\}$ is specified, the intersections of
$P^*$ with all linear subspaces spanned by IP simplices (see section 
\ref{ipps}) are checked for reflexivity.
This is much faster than a complete search for fibration structure but
misses fibrations whose corresponding subspaces are only spanned
by a combination of two or more IP simplices.
The codimension of the IP simplices is restricted to
$1\le$codim$\le\#$ (in contrast to our usual policy this option has a
side effect on {\tt -P} if combined with that option). The output has
the same structure as the output of ${\tt nef.x -F*}$. For further
details see Section~\ref{sec-nef-F*}. 

If two numbers are specified, all fibration structures of the given
type are computed.
For {\tt -11}, {\tt -22} and {\tt -33} all reflexive sections of codimension 1, 2 and 3,
respectively, are constructed.
For {\tt -12} and {\tt -23} all reflexive subpolytopes of codimension 1 and 2 that
themselves contain a reflexive subpolytope with relative codimension 1 are
constructed. The output is a polytope in the $N$ lattice whose choice of
bases and whose order of points reflects the fibration structure.
For example, CY threefolds that are both K3 and elliptically fibered are
found by applying {\tt poly.x -12} to reflexive 4-polytopes, and fourfolds of 
that type are found by applying {\tt poly.x -23} to reflexive 5-polytopes.

\begin{verbatim}
palp$ poly.x -12
Degrees and weights  `d1 w11 w12 ... d2 w21 w22 ...'
  or `#lines #columns' (= `PolyDim #Points' or `#Points PolyDim'):
84 1 1 12 28 42
4 25  Em:7 3 n:7 3  Km:24 4 n:24 4  M:680 5 N:26 5  p=13bgjn256789...
 1  0 -2 -1  0 -1  0 -14 -12 -10 -8 -6 -4 -9 -7 -5 -3 -4 -2 -7 -5 ...
 0  1 -3 -2 -1 -1  0 -21 -18 -15 -12 -9 -6 -14 -11 -8 -5 -7 -4 -10...
 0  0  0  0  0  0  1 -6 -5 -4 -3 -2 -1 -4 -3 -2 -1 -2 -1 -3 -2 -1 ...
 0  0  0  0  0  0  0  0  0  0  0  0  0  0  0  0  0  0  0  0  0  0 ...
\end{verbatim}
Here {\tt Em:7 3 n:7 3} refers to data of the elliptic fiber $E$ lying within 
the K3 fiber indicated by {\tt Km:24 4 n:24 4}. 
The 25 non-zero points of the $N$ lattice polytope are displayed in 
an order (first points from the subspace corresponding to $E$,
then points in the K3 but not in $E$, finally all other points)
and a coordinate system (only the first 2 coordinates non-vanishing for $E$,
last coordinate vanishing for the K3) that reflect the corresponding nesting 
of polytopes.
The sequence after {\tt p=} indicates what permutation with respect to the original 
order of these points led to this representation (this is only interesting
if the points were entered directly with the {\tt -D} option).
Further examples on how to use the fibration options can be found in section
4.2 of \cite{Kreuzer:2002uu}.

\opt{-A}\label{anf}
Given an arbitrary lattice polytope $P$, {\tt poly.x -A}
computes its `affine normal form', i.e. a polytope in ${\mathbb Z}^n$ 
affinely isomorphic to $P$, such that the normal forms of any two
polytopes $P$ and $Q$
coincide if and only if $P$ and $Q$ are related by an affine lattice 
isomorphism
(cf. {\tt -N}, which performs the same task w.r.t. linear rather than affine 
transformations).

\opt{-B}
For a given polytope its volume (normalized such that the standard simplex 
has volume 1) and the coordinates of its barycentre are displayed. 

If an integer $n$ is specified after {\tt -B}, the polytope is interpreted as the 
origin and the first level of a Gorenstein cone (for a discussion see
Section~\ref{sec:nef-part}. The points of the cone up to  
level $n$ are computed and displayed together with information on the type 
of face of the cone they represent (w.r.t. codimension, i.e. the origin has 
maximal codimension and points interior to the cone have {\tt
  cd=0}). As example, we consider the Gorenstein cone with support
polytope given by the dual of the Newton polytope of the quintic
hypersurface in $\mP^4$.
\begin{verbatim}
palp$ poly.x -B2
Degrees and weights  `d1 w11 w12 ... d2 w21 w22 ...'
  or `#lines #columns' (= `PolyDim #Points' or `#Points PolyDim'):
5 6
Type the 30 coordinates as dim=5 lines with #pts=6 columns:
1  1  1  1  1  0
1  0  0  0 -1  0
0  1  0  0 -1  0
0  0  1  0 -1  0
0  0  0  1 -1  0
vol=5, baricent=(5,0,0,0,0)/6
IPs:
 2 -2 -2 -2 -2  cd=4
 2 -1 -1 -1 0  cd=3
 2 0 0 0 2  cd=4
 2 -1 0 -1 -1  cd=3
 2 0 1 0 1  cd=3
 2 0 2 0 0  cd=4
 2 -1 -1 0 -1  cd=3
 2 0 0 1 1  cd=3
 2 0 1 1 0  cd=3
 2 0 0 2 0  cd=4
 2 -1 -1 -1 -1  cd=0
 2 0 0 0 1  cd=0
 2 0 1 0 0  cd=0
 2 0 0 1 0  cd=0
 2 0 0 0 0  cd=0
 2 0 -1 -1 -1  cd=3
 2 1 0 0 1  cd=3
 2 1 1 0 0  cd=3
 2 1 0 1 0  cd=3
 2 1 0 0 0  cd=0
 2 2 0 0 0  cd=4
 1 -1 -1 -1 -1  cd=4
 1 0 0 0 1  cd=4
 1 0 1 0 0  cd=4
 1 0 0 1 0  cd=4
 1 0 0 0 0  cd=0
 1 1 0 0 0  cd=4
 0 0 0 0 0  cd=5
\end{verbatim}

\opt{-F}
Each facet of a polytope is displayed by listing its vertices in some basis 
for the sublattice carrying the facet. The result is not the affine normal 
form of the facet, however.

\opt{-G}
If the input polytope $P$ is a non-trivial multiple $P=gQ$ of another 
lattice polytope $Q$, the maximal proportionality factor $g$ and $P$ are
 displayed.

\opt{-L}
This option results in the same calculations as the {\tt -l} option (see above), 
but gives a more detailed output including the Witten index.

\opt{-U}
If\footnote{We are very grateful to Benjamin Nill for providing infomation for the {\tt -U}, {\tt -C} and {\tt -E} options.} {\tt -U} is specified without a number following (otherwise see below), it is 
computed whether the $N$ lattice polytope has only simplicial facets; if this 
is not the case no further computations are performed on the polytope and no 
output results from it.

\opt{-U1}
Like {\tt -U} without number, but the facets now have to be unimodular (i.e. of 
volume 1). This corresponds to the case of Fano varieties.

\opt{-U5}
This option is related to the classification of Fano polytopes (unimodular
and simplicial reflexive polytopes) up to dimension five \cite{Kreuzer:2007zy}. 
The corresponding data supplement 
\cite{math0702890}
contains examples of the usage of {\tt -U5}.

\opt{-C1}
This option was implemented for a search of Calabi-Yau threefolds which are 
related to previously known ones via conifold transitions
\cite{Batyrev:2008rp}. 
The use of PALP to produce these results is described at
\cite{math08023376}.
The input should be a 4-dimensional reflexive polytope $P$ such that $P^*$ 
has only basic triangles or squares as two-dimensional faces. In this case
the associated generic CY-hypersurface has isolated conifold singularities. 
The option {\tt poly.x -C1} checks whether this three-dimensional CY is smoothable.

\opt{-C2}
The option {\tt -C2} was used to generate three-dimensional Fano
hypersurfaces with conifold singularities~\cite{Batyrev:2012ab}. For
the corresponding data see~\cite{mathFano}.
The input must be a reflexive 4-polytope that is divisible by 2, 
in which case the generic hypersurface associated to $P/2$ is a 
three-dimensional Fano variety.
A list of (hopefully all) such polytopes can be found at Max Kreuzer's site
\cite{Reflexive2x}.
It is not clear how he obtained this list. 
Possibly he piped an extraction of the complete database of reflexive 
4-polytopes through {\tt poly.x -G} and then used some script to eliminate 
polytopes that are odd multiples of other polytopes.
The program requires the input to be in matrix (rather 
than weight) format.

Example: {\tt poly.x -C2} with input `8 1 1 1 1 4' crashes, but with the 
corresponding matrix input one gets

\begin{verbatim}
palp$ poly.x -C2
Degrees and weights  `d1 w11 w12 ... d2 w21 w22 ...'
  or `#lines #columns' (= `PolyDim #Points' or `#Points PolyDim'):
4 5
Type the 20 coordinates as dim=4 lines with #pts=5 columns:
   -1    7   -1   -1   -1
   -1   -1    7   -1   -1
   -1   -1   -1    7   -1
   -1   -1   -1   -1    1
pic=1  deg=64  h12= 0  rk=0 #sq=0 #dp=0 py=1  F=5 10 10 5 #Fano=1
4 5  Vertices of P* (N-lattice)    M:201 5 N:7 5
 1  0  0  0 -1
 0  0  1  0 -1
 0  1  0  0 -1
 0  0  0  1 -4
P/2: 36 points (5 vertices) of P'=P/2 (M-lattice):
P/2:  0  4  0  0  0   0  0  0  0  0  0  0  0  0  0  0  0  1  1  1 ... 
P/2:  0  0  4  0  0   1  2  3  0  1  2  3  0  1  2  0  1  0  1  2 ... 
P/2:  0  0  0  4  0   0  0  0  1  1  1  1  2  2  2  3  3  0  0  0 ... 
P/2:  0  0  0  0  1   0  0  0  0  0  0  0  0  0  0  0  0  0  0  0 ... 
\end{verbatim}

\opt{-E}
{\tt poly.x -E} checks several symmetry properties of a reflexive polytope related to the set of roots of the associated toric variety. These are of interest with respect to the existence of Einstein-Kaehler metrics. Here, a root is a lattice point in the interior of a facet of the reflexive polytope (in the $M$ lattice). Centrally-symmetric roots are called semisimple. If all roots are semisimple, then {\tt ssroot=1}. Only in this case PALP proceeds to check and display the following conditions: if the barycenter is $0$ ({\tt bary=1}), the sum of lattice points is $0$ ({\tt $\#$Psum=1}), the sum of lattice points in each multiple is $0$ ({\tt $\#$kPsum=1}), the group of lattice automorphisms has only the origin as a fixed point ({\tt $\#$symm=1}). These conditions can be found explained in Chapter 5 (in particular, sections 5.5 and 5.6) of the dissertation of Benjamin Nill \cite{nillthesis}.
This option cannot be combined with others; if {\tt -E} is specified, all other options are ignored.
%%%%%%%%%%%%%%%%%%%%%%%%%%%%%%%%%%%%%%%%%%%%%%%%%%%%%%%%%%%%%%%%%%%%%%%%%%%%%%
\section{{\tt cws.x}}
\label{sec-cws}
\subsection{General description of {\tt cws.x}}

\cws\ is concerned mainly with the first steps of the algorithm 
of \cite{Kreuzer:1995cd,Skarke:1996hq,Kreuzer:2000qv} for the classification of
reflexive polytopes in a given dimension.
In particular it contains an implementation of the algorithm 
\cite{Skarke:1996hq} for the classification of weight systems, and routines for
combining these weight systems into CWS.
As always, rough information can be obtained with the help screen.

\begin{verbatim}
palp$ cws.x -h
This is `cws.x': create weight systems and combined weight systems.
Usage:    cws.x -<options>; 
          the first option must be `w', `c', `i', `d' or `h'.
Options:
-h        print this information
-f        use as filter; otherwise parameters denote I/O files
-w# [L H] make IP weight systems for #-dimensional polytopes.
          For #>4 the lowest and highest degrees L<=H are required.
    -r/-t make reflexive/transversal weight systems (optional).
-c#       make combined weight systems for #-dimensional polytopes.
          For #<=4 all relevant combinations are made by default,
          otherwise the following option is required:
    -n[#] followed by the names wf_1 ... wf_# of weight files
          currently #=2,3 are implemented.
     [-t] followed by # numbers n_i specifies the CWS-type, i.e.
          the numbers n_i of weights to be selected from wf_i.
          Currently all cases with n_i<=2 are implemented.
 -i       compute the polytope data M:p v [F:f] N:p [v] for all IP
          CWS, where p and v denote the numbers of lattice points
          and vertices of a dual pair of IP polytopes; an entry  
          F:f and no v for N indicates a non-reflexive `dual pair'.
 -d#      compute basic IP weight systems for #-dimensional reflexive
          Gorenstein cones;
     -r#  specifies the index as #/2.
 -2       adjoin a weight of 1/2 to the input of the weight system.
 -N       make CWS for PPL in N lattice.
\end{verbatim}
There is also a second help screen that can be called with the option {\tt -x}
(see below).

\subsection{Options of {\tt cws.x}}
\subsubsection{-w[number]}
The behaviour of {\tt cws.x -w\#} depends crucially on \#.

If $\#\le 4$ all weight systems corresponding to \#-dimensional IP-simplices 
are determined by executing the algorithm of \cite{Skarke:1996hq}:

\begin{verbatim}
palp$ cws.x -w2
3  1 1 1  rt
4  1 1 2  rt
6  1 2 3  rt  #=3  #cand=3
\end{verbatim}
The algorithm determines candidates for weight systems
and prints them if they lead to polytopes with the IP property (see section
\ref{ipps});
this holds for all 3 candidates, as {\tt \#=3  \#cand=3} indicates.
If such a weight system gives rise to a reflexive polytope
(which is always the case in dimension $\le 4$ \cite{Skarke:1996hq}) 
this is indicated by an {\tt r};
if the (possibly singular) weighted projective space corresponding to the 
weight system obeys the `transversality condition' that the Calabi--Yau 
hypersurface equation introduces no additional singularities, this is 
indicated by a {\tt t}.

If $\# > 4$, one has to enter a lower and an upper bound for the degrees of 
the weight systems. {\tt cws.x -w\#} then examines all possible such systems 
and displays the ones that define polytopes with the IP property.

If an extra option of {\tt -r} or {\tt -t} is specified, the output contains only
the reflexive or transverse weight systems, respectively.
Just try {\tt cws.x -w5 5 8}, {\tt cws.x -w5 5 8 -r} and 
{\tt cws.x -w5 5 8 -t} to see how this works.

\subsubsection{-c[number]}
Now the output contains combined weight systems (CWS).
Again all of them are created if the number after {\tt -c} is $\le 4$
(try {\tt cws.x -c3}).
Otherwise weight systems that are read from files are combined.
We apologize for not being able to give information beyond the one given in
the help screen.

\subsubsection{-i}
In this case polytope input is required.
The output is like that of {\tt poly.x -g}, but suppressed for polytopes
without the IP property. This can be useful to filter a list of
CWS for the IP property.

\subsubsection{-d[number] [-r[number]]}
The so-called basic weight systems for reflexive Gorenstein cones of a given 
dimension (the number after {\tt -d}) and a given index are computed; 
if {\tt -r} is used, the index is \emph{half} the number after {\tt -r}, otherwise
the index is 1 by default. 
See \cite{Skarke:2012zg,cydata} for more details.

\subsubsection{-2}
{\tt cws.x -2 < infile > outfile} writes the list of weight systems 
in {\tt infile} to {\tt outfile}, 
but with a weight of $1/2$ adjoined to each input weight system; this is useful 
because the {\tt -d}--option produces only weight systems without weights of $1/2$.
\subsubsection{-N}
\label{sec:cws-n}
Given a polytope as matrix input, this option reconstructs the CWS of a
%polytope isomorphic to the given one. 
reflexive pair $(P,P^*)$ with $P^*$ isomorphic to the input polytope.
The option only accepts matrix
input, otherwise it returns the message {\tt Only PPL-input in
  Npoly2cws!}. It is useful to verify whether a polytope is
actually defined on a sublattice of finite index. As a simple example
consider the polytope consisting of only the vertices of the Newton
polytope of the quintic hypersurface in $\mP^4$
\begin{verbatim}
palp$ cws.x -N
Degrees and weights  `d1 w11 w12 ... d2 w21 w22 ...'
  or `#lines #columns' (= `PolyDim #Points' or `#Points PolyDim'):
4 5
Type the 20 coordinates as dim=4 lines with #pts=5 columns:
-1 -1 -1 -1  4
-1 -1 -1  4 -1
-1 -1  4 -1 -1
-1  4 -1 -1 -1
5 1 1 1 1 1 /Z5: 4 1 0 0 0 /Z5: 4 0 1 0 0 /Z5: 4 0 0 1 0 
\end{verbatim}
This option is in some sense inverse to {\tt poly.x} without any options,
see Section~\ref{sec-poly}, and hence should viewed as part of {\tt poly.x} which
for historical reasons ended up in {\tt cws.x}.

\subsubsection{-x}
A further help screen with additional options is displayed:
\begin{verbatim}
palp$ cws.x -x
This is `cws.x': -x gives undocumented extensions:
              -ip    printf PolyPointList
              -id    printf dual PolyPointList
              -p#    [infile1] [infile2] makes cartesian product
                     of Vertices. # dimensions are identified.
              -S     count simplex points for weight system
              -L     count using LattE (-> count redcheck cdd)
\end{verbatim}
As {\tt -x} refers to `experimental' and none of the authors is familiar with 
them, we leave it to the reader to play with them
and perhaps find useful applications. It would be greatly appreciated
if any insights gained in this way were communicated via the PALP Wiki \cite{wiki}.
%%%%%%%%%%%%%%%%%%%%%%%%%%%%%%%%%%%%%%%%%%%%%%%%%%%%%%%%%%%%%%%%%%%%%%%%%%%%%
\section{\tt class.x}\label{sec-class}
\subsection{General description of {\tt class.x}}
\class\ implements the actual creation of the complete lists of reflexive
polytopes in dimensions up to 4 (see \cite{Kreuzer:1998vb,Kreuzer:2000xy}).
It contains routines for determining all subpolytopes of a 
given polytope in an efficient manner, for finding all lattices on which a 
given polytope is reflexive, and for several other tasks relevant to the
classification.
In particular, as the resulting numbers of polytopes tend to get very large, 
\class\ can encode a corresponding list in various binary formats.
This can take the form of a single binary file or 
a collection of binary files to which we refer as a database (even though it 
is not a relational database in the usual sense).
Actually there are two types of database: the first type is used in the 
classification procedure and contains only information on the polytopes 
themselves, whereas the second type also contains Hodge number information;
the latter is accessed by the website \cite{cydata}.
\class\ contains routines for converting the various formats and for
creating set theoretic unions or differences of the corresponding lists.

It is very important to use \class\ with the appropriate parameter--setting
of {\tt POLY\_Dmax = 4}.
To get a good feeling for how \class\ functions, we recommend to use it
for rederiving the classification of reflexive 3-polytopes following
section 3.2 of \cite{Kreuzer:2002uu}. The available options are listed in the help screen:
\begin{verbatim}
palp$ class.x -h
This is  `class.x', a program for classifying reflexive polytopes
Usage:     class.x  [options] [ascii-input-file [ascii-output-file]]
Options:   
-h          print this information
-f or -     use as filter; otherwise parameters denote I/O files
-m*         various types of minimality checks (* ... lvra)
-p* NAME    specification of a binary I/O file (* ... ioas)
-d* NAME    specification of a binary I/O database (DB) (* ... ios)
-r          recover: file=po-file.aux, use same pi-file
-o[#]       original lattice [omit up to # points] only
-s*         subpolytopes on various sublattices (* ... vphmbq)
-k          keep some of the vertices
-c, -C      check consistency of binary file or DB
-M[M]       print missing mirrors to ascii-output
-a[2b], -A  create binary file from ascii-input
-b[2a], -B  ascii-output from binary file or DB
-H*         applications related to Hodge number DBs (* ...cstfe)

Type one of [m,p,d,r,o,s,c,M,a,b,H] for help on options,
`g' for general help, `I' for general information on I/O or `e' 
to exit:
\end{verbatim}
As the last lines show, this help screen is interactive, i.e.
by typing one of the option letters within the help program one can obtain 
further information.
For this reason, and also because there are probably not too many potential
users of this program, we shall be rather brief in the following descriptions.

\subsection{Options of \class}

In the following we assume that 
the reader also consults the help text 
of {\tt class.x -h} followed by the letter for the corresponding option;
otherwise our explanations may not make much sense or 
important details can be missed.

\subsubsection{-h}
The interactive help screen is displayed (see above).

\subsubsection{-f}
The filter flag turns off the prompt for input from the screen.

\subsubsection{-ml, -mv, -mr, -ma}
In \cite{Kreuzer:1995cd,Kreuzer:1998vb,Kreuzer:2000xy} several 
versions of the concept of a minimal polytope were defined 
(see \cite{Kreuzer:2000qv} for a summary).
Given polytope input, it is determined whether these definitions are satisfied.
This is the only option except {\tt -h} for which both input and output are ascii.

\subsubsection{-pi, -pa, -ps, -po}
{\tt -p*} specifies the name of a binary file.
Depending on the second letter, this may be an input file, possibly such that 
the list it encodes should be added to or subtracted from some other list,
or an output file.

\subsubsection{-di, -ds, -do}
Like {\tt -p*}, but now the name of a database is specified.

\subsubsection{-r}
This option was used in the classification of reflexive 4-polytopes to
recover intermediate results after computer crashes.

\subsubsection{-o, -o[number], -oc}
{\tt -o} instructs \class\ to ignore polytopes that are
reflexive only on a sublattice.
With a number following, only up to that number of points are omitted in 
the search for subpolytopes 
(for an example see section 3.1 of \cite{Kreuzer:2002uu}).
{\tt -oc} modifies the behaviour of the previous option {\tt -r}.

\subsubsection{-sh, -sp, -sv, -sm, -sb, -sq}
A given polytope $\Delta$ may be reflexive w.r.t. several distinct 
lattices if the lattice $M_{\rm coarsest}$ generated by the vertices of 
$\Delta$ is not dual to the lattice $N_{\rm coarsest}$ generated by the 
vertices of $\Delta^*$ (see e.g. \cite{Kreuzer:2000qv}).
{\tt -s*} finds such cases in the database; various versions called by 
different letters after {\tt -s} correspond to lattices generated by different
types of subsets of the set of lattice points of $\Delta$.
For example Calabi-Yau hypersurfaces that are free quotients result from
sublattices generated by all points except those interior to facets.

\subsubsection{-k}
With {\tt -k} (for `keep') one can specify some vertices of an input polytope 
$\Delta$ such that \class\ finds all reflexive subpolytopes of $\Delta$
that still contain these vertices.

\subsubsection{-c, -C}
These options perform consistency checks on a given file or database.
This is very useful if one suspects that something may have gone wrong 
with the data.

\subsubsection{-M}
A list given in binary format is checked for mirror symmetry.
Those polytopes that would be required to make the list mirror symmetric 
are displayed.

\subsubsection{-a, -A}
Ascii input is converted to binary file format.
{\tt -a} corresponds to the standard normal form (cf.~section \ref{normalform})
and {\tt -A} to the affine normal form (cf.~section \ref{anf}).

\subsubsection{-b, -B}
Binary input is converted to ascii.
{\tt -b} is the standard version and {\tt -B} should be used for lists created
with {\tt -A}.

\subsubsection{-Hc, -Hs, -Ht, -Hf, -He}
Options of this type are related to the second type of data base which 
contains data on Hodge numbers. They work only in dimension 4.
%%%%%%%%%%%%%%%%%%%%%%%%%%%%%%%%%%%%%%%%%%%%%%%%%%%%%%%%%%%%%%%%%%%%%%%%%%%%%%
\section{{\tt nef.x}}
\label{sec-nef}
\subsection{General Description of {\tt nef.x}}
{\tt nef.x} is a package to analyze nef partitions of a reflexive
polytope. Such nef partitions determine complete intersections of
Calabi-Yau type in toric varieties of, in principle, arbitrary
codimension. Given a  reflexive polytope in terms of a combined weight
system (cf. Section~\ref{sec:polytope-input}) or a list of points the main
objective of the program is to determine the nef partitions and the
Hodge numbers of the corresponding Calabi-Yau
varieties. Further features include the calculation of the corresponding
reflexive Gorenstein cones as well as information about the fibration
structure. A short summary of the available options can be obtained
from the help screen: 
\begin{verbatim} 
palp$ nef.x -h
This is './nef.x':  calculate Hodge numbers of nef-partitions
Usage:   ./nef.x <Options>
Options:
-h        prints this information
-f or -   use as filter; otherwise parameters denote I/O files
-N        input is in N-lattice (default is M)
-H        gives full list of Hodge numbers
-Lv       prints L vector of Vertices (in N-lattice)
-Lp       prints L vector of Points (in N-lattice)
-p        prints only partitions, no Hodge numbers
-D        calculates also direct products
-P        calculates also projections
-t        full time info
-cCODIM   codimension (default = 2)
-Fcodim   fibrations up to codim (default = 2)
-y        prints poly/CWS in M lattice if it has nef-partitions
-S        information about #points calculated in S-Poly
-T        checks Serre-duality          
-s        don't remove symmetric nef-partitions   
-n        prints polytope only if it has nef-partitions
-v        prints vertices and #points of input polytope in one
          line; with -u, -l the output is limited by #points:
    -uPOINTS  ... upper limit of #points (default = POINT_Nmax)
    -lPOINTS  ... lower limit of #points (default = 0)
-m        starts with [d  w1 w2 ... wk d=d_1 d_2 (Minkowski sum)
-R        prints vertices of input if not reflexive
-V        prints vertices of N-lattice polytope
-Q        only direct products (up to lattice Quotient)
-gNUMBER  prints points of Gorenstein polytope in N-lattice
-dNUMBER  prints points of Gorenstein polytope in M-lattice
      if NUMBER = 0 ... no            0/1 info
      if NUMBER = 1 ... no redundant  0/1 info (=default)
      if NUMBER = 2 ... full          0/1 info
-G        Gorenstein cone: input <-> support polytope
\end{verbatim}
Note that for examples of codimension greater than two calculation
times can become very long and parameters specified in the files {\tt
  Global.h} and {\tt Nef.h} may have to be suitably chosen upon
compilation (cf. Section~\ref{error-handling}).

In this section we will begin with a brief reminder of the notion of a nef
partition. Then we will describe in detail the standard output of {\tt nef.x} when called without any options. Finally we will discuss each option in detail and demonstrate its functionality in examples. Further examples and details can be found at the PALP Wiki \cite{wiki}.
%%%%%%%%%%%%%%%%%%%%%%%%%%%%%%%%%%%%%%%%%%%%%%%%%%%%%%%%%%%%%%%%%%%%%%%%%%%%%%%
\subsection{Nef partitions and reflexive Gorenstein cones}
\label{sec:nef-part}

In this subsection, we give a brief summary of nef ('nef' stands for numerically effective) partitions and
related notions. For details,
see~\cite{borisov,Batyrev:1994pg,Kreuzer:2001fu,Batyrev:2007cq}. 

Consider a dual pair of $d$-dimensional reflexive polytopes
$\Delta\subset M_\mR,\Delta^*\subset N_\mR$. A
partition $V=V_0\cup\dots\cup V_{r-1}$ of the set of vertices of
$\Delta^*$ into disjoint subsets $V_0,\dots,V_{r-1}$ is called a
nef partition of length $r$ if there exist $r$ integral upper convex
$\Sigma(\Delta^*)$-piecewise linear support functions
$\phi_l:N_\mR\rightarrow \mR$, $l=0,\dots,r-1$ such that
\begin{equation} 
  \phi_l(v)=
  \begin{cases} 
    1 & \text{if } v \in V_l, \\ 
    0 & \text{otherwise.}
  \end{cases}
  \label{eq:supportfunction}
\end{equation} 
Each $\phi_l$ corresponds to a divisor $D_{0,l}=\sum_{v \in V_l} D_{v}$
on the toric variety $\mP_{\Delta^*}$ associated to $\Delta^*$, and
the intersection of all these divisors  
\begin{equation}
  \label{eq:CICY}
  X=D_{0,0}\cap\dots\cap D_{0,r-1}
\end{equation}
defines a family $X\subset \mP_{\Delta^*}$ of Calabi-Yau complete
intersections of codimension $r$.   

Moreover, each $\phi_l$ corresponds to a lattice
polytop $\Delta_l$ defined as
\begin{equation*}
  \label{eq:di} 
  \Delta_l=
  \Big\{
  x\in M_\mR
  \,\Big|\,
  (x,y)\geq -\phi_l(y)
  \quad \forall y\in N_\mR
  \Big\}
  .
\end{equation*} 
The sum of the functions $\phi_l$ is equal to the
support function of the anticanonical divisor $K_{\mP_{\Delta^*}}^{-1}$ and, therefore, the
corresponding Minkowski sum is $\Delta_0 + \dots + \Delta_{r-1} = \Delta$.
Moreover, the knowledge of the decomposition $V=V_0\cup\dots\cup V_{r-1}$
is equivalent to that of the set of supporting polytopes
$\Pi(\Delta)=\{\Delta_0,\dots,\Delta_{r-1}\}$, and therefore this data is
often also called a nef partition. 

For a given nef partition $\Pi(\Delta)$ the polytopes\footnote{The
  brackets $\big< \cdots \big>$ denote the convex hull.} 
\begin{equation*}
  \nabla_{l'}=
  \big\langle
  \{0\}\cup V_{l'}
  \big\rangle
  ~\subset N_\mR
\end{equation*}
define again a nef partition
$\Pi^\ast(\nabla)=\{\nabla_0,\dots,\nabla_{r-1}\}$ such that the Minkowski
sum $\nabla=\nabla_0+\dots+\nabla_{r-1}$ is a reflexive
polytope. Then, its dual $\nabla^*$ is also reflexive, and
$\Pi^\ast(\nabla)$ is called the dual nef partition. This is the
combinatorial manifestation of mirror symmetry in terms of dual 
pairs of nef partitions of $\Delta^*$ and $\nabla^*$, which we
summarize in the following diagram
\begin{equation}
  \label{eq:polyDN}
  \begin{array}{rcl}
    M_{\mathbb{R}}\qquad\qquad && \qquad\qquad N_{\mathbb{R}}\\
        & & \\
    \Delta=\Delta_0+\ldots+\Delta_{r-1} && \Delta^*=\langle\nabla_0,\ldots,\nabla_{r-1}\rangle\\
        & & \\
        &~~~~ (\Delta_l,\nabla_{l'})\ge-\delta_{l\,l'} ~~~~&\\
        & & \\
    \nabla^*=\langle\Delta_0,\ldots,\Delta_{r-1}\rangle && \nabla=\nabla_0+\ldots+\nabla_{r-1}
  \end{array}
\end{equation}
In the horizontal direction, we have the duality between the lattices
$M$ and $N$ and mirror symmetry goes from the upper right to the lower
left. The complete intersections $X\subset \mP_{\Delta^*}$ and $\widecheck{X}
\subset \mP_{\nabla^*}$ associated to the dual nef partitions are then
mirror Calabi-Yau varieties.

There are two constructions to build new nef partitions from old ones:
projections and direct products.
Given a nef partition $V=V_0\cup\dots\cup V_{r-1} $ 
where one of the subsets $V_l$, say $V_0$,
consists of a single vertex $v$, the nef condition implies that the
projection $\Delta^*_v$ of $\Delta^*$ along $v$ is
reflexive. Moreover, by~(\ref{eq:CICY}) the Calabi--Yau complete
intersection $X$ is given by $D \cap X'$ with $X' = \bigcap_{l=1}^{r-1}
D_{0,l}$. Since $D$ can only intersect the toric divisors that 
correspond to points bounding the reflexive projection along $v$, the
variety $X$ is isomorphic to the variety $X' \subset 
\mP_{\Delta^*_v}$, where $\mP_{\Delta^*_v}$ is obtained from the
projection $\Delta^*_v$. In~\cite{Klemm:2004km} such nef partitions
were called trivial. 
In {\tt nef.x} they are labeled by ${\tt P}$ for projection, 
see Section~\ref{sec-nef-P}.

Suppose we are given two lattices $M^{(1)},M^{(2)}$ and two reflexive polytopes 
$\Delta^{(1)}\subset M^{(1)}_{\mR}$, $\Delta^{(2)}\subset M^{(2)}_{\mR}$ %with
such that $(\Delta^{(1)})^*$ and  $(\Delta^{(2)})^*$ admit
nef partitions $V^{(1)} = \bigcup_l V_l^{(1)}$ and 
$V^{(2)} = \bigcup_l V_l^{(2)}$, respectively. 
Then $\Delta= \Delta^{(1)} %\oplus 
\times\Delta^{(2)}$ is
reflexive with respect to $M=M^{(1)}\oplus M^{(2)}$ and dual to $\Delta^*$
whose set of vertices $V$ is $\{ (v^{(1)},0) \,|\,
v^{(1)}\in V^{(1)}\}  \cup \{ (0,v^{(2)}) \,|\, v^{(2)}
\in V^{(2)}\}$. 
$V$ admits a nef partition induced from the
nef partitions $V^{(1)}$ and $V^{(2)}$. Such a nef partition is
called a direct product since the corresponding Calabi--Yau complete
intersection $X$ is a direct product $X=X^{(1)} \times
X^{(2)}$ in $\mP_{\Delta^*} = \mP_{(\Delta^{(1)})^*} \times
\mP_{(\Delta^{(2)})^*}$. 

One can reformulate the duality of nef partitions in terms of
reflexive Gorenstein cones as follows. We extend the lattices $M$ and
$N$ to $\widetilde M = \mZ^r \oplus M$ and $\widetilde N = \mZ^r \oplus N$ 
and set $\widetilde d = d+r$.

A $\widetilde d$--dimensional rational polyhedral cone $C$
in $\widetilde M_{\mR}$ is called Gorenstein 
if $C\cap(-C)=\{0\}$, there 
exists an element $n_C \in \widetilde N_{\mR}$
such that $\langle x,n_C\rangle > 0$ for any nonzero $x \in C$, and all
vertices of the $(\widetilde d -1)$--dimensional convex polytope
\[\Delta(C)=\{x\in C \,|\, \langle x,n_C\rangle=1\}\] 
belong to
$\widetilde M$. 
The polytope $\Delta(C)$ is called 
the support 
of $C$. 
Conversely, any $(\widetilde d -1)$--dimensional lattice polytope $\Lambda$ 
determines a $\widetilde d$--dimensional Gorenstein cone $C(\Lambda)$ as the cone 
over $\Lambda$ with apex at lattice distance 1 from the hyperplane carrying $\Lambda$;
obviously $\Delta(C(\Lambda))=\Lambda$.
For any $m \in  C \cap  \widetilde M$, we define the degree of $m$ as $\deg m
= \langle m, n_C  \rangle$. 

A Gorenstein cone $C$ is called reflexive if the dual
cone \[\widecheck C  =
\{y \in \widetilde N_\mR \,|\, \langle x,y \rangle \geq 0\; \forall\,
x\in C\}\] is also Gorenstein, i.e., there exists $m_{\widecheck C} \in \widetilde M$
such that $\langle m_{\widecheck C},y\rangle > 0$ for all $y \in  \widecheck C
\setminus \{0\}$, and all vertices of the 
support 
$\Delta ( \widecheck
C ) = \{ y \in \widecheck C \,|\, \langle m_{\widecheck C} , y \rangle = 1 \}$
belong to $\widetilde N$. We will call the integer $r = \langle
m_{\widecheck C} , n_C\rangle$ the index of $C$ (or $\widecheck C$).

Any
nef partition $\Pi(\Delta)=\{\Delta_0,\dots,\Delta_{r-1}\}$
of length $r$ of a reflexive polytope $\Delta$ 
determines
a $\widetilde d$--dimensional 
dual pair of reflexive Gorenstein cones
$C = C(\Delta_1,\dots,\Delta_r)\subset\widetilde M_{\mR}$,
$\widecheck C = \widecheck C(\nabla_1,\dots,\nabla_r)\subset\widetilde N_{\mR}$
of index $r$ by 
\begin{align}
C = \{ (\lambda_1, \dots ,\lambda_r,\lambda_1x_1 + \dots + \lambda_r x_r) \in
\widetilde M_{\mR} \,|\, \lambda_i \geq 0, x_i \in \Delta_i,
i=1,\dots, r\},\nonumber \\
\widecheck C = \{ (\mu_1, \dots ,\mu_r,\mu_1x_1 + \dots + \mu_r x_r) \in
\widetilde N_{\mR} \,|\, \mu_i \geq 0, x_i \in \nabla_i,
i=1,\dots, r\}.\nonumber
\end{align}
There are, 
however, reflexive Gorenstein cones that do not come from nef partitions.

A reflexive Gorenstein cone admits a representation in terms of the
points of the underlying reflexive polytope as follows. Given a 
point
$p \in \nabla_l$, the 
corresponding point $\widetilde{p} \in
\widecheck C(\nabla_1,\dots,\nabla_r)$ is given as
\begin{equation}
  \widetilde p = (\phi_0(p),\dots,\phi_{r-1}(p),p).
  \label{eq:ptilde}
\end{equation}
where $\phi_l$ is the support function~(\ref{eq:supportfunction}).
To see that the two descriptions of $\widecheck C$ are equivalent, note
that both correspond to a cone whose support has vertices
\begin{equation}
(e_{i(1)},v_1),\ldots,(e_{i(n)},v_n),(e_1, 0_N),\ldots,(e_r,0_N),
  \label{eq:supportverts}
\end{equation}
where $\{e_i\}$ is the standard basis of $\mZ^r$, $i(k)$ is the number 
such that $v_k\in V^{(i(k))}$ and $0_N$ is the origin in the N--lattice.

The Hodge numbers of a Calabi--Yau manifold $X$ defined by means of a
nef partition as in~(\ref{eq:CICY}) depend only on the structure of
the corresponding reflexive Gorenstein cone in a manner described in
\cite{Batyrev:1995ca,Kreuzer:2001fu}.  
The corresponding formulas rely heavily on the counting of lattice points.
For any lattice polytope $\Lambda$ let us denote by $\ell(\Lambda)$ the number 
of lattice points of $\Lambda$ and by $\ell^*(\Lambda)$ the number of lattice 
points in the interior of $\Lambda$.
It can be shown that
\begin{equation}
S_\Lambda(t) = (1-t)^{\dim \Lambda + 1}\sum_{k \geq 0} \ell(k\Lambda) t^k
\end{equation}
is a polynomial of degree $d\le \dim \Lambda + 1$;
$S_\Lambda(t)$ is called the Ehrhart polynomial of $\Lambda$.
Similarly one can define a polynomial 
\begin{equation}
T_\Lambda(t) = (1-t)^{\dim \Lambda + 1}\sum_{k \geq 0} \ell^*(k\Lambda) t^k.
\end{equation}
In terms of a Gorenstein cone $C$ over $\Lambda$, with underlying lattice $M_C$,
$S$ and $T$ can be written as
\begin{align}
S(C,t) = (1-t)^{\dim C}\sum_{m \in C \cap M_C} t^{\deg m}, \\
T(C,t) = (1-t)^{\dim C}\sum_{m \in \mathrm{int}(C) \cap M_C} t^{\deg m}.
\end{align}
The two polynomials satisfy a relation which is a consequence of Serre duality,
\begin{equation}
S(C,t) = t^{\dim C}\,T(C,t^{-1}), \label{eq:Serre}
\end{equation}
which provides a stringent test on any results involving lattice point counting.
For the computation of Hodge numbers, the S-- and T-- polynomials for all the
faces of $C(\Delta)$ as well as a polynomial called $B$, which is related 
to the poset structure of $C(\Delta)$, are required.

\subsection{Standard output}
\label{sec:stdout}
In this subsection we will explain in detail how to interpret the output of 
{\tt nef.x} when called without any options.

The standard output slightly depends on whether the reflexive polytope
is input as a combined weight system
or as a collection of points.
If the polytope was entered as a collection of points, the first line
of the output takes the following form: 
\begin{verbatim}
M:# # N:# # codim=# #part=#
\end{verbatim}
Note that the input polytope 
is interpreted as $\Delta\subset M_\IR$ unless the option {\tt -N} (cf. Section
\ref{sec-nef-N}) is used,
while any output of a polytope in matrix format refers to its dual 
$\Delta^*\subset N_\IR$ except for the option {\tt -y} (cf. Section~\ref{sec-nef-y}. 
If the input is a CWS, the line starts with the CWS repeated before the 
letter {\tt  M}.
\begin{verbatim}
# M:# # N:# # codim=# #part=#
\end{verbatim}
where the first {\tt \#} stands for the sequence of numbers describing the
CWS. 
The two numbers {\tt \#} after
{\tt M} correspond to the numbers of 
lattice points and vertices of $\Delta\subset M_\IR$ and the
numbers {\tt  \#} after {\tt N} correspond to the numbers of lattice 
points and vertices of $\Delta^*\subset N_\IR$, 
respectively. 
The number $r$ in
{\tt codim=r} is the length of the nef partition, i.e. the codimension
of the corresponding Calabi--Yau complete intersection. The
default value is 2, otherwise it is specified by the option {\tt -c*}
described in Section~\ref{sec-nef-c*}. The number $n$ in {\tt
  \#part=n} is the number of all the nef partitions that {\tt nef.x}
has found, up to symmetries of the underlying lattice. If the
symmetries of the underlying lattice should not be taken into account,
use the option {\tt -s} (cf. Section \ref{sec-nef-s}).  

The subsequent lines contain the information about the various nef
partitions. Note that the standard output suppresses the output of nef
partitions which are equivalent under symmetries of the
CWS. 
If the codimension is 2 the output line containing the information on
a particular nef partition takes the following form: 
\begin{verbatim}
H:# [#] P:# V:#  #      #sec #cpu
\end{verbatim}
The numbers {\tt \#} after {\tt H}: are the Hodge numbers $h^{1,i}(X),\:
i=1,\dots,d-1$, where $d\,$ is the dimension of the Calabi-Yau
manifold $X\,$ obtained via~(\ref{eq:CICY}).

The number {\tt \#} in the square brackets {\tt [\#]} is the Euler
number of $X\,$. If $h^{0,i}(X) \not= 0$ for some $i=1,\dots,d-1$ the
Calabi-Yau manifold factorizes. See the option {\tt -D}
(Section \ref{sec-nef-D}) for this case. In any case, the full Hodge
diamond is displayed with the option {\tt -H} (Section
\ref{sec-nef-H}). 

The number {\tt \#} after {\tt P}: is a counter specifying the
nef partition. It runs from $0$ to $n - 1$. Note that nef partitions
corresponding to direct products and projections to nef partitions of
lower length are omitted by default. To display them use the options
{\tt -D} (cf. Section~\ref{sec-nef-D}), {\tt -Q}
(cf. Section~\ref{sec-nef-Q}) for direct products and {\tt -P} (cf. Section
\ref{sec-nef-P}) for projections.

The sequence of numbers {\tt \#} separated by a single space after {\tt
  V}: corresponds to the vertices that belong to the first part $V_0$ of the
nef partition. Note that the vertices are counted starting from
$0$. These numbers only make sense if the options {\tt -n}
(cf. Section~\ref{sec-nef-n}), {\tt -Lv}
(cf. Section \ref{sec-nef-Lv}) or {\tt -Lp}
(cf. Section~\ref{sec-nef-Lp}) are used. The vertices that belong to
the second part $V_1$ of the nef partition are not displayed, since
they are simply the remaining ones. If the polytope entered also has
points that are not vertices and if the option {\tt -Lp} is used, then
the second sequence of numbers {\tt \#} that is separated from the
first sequence by two spaces corresponds to the non-vertex points that
belong to the first part $V_0$. For representations of the nef
partition in terms of the Gorenstein cone see the option {\tt -g*}
(cf. Section~\ref{sec-nef-g*}). 

The number {\tt \#} before {\tt sec} indicates the time that was
needed to compute this partition. The number {\tt \#} before {\tt cpu}
indicates the number of CPU seconds that were needed to compute the
Hodge numbers. For determining the nef partitions without computing
the Hodge numbers see the option {\tt -p} (cf. Section \ref{sec-nef-p}).

If the length $r$ is bigger than $2$ the lines containing the information about the various nef partitions take the following form:
\begin{verbatim}
H:# [#] P:# V0:#  # V1:#  # ... V(r-2):#  #      #sec #cpu
\end{verbatim}
Now, there are $r - 1$ expressions of the form {\tt Vi:\#  \#}, where
$i$ runs from $0$ to $r - 2$ each representing a part $V_i$ of the nef
partition. The points and vertices in each $V_i$ are listed in the
same order as in the codimension two case.

The final line of the output always takes the following form:
\begin{verbatim}
np=# d:# p:#    #sec     #cpu
\end{verbatim}
The numbers {\tt\#}  after {\tt d:}, {\tt p:},  {\tt np=} are the numbers
of nef partitions which are direct products, projections, and neither
of the two, respectively. The total of the three numbers adds up to $n$, the total
number of nef partitions as indicated in the first line after {\tt
  \#part=}.

The following example illustrates the standard output of {\tt nef.x}. We
consider complete intersections of codimension 2 in
$\mathbb{P}^2\times\mathbb{P}^1\times\mathbb{P}^2$ discussed in
\cite{Braun:2007vy}. Let $e_1,\dots,e_5$ be the standard basis of $\mathbb{R}^5$. 
We define the polytope $\Delta^* \subset N$ by the vertices
$v_0,\dots,v_7$ given by
\begin{equation*}
  \begin{aligned}
    v_0 &= e_1,& v_1 &= e_2,& v_2 &= -e_1-e_2,& v_3 &= e_3, \\
    v_4&=-e_3,& v_5&=e_4,& v_6&=e_5,&v_7&=-e_4-e_5.
  \end{aligned}
\end{equation*}
By elementary toric geometry, we see that $\mP_{\Delta^*} =
\mathbb{P}^2\times\mathbb{P}^1\times\mathbb{P}^2$ and the combined
weight system can be read off from the linear relations
\begin{equation*}
  \begin{aligned}
    v_0 + v_1 + v_2 &=0, & v_3+v_4 &=0, & v_5 + v_6 + v_7 &=0.
  \end{aligned}
\end{equation*}
First, we enter the polytope by giving this combined weight system 
\begin{verbatim}
palp$ nef.x 
Degrees and weights  `d1 w11 w12 ... d2 w21 w22 ...'
  or `#lines #colums' (= `PolyDim #Points' or `#Points PolyDim'):
3 1 1 1 0 0 0 0 0  2 0 0 0 1 1 0 0 0  3 0 0 0 0 0 1 1 1
3 1 1 1 0 0 0 0 0  2 0 0 0 1 1 0 0 0  3 0 0 0 0 0 1 1 1 
M:300 18 N:9 8  codim=2 #part=15
H:19 19 [0] P:0 V:2 4 6 7       1sec  0cpu
H:9 27 [-36] P:2 V:3 4 6 7       1sec  0cpu
H:3 51 [-96] P:3 V:3 5 6 7       1sec  1cpu
H:3 75 [-144] P:4 V:3 6 7       1sec  0cpu
H:3 51 [-96] P:6 V:4 5 6 7       2sec  1cpu
H:3 51 [-96] P:7 V:4 5 6       1sec  1cpu
H:6 51 [-90] P:8 V:4 6 7       1sec  1cpu
H:3 75 [-144] P:9 V:4 6       1sec  1cpu
H:3 60 [-114] P:10 V:5 6 7       2sec  1cpu
H:3 69 [-132] P:11 V:5 6       1sec  1cpu
H:3 75 [-144] P:12 V:6 7       1sec  0cpu
np=11 d:2 p:2    0sec     0cpu
\end{verbatim}
Equivalently, we can use the option {\tt -N} and enter the points of
the polytope $\Delta^*$ of the normal fan of
$\mathbb{P}^2\times\mathbb{P}^1\times\mathbb{P}^2$: 
\begin{verbatim}
palp$ nef.x -N
Degrees and weights  `d1 w11 w12 ... d2 w21 w22 ...'
  or `#lines #colums' (= `PolyDim #Points' or `#Points PolyDim'):
5 8
Type the 40 coordinates as dim=5 lines with #pts=8 colums:
 1  0 -1  0  0  0  0  0
 0  1 -1  0  0  0  0  0
 0  0  0  1 -1  0  0  0
 0  0  0  0  0  1  0 -1
 0  0  0  0  0  0  1 -1
M:300 18 N:9 8  codim=2 #part=15
H:3 51 [-96] P:0 V:2 3 4 7       1sec  1cpu
H:3 51 [-96] P:1 V:2 4 6 7       2sec  1cpu
H:3 60 [-114] P:2 V:2 4 7       2sec  1cpu
H:3 51 [-96] P:3 V:2 6 7       1sec  1cpu
H:3 69 [-132] P:4 V:2 7       1sec  1cpu
H:9 27 [-36] P:5 V:3 4 6 7       1sec  0cpu
H:3 75 [-144] P:6 V:3 4 7       0sec  0cpu
H:19 19 [0] P:8 V:4 5 6 7       1sec  0cpu
H:6 51 [-90] P:9 V:4 6 7       1sec  1cpu
H:3 75 [-144] P:10 V:4 7       1sec  0cpu
H:3 75 [-144] P:13 V:6 7       1sec  1cpu
np=11 d:2 p:2    0sec     0cpu
\end{verbatim}
Note that both the points and the nef partitions are given in
different orders. 
The polytope $\Delta^* \subset N_{\mR}$ has $9$ points, $8$ vertices
and the interior point, while the dual polytope $\Delta \subset
M_{\mR}$ has $300$ points, $18$ of which are vertices. The codimension
is $2$ and there are $15$ nef partitions.
There are $11$ nef partitions listed, furthermore there are $2$ nef partitions which are direct products, and $2$ which are projections. 
According to the output the nef partitions e.g. $0$ and $8$ are given as follows (with the Hodge numbers and the Euler number of the corresponding Calabi-Yau 3-fold $X$):
\begin{equation}
\begin{aligned}
0:\;& V_0=\langle v_2, v_3, v_4, v_7 \rangle,\quad V_1=\langle v_0, v_1, v_5, v_6 \rangle \\
& h^{1,1}(X) = 3,\, h^{2,1}(X) = 51,\, \chi(X) = -96.\\
\dots &\\
8:\;& V_0=\langle v_4, v_5, v_6, v_7 \rangle,\quad V_1=\langle v_0,
v_1, v_2, v_3 \rangle,\\ & h^{1,1}(X) = 19,\, h^{2,1}(X) = 19,\,
\chi(X) = 0.\nonumber\\
\dots
\end{aligned}\label{eq:2}
\end{equation}
%%%%%%%%%%%%%%%%%%%%%%%%%%%%%%%%%%%%%%%%%%%%%%%%%%%%%%%%%%%%%%%%%%%%%%%%%%%%%%%
\subsection{Options of {\tt nef.x}} 
In this subsection we will explain all the options of {\tt nef.x} in
the order of their appearance in the help screen. Here is a rough guide in terms of specific topics:
\begin{itemize}
  \item Polytope structure: {\tt -N}, {\tt -Lv}, {\tt -Lp},  {\tt -v},
    {\tt -R}, {\tt -V}
  \item Input control: {\tt -N}, {\tt -c*}, {\tt -m}
  \item Structure of nef partitions: {\tt -D}, {\tt -p}, {\tt -P},
    {\tt -s}, {\tt -m} 
  \item Hodge numbers:  {\tt -H}, {\tt -t}, {\tt -S}, {\tt -T}
  \item CWS: {\tt -N}, {\tt -Lv}, {\tt -Lp}, {\tt -m}
  \item Fibrations: {\tt -F*}
  \item Gorenstein cone: {\tt -g*}, {\tt -d*}, {\tt -S}, {\tt -T},
    {\tt -G}
  \item Diagnostics: {\tt -t}, {\tt -S}, {\tt -T}
  \item Polytope statistics: {\tt -y}, {\tt -n}, {\tt -v}, {\tt -R} 
\end{itemize}

\subsubsection{-h}
\label{sec-nef-h}
This option prints the help screen.
%%%%%%%%%%%%%%%%%%%%%%%%%%%%%%%%%%%%%%%%%%%%%%%%%%%%%%%%%%%%%%%%%%%%%%%%%%%%%
\subsubsection{-f or -}
\label{sec-nef-f}
This option switches off the prompt for the input data. This is useful for building pipelines.
%%%%%%%%%%%%%%%%%%%%%%%%%%%%%%%%%%%%%%%%%%%%%%%%%%%%%%%%%%%%%%%%%%%%%%%%%%%%%
\subsubsection{-N}
\label{sec-nef-N}
The option {\tt -N} interprets the input polytope in the N-lattice
instead of the M-lattice. The following example of a complete intersection of degree $(2,2)$ in $\mathbb{P}^3$ illustrates the difference. In order to point out the difference we display the points in the two lattices with the option {\tt -Lv}.
\begin{verbatim}
palp$ nef.x -Lv
Degrees and weights  `d1 w11 w12 ... d2 w21 w22 ...'
  or `#lines #colums' (= `PolyDim #Points' or `#Points PolyDim'):
3 4 
Type the 12 coordinates as dim=3 lines with #pts=4 colums:
   -1    0    0    1
   -1    0    1    0
   -1    1    0    0
M:5 4 N:35 4  codim=2 #part=0
3 4 Vertices in N-lattice:
   -1   -1   -1    3
   -1   -1    3   -1
   -1    3   -1   -1
--------------------
    1    1    1    1  d=4  codim=0
np=0 d:0 p:0    0sec     0cpu
\end{verbatim}
Without the option {\tt -N}, the output 
polytope with 35 points and no nef partition is the dual of the input polytope. 
\begin{verbatim}
palp$ nef.x -Lv -N
Degrees and weights  `d1 w11 w12 ... d2 w21 w22 ...'
  or `#lines #colums' (= `PolyDim #Points' or `#Points PolyDim'):
3 4 
Type the 12 coordinates as dim=3 lines with #pts=4 colums:
   -1    0    0    1
   -1    0    1    0
   -1    1    0    0
M:35 4 N:5 4  codim=2 #part=2
3 4 Vertices in N-lattice:
   -1    0    0    1
   -1    0    1    0
   -1    1    0    0
--------------------
    1    1    1    1  d=4  codim=0
H:[0] P:0 V:2 3   (2 2)     0sec  0cpu
np=1 d:0 p:1    0sec     0cpu
\end{verbatim}
With the option {\tt -N}, the output polytope is the same as input
polytope with $4$ points and the expected nef partition corresponding
to the complete intersection of degree $(2,2)$ in $\mathbb{P}^3$. Note
that the order of the points in the output is the same as in the
input. This last feature is the main advantage of the option {\tt
  -N}. 
The reason is that the basis chosen does not respect the order given by the
combined weight system that was entered. This can be extremely inconvenient at times. The option {\tt -N} provides a way to work around this issue: first use the option {\tt -Lv} to obtain the vertices for a given CWS. Then reorder them so that the basis of linear relations complies with the input and enter the reshuffled vertices into {\tt nef.x} using the option {\tt -N}. This will force the linear relations chosen by {\tt nef.x} to be the same as the CWS.

\subsubsection{-H}
\label{sec-nef-H}
The option {\tt -H} replaces the output lines starting with {\tt H}:
with the full Hodge diamond of the corresponding partition. Note that
the information about the nef partitions is omitted. The following
example of codimension $2$ complete intersections in $\mathbb{P}^7$
illustrates this option (increase {\tt POLY\_Dmax} to $7$):
\begin{verbatim}
palp$ nef.x -H
Degrees and weights  `d1 w11 w12 ... d2 w21 w22 ...'
  or `#lines #colums' (= `PolyDim #Points' or `#Points PolyDim'):
7 1 1 1 1 1 1 1
7 1 1 1 1 1 1 1 M:1716 7 N:8 7  codim=2 #part=3


                               h 0 0   
                          h 1 0      h 0 1   
                     h 2 0      h 1 1      h 0 2   
                h 3 0      h 2 1      h 1 2      h 0 3   
           h 4 0      h 3 1      h 2 2      h 1 3      h 0 4   
                h 4 1      h 3 2      h 2 3      h 1 4   
                     h 4 2      h 3 3      h 2 4   
                          h 4 3      h 3 4   
                               h 4 4   


                                  1
                             0         0
                        0         1         0
                   0         0         0         0
              1       237       996       237         1
                   0         0         0         0
                        0         1         0
                             0         0
                                  1
     16sec  15cpu
[analogous output for a second nef partition]
\end{verbatim}

\subsubsection{-Lv}
\label{sec-nef-Lv}
The option {\tt -Lv} prints the vertices 
of 
$\Delta^*\subset N_\IR$ and the non-negative
linear relations among them in addition to the standard output. If
only the vertices should be printed use the option {\tt -V} in Section
\ref{sec-nef-V}. The output takes the following form: The part before
the dashed line reads:
\begin{verbatim}
 D n Vertices in N-lattice:
    #    #      ...      #    #
    .    .      ...      .    .  
    .    .      ...      .    . 
    #    #      ...      #    #
\end{verbatim}
The first line means that  
$\Delta^*$ has dimension $D\,$
and is given by $n\,$ vertices which are the columns of the subsequent
$D \times n$ array of numbers {\tt \#}. 

Below the dashed line the non-negative
linear relations among these vertices are indicated as follows: 
Let $ v_0, \dots, v_{n-1}$ denote the $n\,$ vertices corresponding to the 
$n\,$ columns above the dashed line. 
Any IP simplex (cf. Section~\ref{ipps}) with vertices in 
$\{v_0, \dots, v_{n-1}\}$ determines a linear relation
$\sum_{i=0}^{n-1} l_i\,v_i = 0$, 
with $l_i\,$ that are positive for the vertices of the IP simplex and $0$ 
otherwise.
It results in
an output line of the form
\begin{verbatim}
l_0 l_1 ... l_{n-1} d=l codim=c
\end{verbatim}
where $l = \sum_{i=0}^{n-1} l_i$ is the degree of the linear relation and
$c$ is the codimension of the IP simplex.
In other words, these lines give 
the set of generators of the cone of non-negative linear relations within
the ($n-D$)--dimensional
vector space of linear relations among the vertices. 
This set is completely fixed by the
order of the vertices, and the conditions that each vector, i.e. each
linear relation, is positive and primitive. 

The information contained in these lines can be very useful in
conjunction with the option {\tt -F*} (cf. Section \ref{sec-nef-F*}). 
 To suppress them see the option {\tt -V} (cf. Section~\ref{sec-nef-V}). 

Besides the standard output the degrees of the nef
partition with respect to the linear relations are inserted in the
output lines containing the information about the nef
partitions as follows. Consider a nef partition of length $r$ defined by $r$ collections of vertices $V_0,\dots, V_{r-1}$ such that every vertex is a member of some collection $V_j$. The (multi)degree of the nef partition $\{V_0,\dots,V_{r-1}\}$ with respect to the linear relation $\sum_{i=0}^{n-1} l_i\,v_i = 0$ is the vector $(d_0,\dots,d_{r-1})$ where $d_j = \sum_{i:v_i \in V_j} l_i$. 
Note that $\sum_{j=1}^r d_j = l$, the degree of the linear relation. The degrees $(d_0,\dots,d_{r-1})$ are the degrees of the polynomials defining the complete intersection. If the codimension is 2 the output lines describing the nef partitions take the following form
\begin{verbatim}
H:# [#] P:# V:#  #   (d10 d11) ... (dn0 dn1)    #sec #cpu
\end{verbatim}
or if the codimension $r$ is bigger than $2$
\begin{verbatim}
H:# [#] P:# V0:#  # V1:#  # ... V(r-2):#  #
      (d10 ... d1(r-1)) ... (dn0 ... dn(r-1)) #sec #cpu
\end{verbatim}
The additional data is {\tt(d10 d11) ... (dn0 dn1)}  and {\tt (d10
  ... d1(r-1)) ... (dn0 ... dn(r-1))}, respectively, where $n$ is the
number of linear relations. If $d_i = (d_{i0},\dots,d_{i,r-1})$ are
the degrees with respect to the $i$-th linear relation, then {\tt
  di0}~$ = d_{i0},\ldots$, {\tt di(r-1)}~$ = d_{i,r-1}$.
The following example of a codimension $2$ complete intersection taken from \cite{Klemm:2004km} illustrates this option:
\begin{verbatim}
palp$ nef.x -Lv
Degrees and weights  `d1 w11 w12 ... d2 w21 w22 ...'
  or `#lines #colums' (= `PolyDim #Points' or `#Points PolyDim'):
5 1 1 1 1 1 0 0  4 0 0 0 1 1 1 1 
5 1 1 1 1 1 0 0  4 0 0 0 1 1 1 1 M:378 12 N:8 7  codim=2 #part=8
5 7 Vertices in N-lattice:
    0   -1    0    1    0    0    0
    0   -1    1    0    0    0    0
   -1    0    0    0    0    0    1
   -1    1    0    0    1    0    0
   -1    1    0    0    0    1    0
-----------------------------------
    1    1    1    1    0    0    1  d=5  codim=1
    1    0    0    0    1    1    1  d=4  codim=2
H:2 64 [-124] P:0 V:0 6   (2 3) (2 2)     1sec  0cpu
[standard output for the remaining Hodge data and nef partitions]
\end{verbatim}
From the output we deduce that the $7$ vertices of the $5$-dimensional polytope satisfy the following 
linear relations:
\begin{equation*}
v_0 + v_1 + v_2 + v_3 + v_6 = 0,\quad  v_0 + v_4 + v_5 + v_6 = 0.
\end{equation*}
The first of these linear relations has degree $5$, the second has
degree $4$. The corresponding IP simplices have codimension $1$ and $2$, respectively.

\subsubsection{-Lp}
\label{sec-nef-Lp}
The option {\tt -Lp} prints all the points of the N-lattice polytope and the
linear relations among them, including those that are not vertices. The output has the same structure as for
the option {\tt -Lv}. The points are ordered such that first the
vertices $\{v_0,\ldots,v_k\}$ are listed, then the points $\{p_{k+1},\ldots,p_{N-2}\}$ which are not vertices and
finally the origin $p_{N-1}$. Note that there will be additional
linear relations including
the points which are neither vertices nor the origin. The following example of a
codimension $2$ complete intersection taken from~\cite{Klemm:2004km}
illustrates this option: 
\begin{verbatim}
palp$ nef.x -Lp
Degrees and weights  `d1 w11 w12 ... d2 w21 w22 ...'
  or `#lines #columns' (= `PolyDim #Points' or `#Points PolyDim'):
5 1 1 1 1 1 0 0  10 2 2 2 2 0 1 1
5 1 1 1 1 1 0 0  10 2 2 2 2 0 1 1 M:378 6 N:8 6  codim=2 #part=4
5 8  Points of Poly in N-Lattice:
   -1    0    0    0    1    0    0    0
   -1    0    1    0    0    0    0    0
   -1    0    0    1    0    0    0    0
   -1    2    0    0    0    0    1    0
   -1    1    0    0    0    1    1    0
----------------------------------------
    2    1    2    2    2    1    0  d=10  codim=0
    1    0    1    1    1    0    1  d=5  codim=1
H:2 86 [-168] P:0 V:1 5  6   (2 8) (1 4)     2sec  2cpu
H:2 68 [-132] P:1 V:2 3 4   (6 4) (3 2)     1sec  0cpu
H:2 68 [-132] P:2 V:3 4   (4 6) (2 3)     1sec  0cpu
np=3 d:0 p:1    0sec     0cpu
\end{verbatim}
The last two points are not vertices. There is one more linear
relation including the point $p_6$. 

\subsubsection{-p}
\label{sec-nef-p}
The option {\tt -p} computes the nef partitions without the
(time-consuming) calculation of Hodge numbers. As an example we
consider the codimension 4 (cf. Section \ref{sec-nef-c*}) complete intersections in $\mP^7$. Note that one must set {\tt POLY\_Dmax} in {\tt Global.h} to at least $10$. 
\begin{verbatim}
palp$ nef.x -c4 -p
Degrees and weights  `d1 w11 w12 ... d2 w21 w22 ...'
  or `#lines #columns' (= `PolyDim #Points' or `#Points PolyDim'):
8 1 1 1 1 1 1 1 1
8 1 1 1 1 1 1 1 1 M:6435 8 N:9 8  codim=4 #part=5
 P:0 V0:2 3  V1:4 5  V2:6 7       0sec  0cpu
np=1 d:0 p:4    0sec     0cpu
\end{verbatim}
The Hodge data in the line containing the partition information is
omitted, and the computation time is 0. Without the option {\tt -p}
this line would look like this:
\begin{verbatim}
H:1 65 [-128] P:0 V0:2 3  V1:4 5  V2:6 7       13127sec  13120cpu
\end{verbatim}
Note the computation time.
%%%%%%%%%%%%%%%%%%%%%%%%%%%%%%%%%%%%%%%%%%%%%%%%%%%%%%%%%%%%%%%%%%%%%%%%%%%%%
\subsubsection{-D}
\label{sec-nef-D}
The option {\tt -D} also prints those nef partitions which are direct
products of lower-dimensional nef partitions. If only direct products
are to be printed use the option {\tt -Q} described in Section~\ref{sec-nef-Q}. 
As an example we consider a codimension $2$ complete intersection in
$\mathbb{P}^2\times \mathbb{P}^2$: 
\begin{verbatim}
palp$ nef.x -D
3 1 1 1 0 0 0 3 0 0 0 1 1 1
3 1 1 1 0 0 0  3 0 0 0 1 1 1 M:100 9 N:7 6  codim=2 #part=5
H:4 [0] h1=2 P:0 V:2 3 5   D     0sec  0cpu
H:20 [24] P:1 V:3 4 5       0sec  0cpu
H:20 [24] P:2 V:3 5       0sec  0cpu
H:20 [24] P:3 V:4 5       0sec  0cpu
np=3 d:1 p:1    0sec     0cpu
\end{verbatim}
The last three nef partitions each describe a K3 surface. The first one is a $T^4=T^2\times T^2$. The extra output triggered by {\tt -D} is:
\begin{verbatim}
H:4 [0] h1=2 P:0 V:2 3 5   D     0sec  0cpu
\end{verbatim}
{\tt h1=2} indicates that the Hodge number $h^{1,0}(T^4)=2\,$. Furthermore the letter {\tt D} indicates that the nef partition is a direct product.
%%%%%%%%%%%%%%%%%%%%%%%%%%%%%%%%%%%%%%%%%%%%%%%%%%%%%%%%%%%%%%%%%%%%%%%%%%%%%
\subsubsection{-P}
\label{sec-nef-P}
The option {\tt -P} also prints nef partitions corresponding to
projections. 
Consider for example a complete intersection of codimension 2 in $\mathbb{P}^3$:
\begin{verbatim}
palp$ nef.x -P
4 1 1 1 1
4 1 1 1 1 M:35 4 N:5 4  codim=2 #part=2
H:[0] P:0 V:2 3       0sec  0cpu
H:[0] P:1 V:3       0sec  0cpu
np=1 d:0 p:1    0sec     0cpu
\end{verbatim}
Compared to the output without {\tt -P} there is one additional line:
\begin{verbatim}
H:[0] P:1 V:3       0sec  0cpu
\end{verbatim}
Let $v_0,\dots,v_3$ denote the vertices of the polytope.
The nef partition {\tt P:0} is then as follows: 
\begin{equation}
0:\; V_0=\langle v_3 \rangle,\quad V_1=\langle v_0, v_1, v_2 \rangle.
\end{equation}
The part $V_0\,$ only contains the vertex $v_3$. Therefore
the equation of the corresponding divisor $D_{0,0}$ in~(\ref{eq:CICY})
reads $x_3=0$. 
The projection $\pi$ of $\Delta^*$ along $v_3$ yields a reflexive
polytope $\Delta_{v_3}^* = \langle \pi v_0,\pi v_1,\pi v_2 \rangle$.  
Thus, we are left with a hypersurface $X' = D_{0,1} \subset \mathbb{P}^2 = \mathbb{P}^3 \cap D_{0,0}$.
If there is a nef partition such that the dual nef partition in the M-lattice
has a summand with only one vertex, then {\tt DP} is displayed in the output\footnote{We thank Benjamin Nill for pointing this out to us.}. 
%%%%%%%%%%%%%%%%%%%%%%%%%%%%%%%%%%%%%%%%%%%%%%%%%%%%%%%%%%%%%%%%%%%%%%%%%%%%%
\subsubsection{-t}
\label{sec-nef-t}
The option {\tt -t} gives detailed information about the calculation
times of the Hodge numbers. The Hodge numbers of a Calabi--Yau complete
intersection are generated by the so called stringy E-function
introduced by Batyrev and Borisov in \cite{Batyrev:1995ca}. The
combinatorial construction of the E-function involves the construction
of a B-polynomial and an S-polynomial defined in
\cite{Batyrev:1995ca}. The option {\tt -t} returns the accumulated
computing times of the respective polynomials. We illustrate this
option with the example of complete intersections of codimension 4 in
$\mP^7$ (cf. Section~\ref{sec-nef-p}). 
\begin{verbatim}
palp$ nef.x -t -c4
Degrees and weights  `d1 w11 w12 ... d2 w21 w22 ...'
  or `#lines #columns' (= `PolyDim #Points' or `#Points PolyDim'):
8 1 1 1 1 1 1 1 1
8 1 1 1 1 1 1 1 1 M:6435 8 N:9 8  codim=4 #part=5
   BEGIN S-Poly     0sec  0cpu
   BEGIN B-Poly     11564sec  11558cpu
   BEGIN E-Poly     13126sec  13119cpu
H:1 65 [-128] P:0 V0:2 3  V1:4 5  V2:6 7       13126sec  13119cpu
np=1 d:0 p:4    0sec     0cpu
\end{verbatim}
This option can be useful for finding at which point in the calculation of the Hodge numbers the program crashes.

\subsubsection{-c*}
\label{sec-nef-c*}
The option {\tt -c*} where {\tt *} is a positive integer $r$ allows to
specify the length of the nef partition and hence the codimension of
the Calabi-Yau complete intersection. The default value for the
codimension is $2$. Note that the computation time can take several
hours for $r=4$ or even days for $r > 4$ and PALP may crash because
the limits such as the number of vertices etc. set in {\tt Global.h} may be exceeded, cf. Section \ref{error-handling}.  
We illustrate this option with complete intersections of codimension $3$ in $\mathbb{P}^2\times \mathbb{P}^2$:
\begin{verbatim}
palp$ nef.x -c3
Degrees and weights  `d1 w11 w12 ... d2 w21 w22 ...'
  or `#lines #columns' (= `PolyDim #Points' or `#Points PolyDim'):
3 1 1 1 0 0 0 3 0 0 0 1 1 1
3 1 1 1 0 0 0  3 0 0 0 1 1 1 M:100 9 N:7 6  codim=3 #part=7
H:[0] P:0 V0:1 3  V1:4 5       1sec  1cpu
H:[0] P:1 V0:2 3  V1:4 5       1sec  0cpu
np=1 d:1 p:5    0sec     0cpu
\end{verbatim}
Also hypersurfaces can be analyzed with {\tt nef.x}. As an example we consider the quintic hypersurface in $\mathbb{P}^4$:
\begin{verbatim}
palp$ nef.x -c1
Degrees and weights  `d1 w11 w12 ... d2 w21 w22 ...'
  or `#lines #columns' (= `PolyDim #Points' or `#Points PolyDim'):
5 1 1 1 1 1
5 1 1 1 1 1 M:126 5 N:6 5  codim=1 #part=1
H:1 101 [-200] P:0      0sec  0cpu
np=1 d:0 p:0    0sec     0cpu
\end{verbatim}
Compare this to the output of {\tt poly.x}:
\begin{verbatim}
palp$ poly.x
Degrees and weights  `d1 w11 w12 ... d2 w21 w22 ...'
  or `#lines #columns' (= `PolyDim #Points' or `#Points PolyDim'):
5 1 1 1 1 1
5 1 1 1 1 1 M:126 5 N:6 5 H:1,101 [-200]
\end{verbatim}

\subsubsection{-F*}
\label{sec-nef-F*}
The option {\tt -F*} yields information about possible toric
fibrations of the toric variety associated to the given reflexive
lattice polytope. The polytopes associated to the toric fibers are restricted to be reflexive. By considering nef partitions for the given lattice polytope this option also computes possible fibrations of the corresponding complete intersection Calabi-Yau manifolds by lower-dimensional complete intersection Calabi-Yau manifolds. For more details see \cite{Kreuzer:2000qv,Klemm:2004km}.
In practice one should always use the option {\tt -F*} in conjunction
with either {\tt -Lv} or {\tt -Lp}. Here {\tt *} is a non-negative
integer $s$ that specifies the maximal codimension $s$ of the fiber polytope. The default value for $s$ is $2$. Note that this codimension does not need to coincide with the codimension of the corresponding complete intersection Calabi-Yau fiber.
Besides the standard output and the output from the options {\tt -Lv} or {\tt -Lp}, the full information about fibration structures is listed below a second dashed line. The output takes the following form:
\begin{verbatim}
----------------------------------------------- #fibrations=#
    _    _    v    v    ...    p    p    p    v  cd=#  m: # # n: # #

    .    .    .    .    ...    .    .    .    .    .      . .    . .
    .    .    .    .    ...    .    .    .    .    .      . .    . .
    .    .    .    .    ...    .    .    .    .    .      . .    . .

    v    p    _    v    ...    v    _    _    p  cd=#  m: # # n: # #
\end{verbatim}
The number {\tt \#} in {\tt \#fibrations=\#} specifies the number of
fibrations by reflexive polytopes up to symmetry that have been
found. Then each of the subsequent lines corresponds to one of these
fibrations. The points of the given polytope are labeled by either
{\tt v}, {\tt p} or {\tt \_}. This label indicates whether the
corresponding point is a vertex ({\tt v}), a non-vertex point ({\tt
  p}) or not a point at all ({\tt \_}) of the fiber polytope. 
The latter correspond to the directions of the toric base.
The non-negative integer {\tt \#} in {\tt cd=\#} specifies
the codimension of the fiber polytope. The two positive integers {\tt
  \# \#} after {\tt m}: specify the numbers of points and vertices of the dual of the fiber polytope, respectively. The two positive integers {\tt
  \# \#} after {\tt n}: specify the numbers of points and vertices of the fiber polytope, respectively.

We illustrate this option with a complete intersection of codimension $2$ with several fibrations. In order to find all fibrations the argument of {\tt -F} must be set to $3$. This is an example where the interpretation of the fibration information depends on the choice of the nef partition.
\begin{verbatim}
palp$ echo "12 4 2 2 2 1 1 0  8 4 0 0 0 1 1 2" | nef.x -f -Lp -F3
12 4 2 2 2 1 1 0  8 4 0 0 0 1 1 2 M:371 12 N:10 7  codim=2 #part=5
5 10  Points of Poly in N-Lattice:
    0    0    0    1    0   -1    0    0    0    0
    0    0    1    0    0   -1    0    0    0    0
   -1    4    0    0    0    0    0    1    2    0
    0   -1    0    0    1    0    0    0    0    0
   -1    2    0    0    0    1    1    1    1    0
--------------------------------------------------
    4    1    2    2    1    2    0    0    0  d=12  codim=0
    4    1    0    0    1    0    2    0    0  d=8  codim=2
    2    0    1    1    0    1    0    0    1  d=6  codim=1
    2    0    0    0    0    0    1    0    1  d=4  codim=3
    1    0    0    0    0    0    0    1    0  d=2  codim=4
--------------------------------------------- #fibrations=3
    v    v    _    _    v    _    v    p    p  cd=2  m: 35  4 n: 7 4
    v    _    v    v    _    v    v    p    v  cd=1  m:117  9 n: 8 6
    v    _    _    _    _    _    v    p    v  cd=3  m:  9  3 n: 5 3
H:4 58 [-108] P:1 V:0 2   (6 6) (4 4) (3 3) (2 2) (1 1)   1sec  0cpu
H:3 65 [-124] P:2 V:0 2 3   (8 4) (4 4) (4 2) (2 2) (1 1) 1sec  0cpu
H:3 83 [-160] P:3 V:3 5   (4 8) (0 8) (2 4) (0 4) (0 2)   1sec  1cpu
np=3 d:0 p:2    0sec     0cpu
\end{verbatim}
There are three fibrations.  The fiber polytope of the second
fibration is of codimension $1$, hence has dimension $5-1=4$. As
usual, we label the vertices and points by
$v_0,\ldots,v_6,p_7,p_8,p_9$. The vertices labeled with {\tt \_} are
$v_1$ and $v_4$, which are all in $V_1$ for all the three nef
partitions. Since we are considering a complete intersection of
codimension $2$, the corresponding Calabi--Yau threefold admits a
fibration by K3 surfaces since the fiber has dimension $4-2=2$. The
linear relation of codimension $1$ and degree $6$ does not involve
$v_1$ and $v_4$, hence it describes the fiber polytope. The degrees of
the nef partitions with respect to this linear relation are given in
the third parentheses in the lines containing the information of the
nef partitions. Hence, the K3 fibers are $\mathbb{P}(2,1,1,1,1)[3,3]$,
$\mathbb{P}(2,1,1,1,1)[4,2]$, and $\mathbb{P}(2,1,1,1,1)[2,4]$,
respectively. Note that the second fibration is an instance of the
situation that a non-vertex point of the polytope becomes a vertex of
the fiber polytope. Here, this is the point $p_8$. 

The fiber polytope of the first fibration is of codimension $2$, hence
has dimension $5-2=3$. Naively, one would expect that the
corresponding Calabi--Yau threefolds admit elliptic fibrations. This
is indeed true for the first two nef partitions where both $V_0$ and
$V_1$ contain vertices belonging to the fiber polytope. 
Repeating the steps of the second fibration above in this case
yields the complete intersection $\mathbb{P}(4,1,1,2)[4,4]$ for both
nef partitions. After discarding the trivial projection to the first
coordinate, they become the hypersurfaces $\mathbb{P}(1,1,2)[4]$. 

For the third nef partition, however, the vertices and points of the
fiber polytope only lie in the part $V_1$ of the nef partition. Hence,
the part $V_0$ reduces the dimension of the base. 
The fiber of the corresponding Calabi-Yau threefold is only of
codimension $1$ in the $3$-dimensional toric fiber, i.e. it is a K3
surface. In fact, the linear relation of codimension $2$ and degree
$8$ involves all points of $V_1$, hence it describes the fiber
polytope. The degrees of the third nef partition with respect to this
linear relation are the second parentheses in the line with {\tt
  P:3}. Hence, the K3 fiber is $\mathbb{P}(4,1,1,2)[8]$. This
phenomenon is further described in \cite{Klemm:2004km}.

Finally, the fiber polytope of the third fibration is of codimension $3$, and hence has dimension $5-3=2$. Naively, one would expect that the corresponding Calabi-Yau threefolds do not admit any fibrations since the codimension is also $2$ and hence the fibers would be points. This is indeed the case for the first two nef partitions. For the third nef partition, the fiber polytope consists of the points $v_0, v_6, p_7$, and $p_8$, all of which lie in $V_1$. Hence, the fiber of the corresponding Calabi-Yau threefold is only of codimension $1$ in the $2$-dimensional toric fiber, i.e. it is an elliptic curve.  The degrees of the third nef partition with respect to the linear relation of codimension $3$ are the fourth parentheses in the line with {\tt P:3}. Hence, the elliptic curve is $\mathbb{P}(2,1,1)[4]$.
%%%%%%%%%%%%%%%%%%%%%%%%%%%%%%%%%%%%%%%%%%%%%%%%%%%%%%%%%%%%%%%%%%%%%%%%%%%%%
\subsubsection{-y}
\label{sec-nef-y}
Depending on the input the option {\tt -y} returns the CWS or the vertices of the M-lattice polytope if there is at least one nef partition. In order to trigger the output this nef partition may also be a projection. If there is no nef partition there is no output.
Depending on the input the following output is given:
\begin{itemize}
\item if there is a nef partition:
\begin{itemize}
\item If the input is a CWS, the CWS is returned along with the polytope data.
\item If the input is a polytope in the M-lattice or N-lattice
  (cf. option {\tt -N} in Section~\ref{sec-nef-N}) the M-lattice polytope is returned.
\end{itemize}
\item if there is no nef partition
\begin{itemize}
\item If the input is a CWS, the CWS is returned without further information about the polytope.
\item If the input is a polytope there is no output.
\end{itemize}
\end{itemize}
As an example consider the codimension 2 complete intersection in
$\mathbb{P}^3$ from Section~\ref{sec-nef-N}.
If we enter the N-lattice polytope we get the following output:
\begin{verbatim}
palp$ nef.x -y -N
Degrees and weights  `d1 w11 w12 ... d2 w21 w22 ...'
  or `#lines #columns' (= `PolyDim #Points' or `#Points PolyDim'):
3 4
Type the 12 coordinates as dim=3 lines with #pts=4 columns:
-1    0    0    1
-1    0    1    0
-1    1    0    0    
3 4 Vertices of Poly in M-lattice:  M:35 4 N:5 4  codim=2 #part=2
   -1   -1   -1    3
   -1   -1    3   -1
   -1    3   -1   -1
\end{verbatim}

\subsubsection{-S}
\label{sec-nef-S}
The option {\tt -S} gives information about the number of points in
the reflexive Gorenstein cone and its dual (cf. options {\tt -g*} and {\tt -d*}
discussed in Sections \ref{sec-nef-g*} and \ref{sec-nef-d*}) for each
nef partition which is not a direct product or a projection. It 
displays the numbers $\ell$ of lattice points and $\ell^*$ of interior lattice 
points in degrees $k \le (\widetilde d + 1)/2$, where $\widetilde d$ is the 
dimension of the 
Gorenstein cone $C$, and the analogous data for the dual cone $\widecheck C$.
These data enter the calculation of the (stringy) Hodge numbers  
via the S-polynomial (hence the name {\tt -S})
as described in Section~\ref{sec:nef-part}.
The output takes the following form. After the first line of the
standard output, there is a part referring to the polytope
$\Delta(\widecheck C)$:
\begin{verbatim}
#points in largest cone:
layer:  1 #p:        l1 #ip:        l*1
  ...   . ...        .  ...         .
layer:  . #p:        .  #ip:        .  
  ...   . ...        .  ...         .
layer:  k #p:        lk #ip:        l*k
\end{verbatim}
where {\tt l1}~$ = \ell(\Delta(\widecheck C) ), \ldots$,
{\tt lk}~$ =\ell(k\Delta(\widecheck C) ),$ {\tt l*1}~$
= \ell^*(\Delta(\widecheck C) ),\ldots$, {\tt l*k}~$
=\ell^*(k\Delta(\widecheck C) )$.
Subsequently there is a second part referring to the polytope $\Delta (C)$.
\begin{verbatim}
#points in largest cone:
layer:  1 #p:        l1 #ip:        l*1
  ...   . ...        .  ...         .
layer:  . #p:        .  #ip:        .  
  ...   . ...        .  ...         .
layer:  k #p:        lk #ip:        l*k
\end{verbatim}
where {\tt l1}~$ = \ell(\Delta(C) ), \ldots$,
{\tt lk}~$ =\ell(\Delta(C) ),$ {\tt l*1}~$
= \ell^*(\Delta(C) ),\ldots$, {\tt l*k}~$
=\ell^*(k\Delta(C) )$. Then the rest of the standard
output concerning the nef partitions follows.

The following example illustrates this option. 
We consider a complete intersection of codimension $2$ in $\mathbb{P}^4$:
\begin{verbatim}
palp$ nef.x -S
Degrees and weights  `d1 w11 w12 ... d2 w21 w22 ...'
  or `#lines #columns' (= `PolyDim #Points' or `#Points PolyDim'):
4 1 1 1 1
4 1 1 1 1 M:35 4 N:5 4  codim=2 #part=2


#points in largest cone:
layer:  1 #p:        6 #ip:        0
layer:  2 #p:       21 #ip:        1
layer:  3 #p:       56 #ip:        6


#points in largest cone:
layer:  1 #p:       20 #ip:        0
layer:  2 #p:      105 #ip:        1
layer:  3 #p:      336 #ip:       20
H:[0] P:0 V:2 3       0sec  0cpu
np=1 d:0 p:1    0sec     0cpu
\end{verbatim}
One of the two nef partitions is a projection and is not analyzed. The
output for the remaining nef partition has two blocks: The first block
counts the numbers of points (after {\tt \#p:}) and points in the
relative interior (after {\tt \#ip:}) of the Gorenstein cone
$\widecheck C \subset \widetilde N$ at degrees $k = 1,2,3$. Hence
\begin{gather*}
\ell(\Delta(\widecheck C) ) =6,\quad \ell(2\Delta(\widecheck C) ) =21,\quad \ell(3\Delta(\widecheck C) ) =56,\\
\ell^*(\Delta(\widecheck C) ) =0,\quad \ell^*(2\Delta(\widecheck C) ) =1,
\quad \ell^*(3\Delta(\widecheck C) ) =6.
\end{gather*}
One can check that the number of points at degree $k=1\,$ indeed
coincides with the number of points in the output of the option {\tt
  -g2}.

The second block gives the same information for the dual Gorenstein
cone $C \subset \widetilde M$. Hence
\begin{gather*}
\ell(\Delta(C) ) =20,\quad \ell(2\Delta(C) ) =105,\quad \ell(3\Delta(C) ) =336,\\
\ell^*(\Delta(C) ) =0,\quad \ell^*(2\Delta(C) ) =1,\quad \ell^*(2\Delta(C) ) =20.
\end{gather*}
The output of the option {\tt -d2} coincides with the number of points at degree $k=1\,$.

\subsubsection{-T}
\label{sec-nef-T}
The option {\tt -T} turns on an explicit check of the relation (\ref{eq:Serre})
relating the S-- and T--polynomials.
Normally the program actually uses that relation to avoid point counting
beyond degree $(\widetilde d + 1)/2$, but with {\tt -T} the counting goes
up to degree $\widetilde d$ and an error message is given if (\ref{eq:Serre})
is violated.
This can be useful if one suspects that the program gives wrong Hodge numbers,
for example because of numerical overflows.
If nothing goes wrong, the only effect is a significantly increased
computation time.
The best way to illustrate this option is by combining it with {\tt -S}.
We consider the same example as in the previous subsection.
\begin{verbatim}
palp$ nef.x -S -T
Degrees and weights  `d1 w11 w12 ... d2 w21 w22 ...'
  or `#lines #colums' (= `PolyDim #Points' or `#Points PolyDim'):
4 1 1 1 1
4 1 1 1 1 M:35 4 N:5 4  codim=2 #part=2


#points in largest cone:
layer:  1 #p:        6 #ip:        0
layer:  2 #p:       21 #ip:        1
layer:  3 #p:       56 #ip:        6
layer:  4 #p:      125 #ip:       21
layer:  5 #p:      246 #ip:       56


#points in largest cone:
layer:  1 #p:       20 #ip:        0
layer:  2 #p:      105 #ip:        1
layer:  3 #p:      336 #ip:       20
layer:  4 #p:      825 #ip:      105
layer:  5 #p:     1716 #ip:      336
H:[0] P:0 V:2 3       0sec  0cpu
np=1 d:0 p:1    0sec     0cpu
\end{verbatim}
Note how now the point counting proceeds up to degree 5.
With these data we can compute the Ehrhart polynomial 
\[
S_{\Delta(\widecheck C)}(t) = (1 - t)^5 (1 + 6t + 21 t^2 + 56 t^3 + 125 t^4 + 246
t^5 + \dots) \]
Since it has degree at most $\widetilde d = 5$, we find 
\[ S_{\Delta(\widecheck C)} = 1 + t + t^2 + t^3.\]
Similarly
\[
T_{\Delta(\widecheck C)}(t) = (1 - t)^5 (t^2 + 6t^3 + 21 t^4 + 56 t^5 + \dots) 
= t^2 + t^3 + t^4 + t^5,\]
and it is clear that (\ref{eq:Serre}) is satisfied.
A similar check can be performed for $C \cap \widetilde
M$ with the data from the second block.
%%%%%%%%%%%%%%%%%%%%%%%%%%%%%%%%%%%%%%%%%%%%%%%%%%%%%%%%%%%%%%%%%%%%%%%%%%%%%
\subsubsection{-s}
\label{sec-nef-s} 
The option {\tt -s} includes all nef partitions in the output, not
just one representative for each class of nef partitions that are
equivalent under symmetries of the CWS. Note that this option does not
print all possible nef partitions as those corresponding to
projections (cf. option {\tt -P} in Section~\ref{sec-nef-P}) or direct
products (cf. option {\tt -D} in Section~\ref{sec-nef-D}) are omitted.
The example we consider is a complete intersection of codimension $2$ in $\mathbb{P}^2\times\mathbb{P}^2$. We add the option {\tt -Lv} in order to print the vertices and the CWS.
\begin{verbatim}
palp$ nef.x -s -Lv
Degrees and weights  `d1 w11 w12 ... d2 w21 w22 ...'
  or `#lines #columns' (= `PolyDim #Points' or `#Points PolyDim'):
3 1 1 1 0 0 0 3 0 0 0 1 1 1
3 1 1 1 0 0 0  3 0 0 0 1 1 1 M:100 9 N:7 6  codim=2 #part=31
4 6 Vertices in N-lattice:
    0    0    0    1    0   -1
    0    0    1    0    0   -1
   -1    0    0    0    1    0
   -1    1    0    0    0    0
------------------------------
    1    1    0    0    1    0  d=3  codim=2
    0    0    1    1    0    1  d=3  codim=2
H:20 [24] P:2 V:4 5   (1 2) (1 2)     0sec  0cpu
H:20 [24] P:4 V:0 5   (1 2) (1 2)     0sec  0cpu
H:20 [24] P:5 V:0 4   (2 1) (0 3)     0sec  0cpu
H:20 [24] P:6 V:0 4 5   (2 1) (1 2)     0sec  0cpu
H:20 [24] P:8 V:1 5   (1 2) (1 2)     1sec  0cpu
H:20 [24] P:9 V:1 4   (2 1) (0 3)     0sec  0cpu
H:20 [24] P:10 V:1 4 5   (2 1) (1 2)     0sec  0cpu
[further Hodge data and nef partitions]
\end{verbatim}
Note that the CWS is symmetric under permutations of the vertices
labeled by $0,1,4$ and those labeled by $2,3,5$. Furthermore there
only exist three pairs of degrees of the complete intersection (up to
exchange within a pair):
$\{(1,2),(1,2)\},\{(0,3),(2,1)\},\{(1,2),(2,1)\}$. Therefore we
conclude that there are only three inequivalent nef partitions. This
is indeed confirmed by calling {\tt nef.x} without the option {\tt -s}.

\subsubsection{-n}
\label{sec-nef-n}
The option {\tt -n} prints the points of the polytope in the N-lattice
only if there is at least one nef partition which does not correspond
to a projection or a direct product. In addition, the first line of
the standard output is printed while the other lines are suppressed. 
As an example we consider a codimension 2 complete intersection in $\mathbb{P}^3$
\begin{verbatim}
palp$ nef.x -n
Degrees and weights  `d1 w11 w12 ... d2 w21 w22 ...'
  or `#lines #columns' (= `PolyDim #Points' or `#Points PolyDim'):
4 1 1 1 1
4 1 1 1 1 M:35 4 N:5 4  codim=2 #part=2
3 5  Points of Poly in N-Lattice:
   -1    0    0    1    0
   -1    0    1    0    0
   -1    1    0    0    0
\end{verbatim}
%%%%%%%%%%%%%%%%%%%%%%%%%%%%%%%%%%%%%%%%%%%%%%%%%%%%%%%%%%%%%%%%%%%%%%%%%%%%%
\subsubsection{-v}
\label{sec-nef-v}
The option {\tt -v} prints the size of the matrix of vertices, the number of points and the vertices of the polytope that has been entered (M-lattice or N-lattice, depending on the input). If the input is the CWS the M-lattice polytope is analyzed. The output is printed in a single line with the character {\tt E} as separator. Furthermore one can limit the output to polytopes whose number of points is constrained by a lower and an upper bound:
\begin{itemize}
\item  {\tt -v -u\#}, where {\tt \#} is an integer $\geq0$, only gives output if the polytope has at most \# points. The default value is the parameter {\tt POINT\_Nmax} which fixes the maximal number of points of a polytope at compilation.
\item {\tt -v -l\#}, where {\tt \#} is an integer $\geq0$, only gives output if the polytope has at least \# points. The default value is $0$.
\end{itemize}
After closing the program a summary is printed. It contains
information on the number of the examined polytopes which satisfy the
bounds and the number of polytopes with {\tt \#} of points. 

As an example we consider complete intersections of codimension $2$ in $\mathbb{P}^3$ and $\mathbb{P}^2\times\mathbb{P}^2$ with the weight matrices as input and without bounds.
\begin{verbatim}
palp$ nef.x -v
Degrees and weights  `d1 w11 w12 ... d2 w21 w22 ...'
  or `#lines #columns' (= `PolyDim #Points' or `#Points PolyDim'):
4 1 1 1 1
3 4 P:35 E   -1    3   -1   -1E   -1   -1    3   -1E   -1   -1 ...
Degrees and weights  `d1 w11 w12 ... d2 w21 w22 ...'
  or `#lines #columns' (= `PolyDim #Points' or `#Points PolyDim'):
3 1 1 1 0 0 0 3 0 0 0 1 1 1
4 9 P:100 E   -1    2   -1   -1    2   -1   -1    2   -1E
   -1   -1    2   -1   -1    2   -1   -1    2E
   -1   -1   -1    2    2    2   -1   -1   -1E
   -1   -1   -1   -1   -1   -1    2    2    2
Degrees and weights  `d1 w11 w12 ... d2 w21 w22 ...'
  or `#lines #columns' (= `PolyDim #Points' or `#Points PolyDim'):

2  of  2

  35#    1
 100#    1
\end{verbatim}
Since we have entered a CWS the M-lattice polytope is analyzed. Let us discuss the first line of output:
\begin{verbatim}
3 4 P:35 E   -1    3   -1   -1E   -1   -1    3   -1E   -1   -1 ...
\end{verbatim}
The first two numbers indicate the number of rows and columns of the matrix of vertices in the M-lattice polytope. {\tt P:35} indicates that the M-lattice polytope has $35$ points. The vertices of the M-lattice polytope are then written in one line with the separator {\tt E}. The output of the second example is analogous. After we quit PALP by hitting enter without input the following output is given:
\begin{verbatim}
2  of  2

  35#    1
 100#    1
\end{verbatim}
This means that $2$ out of the $2$ polytopes analyzed satisfy the
bounds and that there is one polytope with $35$ points and one with
$100$.

Next, we consider the same example as above but with the upper bound for the number of points set to $50$:
\begin{verbatim}
palp$ nef.x -v -u50
Degrees and weights  `d1 w11 w12 ... d2 w21 w22 ...'
  or `#lines #columns' (= `PolyDim #Points' or `#Points PolyDim'):
4 1 1 1 1
3 4 P:35 E   -1    3   -1   -1E   -1   -1    3   -1E   -1   -1 ...
Degrees and weights  `d1 w11 w12 ... d2 w21 w22 ...'
  or `#lines #columns' (= `PolyDim #Points' or `#Points PolyDim'):
3 1 1 1 0 0 0 3 0 0 0 1 1 1
Degrees and weights  `d1 w11 w12 ... d2 w21 w22 ...'
  or `#lines #columns' (= `PolyDim #Points' or `#Points PolyDim'):

1  of  2

  35#    1
\end{verbatim}
Now the second polytope exceeds the upper bound for the points as it has $100$ points (cf. previous example). There is no output for the second polytope and the summary indicates that only one of the two polytopes analyzed satisfies the bounds.

\subsubsection{-m}
\label{sec-nef-m}
The option {\tt -m} returns a nef partition of length $2$
resulting from a partition $d=d_1+d_2$ of the degree of a weight system.
More precisely, the input data is a single weight 
system $w$ and two positive integers $d_1, d_2$ such that $\sum_i w_i = d_1 +
d_2$. The input format is {\tt \# d=\# \#}, where the first {\tt \#} is
the usual CWS, while the {\tt \#} after {\tt d=} refer to $d_1$ and
$d_2$, respectively.
As always, $w$ specifies $\Delta^{(d)} \subset M$ 
as the Newton polytope of degree $d$.
Furthermore, the degrees $d_1, d_2$ specify Newton polytopes
$\Delta^{(d_1)}, \Delta^{(d_2)}$ from which one obtains the Minkowski sum
$\Delta^{(d_1,d_2)}=\Delta^{(d_1)}+\Delta^{(d_2)} \subseteq \Delta^{(d)}$. 
If $\Delta_1=\Delta^{(d_1)}, \Delta_2=\Delta^{(d_2)}$
define a nef partition $(\nabla_1, \nabla_2)$ of 
the vertices of $(\Delta^{(d_1,d_2)})^*$, then the data of
this nef partition are given in the standard output.

The following example taken from \cite{Klemm:2004km} illustrates this option.
We consider the weighted projective space 
$\mP(1,1,1,1,4,6)$ specified by the weight vector $14\; 1\; 1\;
1\; 1\; 4\; 6$ of degree $d=14$. The polytope $\Delta=\Delta^{(14)}$ is the
Newton polytope of degree 14 monomials in this space. 
We first analyze the toric variety determined by $\Delta^{(14)}$:
\begin{verbatim}
palp$ nef.x -Lv
Degrees and weights  `d1 w11 w12 ... d2 w21 w22 ...'
  or `#lines #colums' (= `PolyDim #Points' or `#Points PolyDim'):
14 1 1 1 1 4 6
14 1 1 1 1 4 6 M:1271 13 N:10 8  codim=2 #part=2
5 8 Vertices in N-lattice:
    0   -1    0    0    0    1    0    0
    0   -1    1    0    0    0    0    0
    0   -1    0    1    0    0    0    0
    0   -4    0    0    1    0   -1   -1
    1   -6    0    0    0    0   -1   -2
----------------------------------------
    6    1    1    1    4    1    0    0  d=14  codim=0
    1    0    0    0    1    0    1    0  d=3  codim=3
    2    0    0    0    1    0    0    1  d=4  codim=3
H:1 149 [-296] P:1 V:3 4 5 7  8   (6 8) (1 2) (2 2)     2sec  1cpu
np=1 d:0 p:1    2sec     1cpu
\end{verbatim}
So $\Delta^{(14)}$ has 1271 lattice points and 13 vertices, and
$(\Delta^{(14)})^*$ is the convex hull of the eight vertices shown in the 
output.
By considering the weight systems below the dashed line one sees that 
$\mP_{(\Delta^{(14)})^*}$ is the blowup of $\mathbb{P}(1,1,1,1,4,6)$ along the 
divisors corresponding to the last two vertices of $(\Delta^{(14)})^*$.
Now we want to use the option {\tt -m} to see whether the partition 
$14 = 2 + 12$ determines a nef partition via the Minkowski sum 
$\Delta^{(2,12)}=\Delta^{(2)} + \Delta^{(12)}$. 
\begin{verbatim}
palp$ nef.x -Lv -m
type degrees and weights [d  w1 w2 ... wk d=d_1 d_2]: 
14 1 1 1 1 4 6  d=2 12
14 1 1 1 1 4 6 d=2 12 M:1270 12 N:11 7  codim=2 #part=2
5 7 Vertices in N-lattice:
    0   -1    0    0    0    1    0
    0   -1    1    0    0    0    0
    0   -1    0    1    0    0    0
    0   -4    0    0    1    0   -2
    1   -6    0    0    0    0   -3
-----------------------------------
    6    1    1    1    4    1    0  d=14  codim=0
    3    0    0    0    2    0    1  d=6  codim=3
 d=12 2H:3 243 [-480] P:0 V:3 5   (2 12) (0 6)     7sec  6cpu
np=1 d:0 p:1    0sec     0cpu
\end{verbatim}
The output indeed yields such a nef partition. 
Since not every monomial of degree 14 is a product
of monomials of degree 2 and 12, the polytope $\Delta^{(2,12)}$ is only a proper
subpolytope of $\Delta^{(14)}$. 
Consequently $\mP_{(\Delta^{(2,12)})^*}$ is obtained from $\mP_{(\Delta^{(14)})^*}$ 
by a further blowup along the vertex $(0,0,0,-2,-3)^T$.

By using the option {\tt -m} in the same way one can find that 
$\Delta^{(6,8)}=\Delta^{(14)}$ and that $14=3+11$ does not give rise
to a nef partition.

\subsubsection{-R}
\label{sec-nef-R}
The option {\tt -R} prints the vertices of the input polytope if it is not reflexive. To illustrate this we enter the CWS of a polytope which is not reflexive:
\begin{verbatim}
palp$ nef.x -R
Degrees and weights  `d1 w11 w12 ... d2 w21 w22 ...'
  or `#lines #columns' (= `PolyDim #Points' or `#Points PolyDim'):
6 3 2 1 0 0 5 0 0 1 1 3
3 7  Vertices of input polytope:
   -1    1    0    0    1    0   -1
    0    1   -1    1    4    4    1
   -1    0    0    0   -1   -1   -1
\end{verbatim}
The same output is given if we enter the polytope itself. Without the option {\tt -R} there is no output if the polytope is not reflexive.

\subsubsection{-V}
\label{sec-nef-V}
The option {\tt -V} prints the vertices of the polytope in the
N-lattice together with the standard output. In contrast to the option
{\tt -Lv} (cf. Section~\ref{sec-nef-Lv}) the information about the linear relations is not
given. Furthermore, in the lines containing the nef partitions the
additional information about the degrees is omitted. The option {\tt
  -V} also works for non-reflexive polytopes. As an example we consider a complete intersection of codimension 2 in $\mathbb{P}^3$:
\begin{verbatim}
palp$ nef.x -V
Degrees and weights  `d1 w11 w12 ... d2 w21 w22 ...'
  or `#lines #columns' (= `PolyDim #Points' or `#Points PolyDim'):
4 1 1 1 1
4 1 1 1 1 M:35 4 N:5 4  codim=2 #part=2
3 4  Vertices of P:
   -1    0    0    1
   -1    0    1    0
   -1    1    0    0
H:[0] P:0 V:2 3       0sec  0cpu
np=1 d:0 p:1    0sec     0cpu
\end{verbatim}
We can also enter the M-lattice polytope to get the same result.
If the polytope is non-reflexive the output is the same as for the
option {\tt -R} (cf. Section~\ref{sec-nef-R}).

\subsubsection{-Q}
\label{sec-nef-Q}
The option {\tt -Q} prints the information about the nef partitions
and the Hodge numbers only if the corresponding complete intersection
is a direct product (cf. option {\tt -D} in Section~\ref{sec-nef-D})
up to lattice quotients. If none of the nef partitions is a direct
product only the numbers of points and vertices in the M- and
N-lattice, together with the codimension and the number of nef
partitions is given.

Consider the complete intersection of codimension $2$ in $\mathbb{P}^2\times\mathbb{P}^2$. As one can check using the option {\tt -D} one of the nef partitions corresponds to a direct product:
\begin{verbatim}
palp$ nef.x -Q
Degrees and weights  `d1 w11 w12 ... d2 w21 w22 ...'
  or `#lines #columns' (= `PolyDim #Points' or `#Points PolyDim'):
3 1 1 1 0 0 0 3 0 0 0 1 1 1
3 1 1 1 0 0 0  3 0 0 0 1 1 1 M:100 9 N:7 6  codim=2 #part=5
H:4 [0] h1=2 P:0 V:2 3 5   D     0sec  0cpu
np=4 d:1 p:0    0sec     0cpu
\end{verbatim}
The N-lattice polytope of $\mathbb{P}^3$ has no nef partition
corresponding to a direct product.  
Then the output looks as follows:
\begin{verbatim}
palp$ nef.x -Q
Degrees and weights  `d1 w11 w12 ... d2 w21 w22 ...'
  or `#lines #columns' (= `PolyDim #Points' or `#Points PolyDim'):
4 1 1 1 1
4 1 1 1 1 M:35 4 N:5 4  codim=2 #part=2
np=2 d:0 p:0    0sec     0cpu
\end{verbatim}

\subsubsection{-g*}
\label{sec-nef-g*}
The option {\tt -g*}, where {\tt *} is an integer $m=0,1,2$, returns
the points of the supports $\Delta(\widecheck C)$ of the
Gorenstein cones $\widecheck C \subset \widetilde N_\mR$ associated to the 
nef partitions of length $r$ of the input polytope $\Delta^* \subset N_\mR$. 
For the notation on Gorenstein cones see Section~\ref{sec:nef-part}. The
default value is $m=1$. The standard output is changed as follows. 
The lines containing the information about the nef partition including the 
Hodge numbers, the parts of the nef partition etc. are suppressed. 
Instead, for each nef partition the points of $\Delta(\widecheck C)$ are
printed in the following form:
\begin{verbatim}
D n Points of PG: (nv=#)
    #    #      ...      #    #
    .    .      ...      .    . 
    .    .      ...      .    . 
    #    #      ...      #    #
\end{verbatim}
The interpretation depends on the integer $m$. 
For $m=2$ the output is the list of points $\widetilde p \in \widecheck C$
as given in~(\ref{eq:ptilde}). 
Note that since the origin $0_N$
belongs to every part of the nef partition, it appears $r$ times, each
time another of the $r$ support functions being equal to $1$. 
For $m=1$ the redundant coordinate $\phi_0(p)$ is omitted in $\widetilde p$ 
and we obtain vectors $\widetilde p'$. 
For $m=0$ all $\phi_i(p)$ are omitted and the resulting $r$-fold
occurrence of $0_N$ is reduced to just a single occurrence; information on 
the nef partition is lost and the output becomes just 
the list of lattice points of $\Delta^*$.
The values of {\tt D}, {\tt n} and the {\tt \#} columns are summarized in 
Table~\ref{tab:g},
\begin{table}[h]
\begin{center}
  {\begin{tabular}{c|c|c|c}
    $m$ & $0$ & $1$ & $2$\\
    \hline
    {\tt D}  & $d$ & $\widetilde d-1$ & $\widetilde d$\\
    {\tt n} & $n$ & $n+r-1$ & $n+r-1$\\
    {\tt \#} column & $p$ & $\widetilde p'$ & $\widetilde p$
  \end{tabular}}\caption{The meaning of the output of the options {\tt -g*}
    and {\tt -d*}}\label{tab:g}
\end{center}
 \end{table}
where $n$ is the number of lattice
points in $\Delta^*$ and $d, r, \widetilde
d$ are as in Section~\ref{sec:nef-part}. The number {\tt \#} in {\tt
  nv=\#} denotes the number of vertices of the cone $\widecheck
C$. The order of the points is first the vertices, then the non-vertex
points with the origin at the end. 

The following example illustrates this option. 
We consider complete intersections of codimension $2$ in $\mathbb{P}^2\times\mathbb{P}^1\times\mathbb{P}^2$ discussed in \cite{Braun:2007vy}.
The nef partitions for this example were discussed in Section~\ref{sec:stdout}. 
With the choice of {\tt m=2} we obtain the information about the partition in terms of the Gorenstein cone. 
Let $v_0,\dots,v_7$ denote the vertices of the polytope in the N-lattice.
\begin{verbatim}
palp$ nef.x -N -g2
Degrees and weights  `d1 w11 w12 ... d2 w21 w22 ...'
  or `#lines #colums' (= `PolyDim #Points' or `#Points PolyDim'):
5 8
Type the 40 coordinates as dim=5 lines with #pts=8 colums:
 1  0 -1  0  0  0  0  0
 0  1 -1  0  0  0  0  0
 0  0  0  1 -1  0  0  0
 0  0  0  0  0  1  0 -1
 0  0  0  0  0  0  1 -1
M:300 18 N:9 8  codim=2 #part=15
7 10 Points of PG: (nv=8)
    1    1    0    0    0    1    1    0    0    1
    0    0    1    1    1    0    0    1    1    0
    1    0   -1    0    0    0    0    0    0    0
    0    1   -1    0    0    0    0    0    0    0
    0    0    0    1   -1    0    0    0    0    0
    0    0    0    0    0    1    0   -1    0    0
    0    0    0    0    0    0    1   -1    0    0
7 10 Points of PG: (nv=8)
    1    1    0    1    0    1    0    0    0    1
    0    0    1    0    1    0    1    1    1    0
    1    0   -1    0    0    0    0    0    0    0
    0    1   -1    0    0    0    0    0    0    0
    0    0    0    1   -1    0    0    0    0    0
    0    0    0    0    0    1    0   -1    0    0
    0    0    0    0    0    0    1   -1    0    0
[output of further nef partitions]
\end{verbatim}
Let us consider the nef partition {\tt P:8} as produced  for instance by the option {\tt -N -Lp}:
\begin{verbatim}
H:19 19 [0] P:8 V:4 5 6 7   (0 3) (1 1) (3 0)     0sec  0cpu
\end{verbatim}
This example was the focus of \cite{Braun:2007vy}. The output of {\tt -g2} for this nef partition is:
\begin{verbatim}
7 10 Points of PG: (nv=8)
    1    1    1    1    0    0    0    0    0    1
    0    0    0    0    1    1    1    1    1    0
    1    0   -1    0    0    0    0    0    0    0
    0    1   -1    0    0    0    0    0    0    0
    0    0    0    1   -1    0    0    0    0    0
    0    0    0    0    0    1    0   -1    0    0
    0    0    0    0    0    0    1   -1    0    0
\end{verbatim}
Since the first four vertices $v_0, v_1, v_2\,$ and $v_3\,$ are in $V_0\,$, we have $\phi_0(v_i) = 1\,$ and $\phi_1(v_i) = 0\,$, hence the corresponding points of the Gorenstein cone take the form $(1,0,v_i)\,$ for $i=0,\dots,3$. The next four vertices $v_4, v_5, v_6\,$ and $v_7\,$ are in $V_1\,$, we have $\phi_0(v_i) = 0\,$ and $\phi_1(v_i) = 1\,$, hence the corresponding points of the Gorenstein cone take the form $(0,1,v_i)\,$ for $i=4,\dots,7$. Finally, the origin always belongs to every part of the nef partition, hence it appears as often as the codimension which here is $r=2\,$. So $p_8 = 0\,$ and $p_9 = 0\,$. Once with $\phi_0(p_8) = 0\,$ and $\phi_1(p_8) = 1\,$ and once with $\phi_0(p_9) = 1\,$ and $\phi_1(p_9) = 0\,$. 

\subsubsection{-d*}
\label{sec-nef-d*}
The option {\tt -d*}, where {\tt *} is an integer $m=0,1,2$, returns
the points of the Gorenstein cones $C \subset \widetilde M_\mR$ associated 
to the nef partitions of length $r$ of the polytope $\nabla^* \subset M_\mR$. 
For the notation on Gorenstein cones see Section~\ref{sec:nef-part}, in
particular~(\ref{eq:polyDN}) for the polytope $\nabla^*$. 
This option can be used to determine the polytope $\nabla^*$ for each of the nef partitions of the given polytope $\Delta^*$. The polytope $\nabla^*$ can then be further analyzed with {\tt poly.x}.

The integer $m$ triggers the same output format as for the option {\tt -g*} in
Section~\ref{sec-nef-g*},
The default value is $m=1$. The option {\tt -d2} automatically sets the flag {\tt -p}.

The following example illustrates this option. We consider complete
intersections of codimension $2$ in
$\mathbb{P}^2\times\mathbb{P}^1\times\mathbb{P}^2$ discussed in
\cite{Braun:2007vy}. For more details on the nef partitions see the
example in Section~\ref{sec:stdout} and Section~\ref{sec-nef-g*}. 
\begin{verbatim}
palp$ nef.x -N -d2
Degrees and weights  `d1 w11 w12 ... d2 w21 w22 ...'
  or `#lines #colums' (= `PolyDim #Points' or `#Points PolyDim'):
5 8
Type the 40 coordinates as dim=5 lines with #pts=8 colums:
 1  0 -1  0  0  0  0  0        
 0  1 -1  0  0  0  0  0
 0  0  0  1 -1  0  0  0
 0  0  0  0  0  1  0 -1
 0  0  0  0  0  0  1 -1
M:300 18 N:9 8  codim=2 #part=15
7 63 Points of dual PG: (nv=27)
   1   0   1   0   0   1   1   1 ...  
   0   1   0   1   1   0   0   0 ...   
  -1   0   1   1   0  -1  -1  -1 ...   
  -1   1  -1   0   1   1   1  -1 ...   [63-8=55 more points]
   0   1   0   1  -1   0   0   0 ...   
  -1   0  -1   0   0   1  -1  -1 ...   
  -1   1   1   1   1  -1  -1   1 ...   
[...]
\end{verbatim}
For each of the $11$ nef partitions of the input polytope $\Delta^*$
we get a $7$-dimensional dual Gorenstein cone $C$ from which the
points of the polytope
$\nabla^*$ can be read off by omitting the first two entries of each
column, cf~(\ref{eq:ptilde}). 
The numbers of points and vertices of $\nabla^*$ depend on which of the
nef partitions is considered.
The nef partition of interest in \cite{Braun:2007vy} was {\tt P:8}. The
corresponding output of {\tt -d2} is
\begin{verbatim}
7 40 Points of dual PG: (nv=12)
   0   0   1   1   1   0   0   0 ...       
   1   1   0   0   0   1   1   1 ...      
   0   0  -1   2  -1   0   0   0 ...     
   0   0   2  -1  -1   0   0   0 ...   [40-8=32 more points]   
   0   1   0   0   0   1   1   0 ...
  -1  -1   0   0   0   2  -1  -1 ...   
  -1   2   0   0   0  -1  -1   2 ...   
\end{verbatim}
We see that the polytope $\nabla^*$ 
has $39$ points 
(the interior
point appears twice) and $12$ vertices. Let $e_1,\dots,e_5$ be the
standard basis of $\mathbb{R}^5$. Let
$\check{v}_0,\dots,\check{v}_{11}$ denote the vertices of the polytope
$\nabla^*$.
with
\begin{align*}
\check{v}_0 &= -e_4-e_5, & \check{v}_1 &= e_3-e_4+2e_5, & \check{v}_2
&= -e_1+2e_2, \\\check{v}_3 &= 2e_1 -e_2, &
\check{v}_4 &= -e_1-e_2, & \check{v}_5 &= e_3+2e_4-e_5, \\ \check{v}_6
&= e_3-e_4-e_5, & \check{v}_7 &= -e_4+2e_5, &
\check{v}_8 &= -e_1+2e_2-e_3, \\\check{v}_9 &= 2e_1 - e_2 -e_3,
&\check{v}_{10} &=-e_1-e_2-e_3, & \check{v}_{11} &= 2e_4-e_5.
\end{align*}
From the first two rows of the above output we can read off the nef partition $\widecheck V = \widecheck V_0 \cup \widecheck V_1$ of $\nabla^*$:
\begin{equation*}
\widecheck{V}_0 = \langle \check{v}_2, \check{v}_3, \check{v}_4, \check{v}_8, \check{v}_9, \check{v}_{10} \rangle, \quad \widecheck{V}_1 = \langle \check{v}_0, \check{v}_1, \check{v}_5, \check{v}_6, \check{v}_7, \check{v}_{11} \rangle.
\end{equation*}
We can check this by feeding the vertices back into {\tt nef.x} with the options {\tt -N} and {\tt -Lv}.
\begin{verbatim}
palp$ nef.x -N -Lv
Degrees and weights  `d1 w11 w12 ... d2 w21 w22 ...'
 or `#lines #colums' (= `PolyDim #Points' or `#Points PolyDim'):
5 12
Type the 60 coordinates as dim=5 lines with #pts=12 colums:
  0   0  -1   2  -1   0   0   0  -1   2  -1   0
  0   0   2  -1  -1   0   0   0   2  -1  -1   0
  0   1   0   0   0   1   1   0  -1  -1  -1   0
 -1  -1   0   0   0   2  -1  -1   0   0   0   2
 -1   2   0   0   0  -1  -1   2   0   0   0  -1
M:24 15 N:39 12  codim=2 #part=2
5 12 Vertices in N-lattice:
   0    0   -1    2   -1    0    0    0   -1    2   -1    0
   0    0    2   -1   -1    0    0    0    2   -1   -1    0
   0    1    0    0    0    1    1    0   -1   -1   -1    0
  -1   -1    0    0    0    2   -1   -1    0    0    0    2
  -1    2    0    0    0   -1   -1    2    0    0    0   -1
------------------------------------------------------------
[linear relations]
H:19 19 [0] P:0 V:1 2 3 4 5 6  13 15 17 ...   [degrees]
H:19 19 [0] P:1 V:2 3 4 8 9 10  16 17 18 ...   [degrees]
np=2 d:0 p:0    0sec     0cpu
\end{verbatim}
We see that the nef partition {\tt P:1} agrees with $\widecheck V = \widecheck V_0 \cup \widecheck V_1$. 

\subsubsection{-G}
\label{-sec-nef-G}
The option {\tt -G} works directly with Gorenstein cones
which need not correspond to nef partitions. 
The input polytope is interpreted as the support polytope $\Delta(C)$ 
of a reflexive Gorenstein cone $C$, cf. Section~\ref{sec:nef-part}. 
The index $r$ of the cone is 2 by default and can be set to different
values with the {\tt -c} option, cf. Section~\ref{sec-nef-c*}.
The standard output contains information on the support polytopes of the
cone and the dual cone and the string--theoretic Hodge 
numbers $h^{ij}$, 
$0\le i, j\le \dim C - 2r$, see~\cite{Batyrev:1995ca}.
If the 
input does not correspond to a reflexive Gorenstein cone of index $r$, 
no Hodge numbers and no N lattice data can be computed;
as usual, the number of facets is displayed instead.
If the input corresponds to a reflexive Gorenstein cone of an index different from $r$, this is treated like a non-reflexive case but with a warning message.
\begin{verbatim}
palp$ nef.x -G
Degrees and weights  `d1 w11 w12 ... d2 w21 w22 ...'
  or `#lines #columns' (= `PolyDim #Points' or `#Points PolyDim'):
4 2
Type the 8 coordinates as #pts=4 lines with dim=2 columns:
0 0
0 1
1 0
1 1
M:4 4 N:4 4 H:[0] h0=0
Degrees and weights  `d1 w11 w12 ... d2 w21 w22 ...'
  or `#lines #columns' (= `PolyDim #Points' or `#Points PolyDim'):
3 1 1 1 1 1 1
3 1 1 1 1 1 1 M:56 6 N:6 6 H:20 [24]
Degrees and weights  `d1 w11 w12 ... d2 w21 w22 ...'
  or `#lines #columns' (= `PolyDim #Points' or `#Points PolyDim'):
7 1 1 1 2 3 3 3
7 1 1 1 2 3 3 3 M:154 18 F:9
Degrees and weights  `d1 w11 w12 ... d2 w21 w22 ...'
  or `#lines #columns' (= `PolyDim #Points' or `#Points PolyDim'):
7 1 1 2 2 2 3 3
7 1 1 2 2 2 3 3 M:116 18 N:9 9 H:2 70 [-136]
Degrees and weights  `d1 w11 w12 ... d2 w21 w22 ...'
  or `#lines #columns' (= `PolyDim #Points' or `#Points PolyDim'):
3 2
Type the 6 coordinates as #pts=3 lines with dim=2 columns:
0 0
1 0
0 1
Warning: Input has index 3, should be 2!   M:3 3 F:3
Degrees and weights  `d1 w11 w12 ... d2 w21 w22 ...'
  or `#lines #columns' (= `PolyDim #Points' or `#Points PolyDim'):
1 1 1 0 0 0 0  2 0 0 1 1 1 1
1 1 1 0 0 0 0  2 0 0 1 1 1 1 M:20 8 N:6 6 H:[0]
\end{verbatim}
As the examples show, weight input {\tt d w1 ... wn} requires
$w_1+\ldots w_n = rd$; in other words the weights $q_i = w_i/d$ add up to $r$ 
rather than to 1 as in the standard case.
See \cite{Skarke:2012zg} for more information on how weight systems determine
Gorenstein cones.

{\tt nef.x -G} cannot be combined with all of the other options.
Nevertheless
{\tt -N} swaps the lattices $M$ and $N$ as usual;
{\tt -H}, {\tt -S}, {\tt -T} work as expected;
{\tt -t} works (without it no time information is given);
{\tt -c*} determines the index;
{\tt -R} displays the vertices of the input polytope if the cone is not reflexive of index $r$;
{\tt -V} displays the vertices of the support polytope of the dual ($N$ lattice) cone;
{\tt -g}, {\tt -d} display the full sets of points of the support
polytopes in the lattices $N$ or $M$, respectively (here no numbers can be specified with these options).
%%%%%%%%%%%%%%%%%%%%%%%%%%%%%%%%%%%%%%%%%%%%%%%%%%%%%%%%%%%%%%%%%%%%%%%%%%%%%%
\section{{\tt mori.x}}
\label{sec-mori}
The main purpose of {\tt mori.x} is
the computation of the Mori cone of toric varieties given by star 
triangulations of reflexive polytopes, which correspond to crepant 
subdivisions of the associated fans. 
The program is able to perform such triangulations
for four--dimensional polytopes with up to three non--vertex points
if the secondary fan is at most three--dimensional. The
program can also be used with a known triangulation as its input
starting from PALP release $2.1$.
This option, which was not contained in PALP 
$2.0$ as described in \cite{Braun:2011ik}, works for arbitrary dimensions.

\subsection{General aspects of {\tt mori.x}}
We distinguish two types of functionalities of {\tt mori.x}.
The first kind yields information about the appropriately resolved ambient space
(see options {\tt -g, -I, -m, -P, -K} below).
This includes the Stanley-Reisner (SR) ideal (with {\tt -g}) as well as
specific information on the geometry
of the lattice polytope that determines the ambient toric variety: incidence
structure of the facets ({\tt -I}), IP-simplices ({\tt -P}) and subdivisions 
of the fan ({\tt -g});
furthermore, the Oda-Park algorithm \cite{oda1991linear, Berglund:1995gd} 
is used to find the Mori cone of the ambient space ({\tt -m}).
The second kind of functionalities deals with the intersection ring 
({\tt -i, -t}) and topological quantities ({\tt -b, -c, -d})
of the embedded hypersurface. 
They are determined with the help of
SINGULAR \cite{DGPS}, a computer algebra system for polynomial computations. 
Correspondingly, 
the options {\tt -b, -i, -c, -t, -d} (as well as {\tt -a, -H}, see below)
need SINGULAR to be installed.

The generators of the Mori cone are given in terms of their intersections with
the toric divisors. For singular toric varieties, the Picard group of Cartier divisors is a non-trivial 
subgroup of the Chow group, which contains the Weil divisors. Hence one can consider the K\"ahler cone, 
which is dual to the Mori cone, as a cone in the vector space
spanned by the elements of either the Picard or the Chow group. The program
{\tt mori.x} only deals with simplicial toric varieties, 
for which the Picard group is always a finite index subgroup of the Chow group
\cite{Fulton:1993,oda1988convex}. 
Hence the Cartier divisors are integer multiples of the Weil divisors and this 
ambiguity does not arise.

Starting with PALP 2.1, {\tt mori.x} affords two distinct modes of operation.
If used with the option {\tt -M} arbitrary reflexive polytopes of any 
dimension can serve as input (at least in principle),
see section \ref{moriM} for details.
Without {\tt -M} the program only works if the input polytope can be
triangulated by {\tt mori.x} or if it does not require triangulation,
as we will outline in the following sections.

As described in \cite{Braun:2011ik}, {\tt mori.x} can perform
star triangulations of certain four--dimensional reflexive
polytopes. This operation was designed for the CY hypersurface case.
Generic CY hypersurfaces avoid point-like 
singularities of the ambient space as well as divisors that correspond to interior points of facets. Consequently, 
the algorithm performs star triangulations only up to such interior points.

Polytopes can be triangulated by subdividing the secondary fans of its
non-simplicial facets \cite{billera1990constructions, GKZ}. 
This triangulation algorithm is implemented in {\tt mori.x} for 
polytopes with up to three points that are neither vertices nor interior to 
the polytope or one of its facets; 
this implies that the secondary fan of any facet can be at most 
three--dimensional.
The program exits with a warning message if the subdivision is not properly 
completed.

As the dimension of the secondary fan corresponding to a facet grows with the number of 
points in the facet, this limitation tends to become relevant for toric varieties for 
which $h^{1,1}$ is large: $h^{1,1}$ increases with the number of points on the polytope and polytopes 
with many points are more likely to have facets containing many points. 

Complete triangulations of arbitrary polytopes can be performed
programs such as TOPCOM
\cite{Rambau:TOPCOM-ICMS:2002}, which is also included in the open source
mathematics software system Sage \cite{sage}. Sage also contains various tools for handling toric
varieties. 
The triangulations performed in {\tt mori.x} are 
attuned to the case of three-dimensional CY hypersurfaces.
This means, in particular, 
that interior points of facets are ignored: one must use {\tt -M} to
avoid this.
For small Picard numbers, {\tt mori.x} is hence faster than programs 
which perform a complete triangulation.

If a polytope of arbitrary dimension has only simplicial facets whose only
lattice points are its vertices and possibly interior points, it does not
require any triangulation.
Hence {\tt mori.x} can also handle such cases without {\tt -M}.

With the option {\tt -H} the program can also analyze arbitrary hypersurfaces 
embedded in the ambient toric varieties. 
It is capable of computing the intersection ring and certain characteristic 
classes. 
Here the omission of interior points of facets, which happens as a consequence
of {\tt mori.x}'s triangulation algorithm, may introduce severe singularities
which often result in non--integer intersection numbers.
There is a warning if there are indeed points interior to facets; in such a case
it is probably better to repeat the computation with the combination {\tt -HM}.

The help screen
provides essential information about all the functionalities of the program:

\begin{verbatim}
palp$ mori.x -h
This is ``mori.x'':
                 star triangulations of a polytope P* in N
                 Mori cone of the corresponding toric ambient spaces
                 intersection rings of embedded (CY) hypersurfaces
Usage:   mori.x [-<Option-string>] [in-file [out-file]]
Options (concatenate any number of them into <Option-string>):
 -h    print this information 
 -f    use as filter
 -g    general output: triangulation and Stanley-Reisner ideal
 -I    incidence information of the facets (ignoring IPs of facets)
 -m    Mori generators of the ambient space
 -P    IP-simplices among points of P* (ignoring IPs of facets)
 -K    points of P* in Kreuzer polynomial form
 -b    arithmetic genera and Euler number
 -i    intersection ring
 -c    Chern classes of the (CY) hypersurface
 -t    triple intersection numbers
 -d    topological information on toric divisors & 
       del Pezzo conditions
 -a    all of the above except h, f, I and K
 -D    lattice polytope points of P* as input (default CWS)
 -H    arbitrary (also non-CY) hypersurface 
       `H = c1*D1 + c2*D2 + ...' input: coefficients `c1 c2 ...'
 -M    Stanley-Reisner ideal and Mori generators with an
       arbitrary triangulation as input; must be combined with -D
Input: 1) standard: degrees and weights 
          `d1 w11 w12 ... d2 w21 w22 ...'
       2) alternative (use -D): `d np' or `np d' 
          (d=Dimension, np=#[points]) and (after newline) np*d 
          coordinates
Output:   as specified by options
\end{verbatim}

Following PALP's notation we refer to the $M$ lattice polytope which 
determines the CY hypersurface as $P$; consequently its dual, which 
gives rise to the fan of the ambient toric variety, is $P^*$. 

As PALP always interprets the input as $P\subset M_\mathbb{R}$ unless some 
option
modifies this behavior, matrix input of $P^*\subset N_\mathbb{R}$ requires 
the option {\tt -D}.
In order to avoid errors, matrix input is not allowed unless this option is set.
If only $P$ but not $P^*$ is known one can use {\tt poly.x -e} to obtain the
latter.

\subsection{Options of {\tt mori.x}}\label{sec:mori-options}
This section contains a detailed description of the options listed
in the help screen. If no flag is specified, the 
program starts with the parameter {\tt -g}. By default, the program considers a CY
hypersurface embedded in the ambient toric variety. The option {\tt -H} has to be used in order to consider
non-CY hypersurfaces. Note that the options {\tt -b, -i, -c, -t, -d, -a, -H} need SINGULAR \cite{DGPS} to be
installed.

Most options of {\tt mori.x} produce output that is related to the points
of $P^*$ in a specific order which can be determined by combining the
desired functionality with the option {\tt -P} (see sec. \ref{moriP} below).
In order to avoid repeating this information for every option,
we now present an example that will be used for many of the options below:
\begin{verbatim}
palp$ mori.x -P
Degrees and weights  `d1 w11 w12 ... d2 w21 w22 ...':
8 4 1 1 1 1 0  6 3 1 0 1 0 1
4 8  points of P* and IP-simplices
   -1    0    0    0    1    3    1    0
    0    0    0    1    0   -1    0    0
   -1    1    0    0    0    3    1    0
    1    0    1    0    0   -4   -1    0
------------------------------   #IP-simp=2
    4    1    0    1    1    1   8=d  codim=0
    3    0    1    1    0    1   6=d  codim=1
\end{verbatim}
The output above the dashed line just means that $P^*$ has the lattice
points\footnote{Here the index 
starts at $1$ instead of $0$ as it is standard in PALP. 
This shift is needed to match the counting of toric divisor classes displayed 
in certain outputs of {\tt mori.x} and hence avoids confusion.  
} 
\begin{equation}
p_1=\left( \begin{array}{r} -1 \\ 0 \\ -1 \\ 1 \end{array} \right), ~~~
p_2=\left( \begin{array}{r} 0 \\ 0 \\ 1 \\ 0 \end{array} \right), ~~~~~~
\cdots, ~~~
p_7=\left( \begin{array}{r} 1 \\ 0 \\ 1 \\ -1 \end{array} \right), ~~~
p_8=\left( \begin{array}{r} 0 \\ 0 \\ 0 \\ 0 \end{array} \right)
\label{standardex}\end{equation}
 and the last two lines encode the facts $4p_1+p_2+p_4+p_5+p_6=0$, 
$3p_1+p_3+p_4+p_6=0$.
Note how the dashed line proceeds only up to $p_6$. 
This is because {\tt mori.x} always ignores the origin, and, if used without 
{\tt -M}, also ignores points that are interior to facets:
$p_7=(p_2+p_4+p_5+p_6)/4=(p_3+p_4+p_6)/3$ lies inside
the facet with vertices $p_2,p_3,p_4,p_5,p_6$.
The reader is invited to check that the same example with {\tt mori.x -PM} 
results in a dashed line below all points except the origin.

\subsubsection{-h} 
This option prints the help screen.

\subsubsection{-f} 
This parameter suppresses the prompt of the command line. This is
useful if one wants to build pipelines or shorten the input;
e.g.~our standard example (\ref{standardex}) can be entered as
\begin{verbatim}
palp$  echo  '8 4 1 1 1 1 0  6 3 1 0 1 0 1' | mori.x -fP
4 8  points of P* and IP-simplices
...
\end{verbatim}

\subsubsection{-g}

This triggers the general output. 
First, the triangulation data of the facets is displayed. The number
of triangulated simplices is followed by the incidence structure of the 
simplices. The incidence information for each simplex is
encoded in terms of a bit sequence (cf. sec.~ \ref{inci}):
there is a digit for each relevant polytope point; a $1$ denotes that the point belongs to the simplex. 
Second, the SR ideal is displayed: the number of elements of the ideal is followed by its elements. Each element is denoted by a 
bit sequence as above. 

\begin{verbatim}
palp$ echo  '8 4 1 1 1 1 0  6 3 1 0 1 0 1' | mori.x -fg
8 Triangulation
110101 111100 101011 101110 100111 111001 001111 011101
2 SR-ideal
010010 101101
9 Triangulation
110101 111100 101011 101110 100111 111001 010111 011011 011110
2 SR-ideal
110010 001101
\end{verbatim}
\noindent The program performs the two possible triangulations of the
facet $\langle 23456 \rangle$, which is the only non--simplicial one 
(see section \ref{moriI}).
The last two bit sequences of the first result describe the simplices 
$\langle \widehat{25} 346\rangle$, 
whereas the second triangulation gives the three simplices 
$\langle 2\widehat{346}5\rangle$
(in this notation the hat indicates that one of the points is dropped).
Nevertheless, the two resolutions give the same CY intersection polynomial.\footnote{This 
fact suggests that the two resolutions give rise to the same CY hypersurface. 
Indeed, simply connected CY threefolds are completely determined up to 
diffeomorphisms by their Hodge numbers, intersection rings and second Chern 
classes \cite{wall1966classification,Batyrev:2008rp}.
}

\subsubsection{-I}\label{moriI}

The incidence structure of the facets of the polytope $P^*$ is
displayed. Interior points of the facets are neglected.
\begin{verbatim}
palp$  echo  '8 4 1 1 1 1 0  6 3 1 0 1 0 1' | mori.x -fI
Incidence: 110101 111100 011111 101011 101110 100111 111001
\end{verbatim}
\noindent 

The incidence data show the intersections 
of $p_1,\ldots,p_6$ (ignoring $p_7$, $p_8$!) 
with the seven facets. The third facet contains the
five points $p_2,\dots,p_6$, hence it is not simplicial and 
needs to be triangulated. See section~\ref{inci} for more details on the 
representation of incidences as bit sequences.

\subsubsection{-m}

The Mori cone generators of the ambient space are displayed in the
form of a matrix.\footnote{As divisors corresponding to interior 
points of facets do not intersect a CY hypersurface, such divisors are neglected in 
the computation of the Mori cone of the ambient space.} Each row corresponds to a generator. 
The entries of each row are the intersections of the generator with the toric
divisor classes. The Oda-Park algorithm is used to compute the generators.  
Furthermore, the incidence structure between the generators of the Mori cone and its facets
is displayed. 
For the standard example this takes the following form.
\begin{verbatim}
palp$ echo  '8 4 1 1 1 1 0  6 3 1 0 1 0 1' | mori.x -fm 
2 MORI GENERATORS / dim(cone)=2 
  3  0  1  1  0  1   I:10
  0  3 -4 -1  3 -1   I:01
2 MORI GENERATORS / dim(cone)=2 
  1  1 -1  0  1  0   I:10
  0 -3  4  1 -3  1   I:01
\end{verbatim}
\noindent 
The Mori cone is two-dimensional, so that its facets can be identified with the
generators. This explains the trivial incidence structure.

Let us consider another simple
example, a hypersurface in $\IP^2\times \IP^1\times \IP^1$.
\begin{verbatim}
palp$ mori.x -m
Degrees and weights  `d1 w11 w12 ... d2 w21 w22 ...':
3 1 1 1 0 0 0 0  2 0 0 0 1 1 0 0  2 0 0 0 0 0 1 1
3 MORI GENERATORS / dim(cone)=3 
   0  0  0  1  1  0  0   I:110
   1  1  0  0  0  0  1   I:101
   0  0  1  0  0  1  0   I:011
\end{verbatim}
\noindent The Mori cone generators can easily be seen to be dual to the hyperplane sections. Now, the Mori cone
is three-dimensional, so that each of its facets contains two generators. 
Let us, for instance, consider the incidence structure between the first generator and the three facets of the Mori cone.
Here, the string \verb+I:110+ tells that the vector lies on the first and second facets but does not intersect the third one.

For an example with a more complicated structure of the Mori cone see section
\ref{moriM}.

\subsubsection{-P}\label{moriP}
First a list of lattice points of $P^*$ is displayed in the following manner.
If {\tt -P} is combined with both {\tt -M} and {\tt -D}, the list is just
the input provided by the user, in the same order except for the fact that 
the lattice origin comes at the end of the list.
In all other cases the complete list of lattice points of $P^*$ is given in
the following order:\\[1mm]
1. vertices (with {\tt -D} in the order provided by the user),\\
2. points not interior to the polytope or its facets,\\
3. points interior to facets,\\
4. the lattice origin.\\[1mm]
Then a dashed line indicates which points are `relevant': all points
except for the origin in the case of {\tt -M}, but not points interior
to facets otherwise.
Finally the IP-simplices with vertices among these relevant points are 
displayed.

The output for the standard example can be found above equation 
(\ref{standardex}).
\del
As a further example consider a $\IP^1$ fibered over $\IP^3$:

\begin{verbatim}
palp$ mori.x -P
Degrees and weights  `d1 w11 w12 ... d2 w21 w22 ...':
5 1 1 1 1 1 0  2 0 0 0 0 1 1
4 7  points of P* and IP-simplices
    1    0    0    0   -1    0    0
    0    1    0    0   -1    0    0
    0    0    1    0   -1    0    0
    0    0    0    1    1   -1    0
------------------------------   #IP-simp=2
    1    1    1    0    1    1   5=d  codim=0
    0    0    0    1    0    1   2=d  codim=3
\end{verbatim}
\enddel
The following example features all types of lattice points:
\begin{verbatim}
palp$ echo '16 8 4 2 1 1' | mori.x -fP
4 9  points of P* and IP-simplices
   -1    0    0    2    0    0    0    1    0
   -1    0    0    1    2    0    1    1    0
    0    0    2   -1    1    1    1    0    0
    0    1    1    0   -1    1    0    0    0
-----------------------------------   #IP-simp=3
    8    1    1    4    2    0    0  16=d  codim=0
    4    0    0    2    1    1    0   8=d  codim=1
    2    0    0    1    0    0    1   4=d  codim=2
\end{verbatim}
$p_1,\ldots,p_5$ are vertices, $p_6, p_7$ further relevant points, but
$p_8=-p_1=(p_2+p_3+4p_4+2p_5)/8$ is interior to the facet spanned by
$p_2,\ldots,p_5$.
Note that the ordering of the CWS input is not obeyed by the output of
lattice points. Once the order is displayed, however, it
is fixed and determines the labeling of toric divisors in any further output.

\subsubsection{-K}
The Kreuzer polynomial\footnote{We named this output format after Maximilian
Kreuzer, who designed it. This is an example of his proverbial ability to
eliminate unnecessary data redundancies and recast essential information in 
condensed form.} of PALP's representation of $P^*$ is displayed. It encodes lattice polytope points in a compact form. The
number of variables equals the dimension of the polytope. Each lattice point gives rise to a Laurent monomial
in which the exponents of the variables are the coordinates. Vertices and non-vertices are distinguished by
coefficients `$+$' and `$-$' respectively. 
Points in the interior of facets are ignored.
As this is closely connected with the way {\tt mori.x} works when used without
{\tt -M}, the combination {\tt -MK} is not allowed.
\begin{verbatim}
palp$ echo  '8 4 1 1 1 1 0  6 3 1 0 1 0 1' | mori.x -fK
KreuzerPoly=t_4/(t_1t_3)+t_3+t_4+t_2+t_1+t_1^3t_3^3/(t_2t_4^4); 
intpts=1;  Pic=2
\end{verbatim}
\noindent 

A comparison with the output for {\tt -P}, which can be found above equation
(\ref{standardex}),
might help for a better understanding of the present option.

Negative coordinates are always displayed by putting the
variables in the denominator.
The number of points in the interior of facets is shown
as \verb+intpts+. The multiplicities of the toric divisors are displayed as
\verb+multd+ if they are greater than one. Furthermore, the Picard 
number of the CY hypersurface is computed and printed as \verb+Pic+. 

\subsubsection{-b}
The zeroth and first arithmetic genera of the hypersurface are determined 
according to the following formulas \cite{hirzebruch1995topological}:
\begin{equation}
\chi_q(X) = \sum_p (-1)^ph^{p,q}(X) = \int_X \mathrm{ch}(\Omega^q(X))
\mathrm{Td}(X),\quad q = 0,1.
\end{equation}
where $\Omega^0(X) = \cO_X$ is the trivial bundle and $\Omega^1(X) =
T^*X$ is the bundle of 1-forms, $\mathrm{ch}$ is the corresponding Chern character,
and $\mathrm{Td}(X)$ is the Todd class of $X$. 

Furthermore the Euler characteristic is displayed.
Here we compute it by means of the intersection polynomial:
\begin{equation}
\chi =\int_X \mathrm{c}_n\,.
\end{equation}
where $\mathrm{c}_n$ is the top Chern class, $n= \dim X$. 

These formulas hold for arbitrary smooth hypersurfaces; in particular, they do not need 
to be CY.
Indeed, if $X$ is CY, its Euler characteristic can also be computed by {\tt poly.x} 
in terms of polytope combinatorics.
Compare the two Euler characteristics for a consistency check.

Consider the $K3$ surface as a simple example:
\begin{verbatim}
palp$ echo '4 1 1 1 1' | mori.x -bf
SINGULAR  -> Arithmetic genera and Euler number of the CY:
chi_0:  2 , chi_1: -20  [ 24 ]
\end{verbatim}
As expected, the Euler characteristic is $24$ and $h^{1,1}=20$. 
Using the example discussed before we find

\begin{verbatim}
palp$ echo '8 4 1 1 1 1 0 6 3 1 0 1 0 1' | mori.x -bf
SINGULAR  -> Arithmetic genera and Euler number of the CY:
chi_0: 0 , chi_1: 126  [ -252 ]
SINGULAR  -> Arithmetic genera and Euler number of the CY:
chi_0: 0 , chi_1: 126  [ -252 ]
\end{verbatim}
Special care is needed in the interpretation of the results for non-CY 
hypersurfaces: 
the triangulation algorithm might fail to make these varieties smooth, in 
which case the formulas above do not hold and hence the output is misleading; 
see the description of the option {\tt -H} for more details.

\subsubsection{-i}
This option displays the intersection polynomial 
in terms of an integral basis of the toric divisors. The
coefficients of the monomials are the triple intersection numbers in this basis. This option can also 
be used together with {\tt -H} to perform this task for non-CY hypersurfaces. 
\begin{verbatim}
palp$ echo  '8 4 1 1 1 1 0  6 3 1 0 1 0 1' | mori.x -fi
SINGULAR -> divisor classes (integral basis J1 ... J2):
d1=J1+3*J2, d2=J1, d3=-J1+J2, d4=J2, d5=J1, d6=J2
SINGULAR -> intersection polynomial:
2*J1*J2^2+2*J2^3
SINGULAR -> divisor classes (integral basis J1 ... J2):
d1=J1+3*J2, d2=J1, d3=-J1+J2, d4=J2, d5=J1, d6=J2
SINGULAR -> intersection polynomial:
2*J1*J2^2+2*J2^3
\end{verbatim}
\noindent \verb+d1+,$\,\dots$, \verb+d6+ denote the toric divisors corresponding
to the lattice points $p_1,\dots,p_6$, cf. eq. (\ref{standardex}). 
There are two independent divisor 
classes. Indeed, {\tt mori.x} expresses the intersection polynomial in terms
of the integral basis $J_1=D_2=D_5$ and $J_2=D_4=D_6$.

\subsubsection{-c}
The Chern classes of the hypersurface (CY or non-CY) are displayed in
terms of an integral basis of the toric divisors:
\begin{verbatim}
palp$ echo '8 4 1 1 1 1 0 6 3 1 0 1 0 1' | mori.x -fc
SINGULAR -> divisor classes (integral basis J1 ... J2):
d1=J1+3*J2, d2=J1, d3=-J1+J2, d4=J2, d5=J1, d6=J2
SINGULAR  -> Chern classes of the CY-hypersurface:
c1(CY)=  0
c2(CY)=  10*J1*J2+12*J2^2
c3(CY)=  -252 *[pt]
SINGULAR -> divisor classes (integral basis J1 ... J2):
d1=J1+3*J2, d2=J1, d3=-J1+J2, d4=J2, d5=J1, d6=J2
SINGULAR  -> Chern classes of the CY-hypersurface:
c1(CY)=  0
c2(CY)=  10*J1*J2+12*J2^2
c3(CY)=  -252 *[pt]
\end{verbatim}

\subsubsection{-t}

The triple intersection numbers of the toric divisors are displayed.
The form of this output\footnote{The pre-compiler command {\tt DijkEQ} in the C file {\tt SingularInput.c} 
controls the symbol `{\tt ->}' in option {\tt -t}.} is 
designed for further use in Mathematica \cite{mathematica}.
\del
As a (mildly non-trivial) example, consider an elliptic threefold over 
$\mathbb{CP}^2$ which has a single section.  
In order to have control over the ordering of the toric divisors
\new
??? - siehe Problem mit Ordnung der Punkte ???
 \endnew
, we first 
execute
\begin{verbatim}
palp$ mori.x -P
Degrees and weights  `d1 w11 w12 ... d2 w21 w22 ...':
18 1 1 1 6 9 0 6 0 0 0 2 3 1
4 10  points of P* and IP-simplices
   -1    0    1    0    0    0    0    0    0    0
   -1    0    0    1    0    0    0    0    0    0
    0    2    0    0   -3    0   -2    1   -1    0
    3    1    0    0   -2    1   -1    1    0    0
------------------------------   #IP-simp=2
    1    9    1    1    6    0  18=d  codim=0
    0    3    0    0    2    1   6=d  codim=2

\end{verbatim}
Before computing the intersection ring, let us state our expectations. 
As the section corresponds to the divisor $D_6$, any intersection between two divisors of the base (the intersection ring of which is
generated by $D_1$, $D_3$ and $D_4$) and $D_6$ must be unity
Furthermore, the exceptional set of the fibre translates
to $D_6$ having a vanishing intersection with $D_2$ and $D_5$ on the Calabi-Yau threefold hypersurface.
This is confirmed by running
\begin{verbatim}
palp$ mori.x -t
Degrees and weights  `d1 w11 w12 ... d2 w21 w22 ...':
18 1 1 1 6 9 0 6 0 0 0 2 3 1
SINGULAR -> triple intersection numbers:
[...]
d5^2*d6->0,
d4*d5*d6->0,
d3*d5*d6->0,
d2*d5*d6->0,
d1*d5*d6->0,
d4^2*d6->1,
d3*d4*d6->1,
d2*d4*d6->0,
d1*d4*d6->1,
d3^2*d6->1,
d2*d3*d6->0,
d1*d3*d6->1,
d2^2*d6->0,
d1*d2*d6->0,
d1^2*d6->1,
[...]
\end{verbatim}
As the output is rather long in this case, we have just displayed the intersections linear in $D_6$.
\enddel
Before computing the intersection ring for our standard example 
(\ref{standardex}), let us state some expectations. 
Inspection of the data of the polytope reveals that it describes a K3 
fibration with the fiber determined by the weight system 6~3~1~1~1.
There are only the two points $p_2, p_5$ outside the corresponding 3--plane,
so each of them must represent the generic fiber with self--intersection 0.
In other words, the self--intersections of $d_2$ and $d_5$ as well as 
$d_2\cdot d_5$ must all vanish.
This is confirmed by the following excerpt from the output:
\begin{verbatim}
echo '8 4 1 1 1 1 0  6 3 1 0 1 0 1' | mori.x -ft
SINGULAR -> triple intersection numbers:
d6^3->2,
d5*d6^2->2,
d4*d6^2->2,
d3*d6^2->0,
d2*d6^2->2,
d1*d6^2->8,
d5^2*d6->0,
d4*d5*d6->2,
d3*d5*d6->2,
d2*d5*d6->0,
d1*d5*d6->6,
d4^2*d6->2,
d3*d4*d6->0,
d2*d4*d6->2,
d1*d4*d6->8,
d3^2*d6->-2,
d2*d3*d6->2,
d1*d3*d6->2,
d2^2*d6->0,
d1*d2*d6->6,
d1^2*d6->30,
d5^3->0,
d4*d5^2->0,
d3*d5^2->0,
d2*d5^2->0,
d1*d5^2->0,
d4^2*d5->2,
d3*d4*d5->2,
d2*d4*d5->0,
d1*d4*d5->6,
d3^2*d5->2,
d2*d3*d5->0,
d1*d3*d5->6,
d2^2*d5->0,
d1*d2*d5->0,
d1^2*d5->18,
[...]
\end{verbatim}

\subsubsection{-d}
This option displays topological data of the toric divisors restricted
to the (CY or non-CY) hypersurface. The Euler characteristics of the toric 
divisor classes and their arithmetic genera are shown. 

Furthermore, in the case of a three-dimensional hypersurface, the program checks the del Pezzo property 
against two necessary conditions and analyses the mutual intersections of the del Pezzo candididates.
The number of del Pezzo candidates is displayed followed by their type in parenthesis; 
furthermore, those among them that do not intersect other del Pezzo candidates are listed. 

For a del Pezzo divisor $S$ of type $n$, the following equations should
hold:
\begin{equation}
\int_S c_1 (S)^2 = 9-n \, , \qquad \int_S c_2(S) = n+3 \quad \Longrightarrow
\quad \chi_0(S)=\int_S \text{Td}(S)=1\, .
\end{equation} 
Here, \text{Td(S)} denotes the Todd class of $S$, which gives the zeroth
arithmetic genus of $S$ upon integration. 
This test also allows to determine the type of the del Pezzo surface in
question. 
A second necessary condition comes from the fact that a del Pezzo surface is a
two-dimensional Fano manifold. 
Hence, the first Chern class of $S$ integrated over all curves on $S$ has to be
positive:
\begin{equation}
 D_i \cap S\cap c_1(S) >0 \qquad\forall D_i:\, D_i\neq S \, ,\quad  \; D_i\cap
S\neq 0 \, .
\end{equation}  
This condition would be sufficient if we were able to access \textit{all} curves
of the hypersurface.
In our construction, however, we can only check for curves induced by toric
divisors. This functionality was added to carry out the analysis of base manifolds 
for elliptic fibrations in \cite{Knapp:2011wk}.

Consider the following example: it is well-known that the del Pezzo surface $dP_6$ can be
realized as a homogeneous polynomial of degree $3$ in $\mathbb{CP}^3$. Hence a Calabi-Yau
hypersurface in a toric variety with CWS
\begin{align}
5\,1\, 1\, 1\, 1\, 1\, 0 \nonumber\\
2\, 0\, 0\, 0\, 0\, 1\, 1 \nonumber
\end{align}
i.e. a $\mathbb{CP}^1$ fibration over $\mathbb{CP}^3$ contains a $dP_6$: setting
the last coordinate $z_6$ to zero forces all terms to be of the form $z_5^2 P_3(z_1,\ldots ,z_4)$, 
where $P_3(z_1,\ldots ,z_4)$ is a homogeneous polynomial of degree $3$ in $z_1,\ldots ,z_4$. We may set
$z_5$ to $1$ by using the second $\mathbb{C}^*$ action, so that the divisor $D_6$ corresponds to 
a homogenous polynomial of degree $3$ in $\mathbb{CP}^3$, i.e. a $dP_6$.

This is confirmed by 
\begin{verbatim}
palp$ mori.x -d
Degrees and weights  `d1 w11 w12 ... d2 w21 w22 ...':
5 1 1 1 1 1 0 2 0 0 0 0 1 1
SINGULAR -> topological quantities of the toric divisors:
Euler characteristics: 46 46 46 9 46 55 
Arithmetic genera: 4 4 4 1 4 5 
dPs: 1 ; d4(6)  nonint: 1 ; d4
\end{verbatim}
Note that \textsc{Palp} has exchanged the ordering of the divisors, so that 
the $dP_6$ is now given by $D_4$. 
This divisor does not intersect any other del Pezzo as it is the 
only del Pezzo candidate in this example.

\subsubsection{-a}
This is a shortcut for {\tt -gmPbictd}.

\subsubsection{-D}
An alternative way to provide the input is to type lattice polytope points
of $P^*$ directly. In this case, one has to use the parameter {\tt -D}. 
Let us reconsider the example of sec.~\ref{moriP}: 
\begin{verbatim}
palp$ mori.x -DP
`#lines #columns' (= `PolyDim #Points' or `#Points PolyDim'):
4 5
Type the 20 coordinates as dim=4 lines with #pts=5 columns:
-1  2  0  0  0
-1  1  2  0  0
 0 -1  1  0  2
 0  0 -1  1  1
4 9  points of P* and IP-simplices
   -1    2    0    0    0    0    0    1    0
   -1    1    2    0    0    0    1    1    0
    0   -1    1    0    2    1    1    0    0
    0    0   -1    1    1    1    0    0    0
-----------------------------------   #IP-simp=3
    8    4    2    1    1    0    0  16=d  codim=0
    4    2    1    0    0    1    0   8=d  codim=1
    2    1    0    0    0    0    1   4=d  codim=2
\end{verbatim}
Note how the order of the vertices corresponds to that of the input
(cf.~section \ref{moriP}).

\subsubsection{-H}\label{moriH}
Using this option, one can specify a (non-CY) hypersurface. The user
determines the hypersurface divisor class $H = \sum_i c_i D_i$ in terms of the
toric divisor classes $D_i$ by typing its coefficients $c_i$. The hypersurface 
can then be analyzed by combining {\tt -H} with other
options, as described above. Just using {\tt -H}, the program runs {\tt -Hb}.

The reader is warned:
smoothness is not guaranteed anymore, so that
the intersection numbers can become fractional. Some choices of the hypersurface
equation may intersect singularities not resolved by the triangulation. 
Consider e.g. the hypersurface determined by the divisor class 
$H =  D_1 + D_6$ in our example (\ref{standardex}).
Remember that the
order in which {\tt mori.x} expects the coefficients of the hypersurface 
divisor class is fixed by the polytope matrix and not by the CWS input. Hence,
the correct input for $H$ is the string \verb+1 0 0 0 0 1+.  
\begin{verbatim}
palp$ mori.x -H
Degrees and weights  `d1 w11 w12 ... d2 w21 w22 ...'
8 4 1 1 1 1 0  6 3 1 0 1 0 1
WARNING: there is 1 facet-IP ignored in the triangulation. 
This may lead to unresolved singularities in the hypersurface. 
Type the 6 (integer) entries for the hypersurface class:
1 0 0 0 0 1
Hypersurface degrees: ( 5  4 )
Hypersurface class: 1*d1 1*d6 
SINGULAR  -> Arithmetic genera and Euler number of H:
chi_0: 29/27 , chi_1: 128/27  [ -22/3 ]
\end{verbatim}
\noindent To calculate these quantities, the program determines the
characteristic classes of the divisors using adjunction. 
It then performs the appropriate integration with the help of the triple
intersection numbers. The fractional results of the arithmetic 
genera and the Euler number in our example indicate that the intersection
polynomial has fractional entries. This happens because the 
program introduces a singularity into the ambient toric variety which
descends to the hypersurface $H$. 
\del
In the first triangulation, the simplex $\langle 1235\rangle$ still 
has volume four. In the second triangulation, there are two simplices with
volume three: $\langle 1 2 3 4\rangle$ and $\langle 1 3 4 5\rangle$. 
Indeed, any hypersurface in the divisor class $H =  D_1 + D_6$ is forced to pass
through the corresponding singularities.
\enddel
It is therefore much better to combine {\tt -H} with {\tt -M}:
\begin{verbatim}
palp$ mori.x -HM
Degrees and weights  `d1 w11 w12 ... d2 w21 w22 ...':
8 4 1 1 1 1 0  6 3 1 0 1 0 1
4 8  
   -1    0    0    0    1    3    1    0
    0    0    0    1    0   -1    0    0
   -1    1    0    0    0    3    1    0
    1    0    1    0    0   -4   -1    0
`#triangulations': 
1
1 triangulations:
12 1111000 1110010 1101010  0111001 0110011 0101011 
1011100 1010110 1001110  0011101 0010111 0001111 
Type the 7 (integer) entries for the hypersurface class:
1 0 0 0 0 1 0
Hypersurface degrees: ( 5  4  1 )
Hypersurface class: 1*d1 1*d6 
SINGULAR  -> Arithmetic genera and Euler number of H:
chi_0:  1 , chi_1: 3  [ -4 ]
\end{verbatim}
As a second example, consider the quadric in $\mathbb{CP}^3$:
\begin{verbatim}
palp$ mori.x -H
Degrees and weights  `d1 w11 w12 ... d2 w21 w22 ...':
4 1 1 1 1
Type the 4 (integer) entries for the hypersurface class:
2 0 0 0 
Hypersurface degrees: ( 2 )
Hypersurface class: 2*d1 
SINGULAR  -> Arithmetic genera and Euler number of H:
chi_0:  1 , chi_1: -2  [ 4 ]
\end{verbatim}
The hypersurface is smooth in this case, so that the arithmetic genera and 
the Euler number are those of $\mathbb{CP}^1\times\mathbb{CP}^1$. 
Of course, one needs to independently check smoothness in order to rely on 
the output, as the integrality of the arithmetic genera alone is not 
sufficient to conclude that the hypersurface is non-singular.

\subsubsection{-M}\label{moriM}
Option {\tt -M}  allows polytopes of (in principle) arbitrary dimensions,
but expects the triangulations to be provided by the user;
it can be combined with any other option except for {\tt -K}.
As the Mori cone is analysed with PALP's routines, the parameter 
{\tt POLY\_Dmax} must possibly be adjusted for this; see sec.~\ref{polydmax}.

This functionality is useful, for instance, when {\tt mori.x} fails to 
triangulate the polytope by itself.
This happens whenever the dimension of the polytope is different 
from four, or when the polytope contains more than three lattice points that 
are neither vertices nor (facet--)IPs.
Fortunately, there are programs capable of efficiently performing complete 
triangulations of arbitrary polytopes \cite{Rambau:TOPCOM-ICMS:2002, sage}; 
the user can redirect their output as an input for {\tt mori.x} to determine 
the Mori generators of the ambient space. 
Other situations where {\tt -M} is useful arise whenever we prefer to keep 
control over the lattice points involved in the triangulation, rather than 
accept {\tt mori.x}'s convention of omitting precisely the interior points of 
facets from a completed list; 
this is particularly relevant if {\tt -M} is combined with {\tt -H}. 

After the usual polytope input, $P^*$ is displayed in the following manner.
If the input is of CWS type, all points of $P^*$ are given in {\tt mori.x}'s
standard order (see section \ref{moriP}).
If matrix input is used via {\tt -D}, the points entered by the user are
displayed again in the same order, but with the origin appended (if the 
origin is accidentally entered somewhere in the point list, it is swapped
with the last point in the list); possible further polytope points are
ignored by the program, hence singularities can be introduced if desired.
Then the user is asked for the number of triangulations to be analysed, and
afterwards each triangulation should be entered as a line starting with the 
number of simplices involved in the triangulation, followed by bit sequences 
encoding these simplices. 
The number of bits in each sequence should be the number of non--zero lattice
points in the displayed list, with 1's indicating that the point belongs to 
the simplex, and 0's otherwise.

An application to our standard example (\ref{standardex}) was already
demonstrated in section \ref{moriH}.
Consider also the following two-dimensional polytope:
\begin{verbatim}
palp$ mori.x -MDgm
`#lines #columns' (= `PolyDim #Points' or `#Points PolyDim'):
2 7
Type the 14 coordinates as dim=2 lines with #pts=7 columns:
    1    0   -1   -1   -1   -1    0
    0    1    2    1    0   -1   -1
2 8  
    1    0   -1   -1   -1   -1    0    0
    0    1    2    1    0   -1   -1    0
`#triangulations': 
1
1 triangulations:
7 1100000 0110000 0011000 0001100 0000110 0000011 1000001
14 SR-ideal
0000101 0001001 0001010 0010001 0010010 0010100 0100001 0100010 
0100100 0101000 1000010 1000100 1001000 1010000
6 MORI GENERATORS / dim(cone)=5
  1 -2  1  0  0  0  0   I:0111011
  0  1 -1  1  0  0  0   I:1101101
  0  0  1 -2  1  0  0   I:1110110
  0  0  0  1 -2  1  0   I:1011111
  0  0  0  0  1 -1  1   I:1100011
  1  0  0  0  0  1 -1   I:0111100
\end{verbatim}
\noindent 
Since $P^*$ is just a polygon there is only one maximal triangulation of its 
boundary.
For this reason we have typed \verb+1+ after \verb+`#triangulations':+.
The triangulation has seven simplices, which are just the line segments
along the circumference of the polygon. 
For instance, the string \verb+0011000+ denotes the third simplex containing
the points denoted by the third and fourth columns of the
polytope matrix, i.e. the points $p_3=(-1,2)$ and $p_4=(-1,1)$.   

As the Mori cone encodes linear relations among seven points in $d=2$, it is 
five-dimensional and has four-dimensional facets; there are seven of 
them, as the lengths of the bit sequences after {\tt I:} indicate.
There are six generators.
Consider the matrix of incidences whose rows are preceded by \verb+I:+.
The second column reads \verb+111011+, i.e. on the second facet of the Mori 
cone lie five generators. 
It is easily checked that they satisfy $m_1+2m_2+m_3=m_5+m_6$.
All other facets contain instead four generators and are hence simplicial.

%%%%%%%%%%%%%%%%%%%%%%%%%%%%%%%%%%%%%%%%%%%%%%%%%%%%%%%%%%%%%%%%%%%%%%%%%%%%%%
\section*{Acknowledgments} This work is based on the legacy of Maximilian Kreuzer who has been an inspiration for all of us. 
We are grateful to Benjamin Nill for providing information on several PALP options and for useful remarks. J. Knapp 
thanks the Vienna University of Technology for hospitality. N.-O. Walliser thanks the \'Ecole Polytechnique, Paris for hospitality. The work of A.P.Braun was supported by the FWF under grant I192. The work of J. Knapp was supported by World Premier International Research Center Initiative (WPI Initiative), MEXT, Japan. 

%\bibliographystyle{utphys}
%\bibliography{PALP-Bibliography_arxiv}

\def\cprime{$'$}
\providecommand{\href}[2]{#2}\begingroup\raggedright\endgroup

\end{document}